\documentclass{amsart}
\usepackage{amscd,amssymb,amsxtra}
\usepackage[mathscr]{eucal}
\usepackage[all,2cell]{xy}\UseAllTwocells

\CompileMatrices

\title{Loop groups and twisted $K$-theory I}

\author{Daniel~S.~Freed}
\address{Department of Mathematics \\ University of Texas
\\Austin, TX}
\email{dafr@math.utexas.edu}

\author{Michael~J.~Hopkins}
\address{Department of Mathematics \\ Massachusetts Institute of
Technology\\Cambridge, MA 02139-4307}
\email{mjh@math.mit.edu}

\author{Constantin~Teleman}
\address{Cambridge\\Cambridge}
\email{c.teleman}

\newtheorem{thm}[equation]{Theorem}
\newtheorem{cor}[equation]{Corollary}
\newtheorem{lem}[equation]{Lemma}
\newtheorem{prop}[equation]{Proposition}

\newtheorem{assumption}[equation]{Assumption}

\ifx\undefined\theoremstyle
	\newtheorem{definition}[equation]{Definition}
	\newenvironment{defin}{\begin{definition}\rm}{\end{definition}}
	\newtheorem{conj}[equation]{Conjecture}
	\newtheorem{rem}[equation]{Remark}
	\newtheorem{rems}[equation]{Remarks}

\else
	\theoremstyle{definition}

	\newtheorem{defin}[equation]{Definition}

	\theoremstyle{remark}

	\newtheorem{rem}[equation]{Remark}

\ifx\undefined\prem	\newtheorem{prem}{\bf Working Remark} \fi

\fi

\newtheorem{eg}[equation]{Example}

\ifx\undefined\pf
\newenvironment{pf}{\bigskip{\em Proof:\/}}{\qed\medskip}
\newenvironment{pf*}[1]{\bigskip{\em #1:\/}}{\qed\medskip}
\fi

\ifx\undefined\numberwithin
\def\numberwithin#1#2{\makeatletter\@ifundefined{c@#1}{\@nocnterrr}{%
  \@ifundefined{c@#2}{\@nocnterr}{%
  \@addtoreset{#1}{#2}%
  \toks@\expandafter\expandafter\expandafter{\csname the#1\endcsname}%
  \expandafter\xdef\csname the#1\endcsname
    {\expandafter\noexpand\csname the#2\endcsname
     .\the\toks@}}}\makeatother}\fi

\ifx\undefined\qed\newcommand{\qed}{\hfil\rule{4pt}{6pt}\bigskip}\fi
\ifx\undefined\\operatorname\newcommand{\operatorname}[1]{\mathop{\mbox{\rm #1}}}\fi
\newcommand{\Hom}{\operatorname{Hom}}

\newcommand{\Map}{\operatorname{Map}}
\newcommand{\Aut}{\operatorname{Aut}}

\newcommand{\tor}{\operatorname{Tor}}

\newcommand{\Q}{{\mathbb Q}}
\newcommand{\Z}{{\mathbb Z}}
\newcommand{\R}{{\mathbb R}}

\newcommand{\zerowidth}[1]{\hbox to 0pt{\hss$\displaystyle #1$\hss}}

\newcommand{\C}{\mathbb C}

\ifx\undefined\eqref\newcommand{\eqref}[1]{\rm (\ref{#1})}\fi

\newcounter{thmItem}

\newcommand{\thmItemref}[1]
 	 {{\rm \ref{#1})}}

\newenvironment{thmList}{\begin{list}%
{\rm \roman{thmItem})}{\usecounter{thmItem}
\setlength{\labelwidth}{2em}
\setlength{\itemindent}{2em}
\setlength{\leftmargin}{0pt}
\setlength{\listparindent}{0pt}
\setlength{\parsep}{0pt}
\setlength{\partopsep}{0pt}
\setlength{\itemsep}{\medskipamount}
\setlength{\topsep}{\medskipamount}
}}{\end{list}}

\newcounter{textItem}

\newenvironment{textList}{\begin{list}%
{\rm (\arabic{textItem})}{\usecounter{textItem}
\setlength{\labelwidth}{2em}
\setlength{\itemindent}{2em}
\setlength{\leftmargin}{0pt}
\setlength{\listparindent}{0pt}
\setlength{\parsep}{0pt}
\setlength{\partopsep}{0pt}
\setlength{\itemsep}{\medskipamount}
\setlength{\topsep}{\medskipamount}
}}{\end{list}}

\newcounter{condItem}

\ifx\undefined\slot
\newcommand{\slot}{\,-\,}
\fi

\newcounter{probi}

\newcounter{problemi}

\newcounter{multListi}

\newcounter{nmultListi}

\newcommand{\gl}[1]{\mathfrak{gl}_{1}}

\DeclareMathOperator{\rep}{R}

\newcommand{\prek}{\check{\mathcal K}}
\newcommand{\sk}{\mathcal K}

\newcommand{\st}{\mathcal F}
\DeclareMathOperator{\ass}{sh}
\DeclareMathOperator{\waff}{\mathit{W_{aff}}}
\DeclareMathOperator{\wa}{\mathit{W_{aff}^{e}}}
\DeclareMathOperator{\na}{\mathit{N_{aff}^{e}}}

\newcommand{\buone}{\mathcal{B}\T}
\newcommand{\buonepm}{\mathcal{B}\tpm}
\newcommand{\bzt}{\mathcal{B}\Z/2}
\newcommand{\tpm}{\T^{\pm}}

\newcommand{\modmod}[1]{/\!/{#1}}
\newcommand{\rel}{\text{rel}}
\newcommand{\id}{\text{Id}}
\newcommand{\T}{\mathbb T}
\newcommand{\LambdaTorsor}{\Lambda^{\tau}}

\DeclareMathOperator{\fredholm}{Fred}
\newcommand{\fred}[1]{\fredholm^{(#1)}}

\DeclareMathOperator{\covering}{Cov}
\newcommand{\cov}[1]{\underline{\covering}_{#1}}
\newcommand{\ocov}[1]{\covering_{#1}}
\DeclareMathOperator{\frakext}{\mathfrak{Ext}}
\newcommand{\extonecat}{\frakext}
\newcommand{\exttwocat}{\underline{\frakext}}
\newcommand{\hext}[1]{\extonecat_{#1}}
\newcommand{\sext}[1]{\exttwocat_{#1}}
\DeclareMathOperator{\shom}{\mathfrak{Hom}}
\DeclareMathOperator{\twistcat}{\mathfrak{Twist}}
\newcommand{\twist}[1]{\mathfrak{Twist}_{#1}}
\newcommand{\cms}[1]{\left[#1\right]}
\DeclareMathOperator{\spin}{Spin}
\DeclareMathOperator{\spinc}{Spin^{c}}
\newcommand{\point}{\text{pt}}
\DeclareMathOperator{\ad}{ad}
\newcommand{\ul}{\underline}
\DeclareMathOperator{\cyl}{cyl}
\DeclareMathOperator{\cl}{C\ell}
\DeclareMathOperator{\pinc}{Pin^{c}}

\DeclareMathOperator{\pfaff}{\mathit{pfaff}}
\newcommand{\dualc}{\Check h}
\newcommand{\efix}{\lambda_{0}}

\begin{document}

\maketitle

\tableofcontents

\section*{Introduction}
\label{sec:introduction}

Equivariant $K$-theory focuses a remarkable range of perspectives on
the study of compact Lie groups.  One finds tools from topology,
analysis, and representation theory brought together in describing the
equivariant $K$-groups of spaces and the maps between them.  In the
process all three points of view are illuminated.  Our aim in this
series of
papers~\cite{freed:_loop_group_twist_k_theor_ii,freed:_loop_group_twist_k_iii}
is to begin the development of similar relationships when a compact
Lie group $G$ is replaced by $LG$, the infinite dimensional group of
smooth maps from the circle to $G$.

There are several features special to the representation theory of
{\em loop} groups.  First of all, we will focus only on the
representations of $LG$ which have ``positive energy.''  This means
that the representation space $V$ admits an action of the rotation
group of the circle which is (projectively) compatible with the action
of $LG$, and for which there are no vectors $v$ on which rotation by
$\theta$ acts by multiplication by $e^{i n\theta}$ with $n<0$.  It
turns out that most positive energy representations are projective,
and so $V$ must be regarded as a representation of a central extension
$LG^{\tau}$ of $LG$ by $U(1)$.  The topological class of this central
extension is known as the {\em level}.  One thing a topological
companion to the representation theory of loop groups must take into
account is the level.

Next there is the fusion product.  Write $\rep^{\tau}(LG)$ for the
group completion of the monoid of positive energy representations of
$LG$ at level $\tau$.  In~\cite{verlinde88:_fusion_d}, Erik Verlinde
introduced a multiplication on $\rep^{\tau}(LG)\otimes\C$ making it
into a commutative ring (in fact a Frobenius algebra).  This
multiplication is called {\em fusion}, and $\rep^{\tau}(LG)\otimes\C$,
equipped with the fusion product is known as the {\em Verlinde
algebra}.  The fusion product also makes $\rep^{\tau}(LG)$ into a ring
which we will call the {\em Verlinde ring}.  A good topological
description of $\rep^{\tau}(LG)$ should account for the fusion product
in a natural way.

The positive energy representations of loop groups turn out to be
completely reducible, and somewhat surprisingly, there are only
finitely many irreducible positive energy representations at a fixed
level.  Moreover, an irreducible positive energy representation is
determined by its lowest non-trivial energy eigenspace, $V(n_{0})$,
which is an irreducible (projective) representation of $G$.  Thus the
positive energy representations of $LG$ correspond to a {\em subset}
of the representations of the compact group $G$.  This suggests that
$G$-equivariant $K$-theory might somehow play a role in describing the
representations of $LG$.  In fact this is the case.  Here is our main
theorem.
\begin{thm}
\label{thm:4} Let $G$~be a connected compact Lie group and $\tau $~a
level for the loop group.  The Grothendieck group $\rep^{\tau}(LG)$ at
level~$\tau $ is isomorphic to a twisted form~$K^{\zeta(\tau)}_G(G)$,
of the equivariant $K$-theory of $G$ acting on itself by conjugation.
Under this isomorphism the fusion product, when it is defined,
corresponds to the Pontryagin product.  The twisting $\zeta(\tau)$ is
given in terms of the level
\[
\zeta(\tau)= \mathfrak g + \dualc + \tau,
\]
where $\dualc$ is the ``dual Coxeter'' twisting.
\end{thm}

Several aspects of this theorem require clarification.  The main new
element is the ``twisted form'' of $K$-theory.  Twisted $K$-theory was
introduced by Donovan and Karoubi~\cite{donovan70:_graded_brauer_k} in
connection with the Thom isomorphism, and generalized and further
developed by Rosenberg~\cite{rosenberg89:_contin}.  Interest in
twisted $K$-theory was rekindled by its appearance in the late
1990's~\cite{minasian97:_k_ramon_ramon,witten98:_d_k} in string
theory.  Our results came about in the wake of this revival when we
realized that the work of the first author~\cite{freed94:_higher} on
Chern-Simons theory for finite groups could be interpreted in terms
twisted $K$-theory.

The twisted forms of $G$-equivariant $K$-theory are classified by the
nerve of the category of invertible modules over the equivariant
$K$-theory spectrum $K_{G}$.  What comes up in geometry though, is
only a small subspace, and throughout this paper the term ``twisting''
will refer twistings in this restricted, more geometric class.  These
geometric twistings of $K_{G}$-theory on a $G$-space $X$ are
classified up to isomorphism by the set
\begin{equation}\label{eq:43}
H^{0}_{G}(X;\Z/2)
\times H^{1}_{G}(X;\Z/2)
\times H^{3}_{G}(X;\Z).
\end{equation}
The component in $H^{0}$ corresponds to the ``degree'' of a $K$-class,
and the fact that the coefficients are the integers modulo $2$ is a
reflection of Bott periodicity.  In this sense ``twistings'' refine the
notion of ``degree,'' though when considering twistings it is
important to remember more than just the isomorphism class.

The tensor product of $K_{G}$-modules makes these spaces of twistings
into infinite loop spaces and provides a commutative group structure
on the sets of isomorphism classes.  The group structure
on~\eqref{eq:43} is the product of $H^{0}_{G}(X;\Z/2)$ with with
the extension of $H^{1}_{G}(X;\Z/2)$ by $H^{3}_{G}(X;\Z)$ with cocycle
$\beta(x\cup y)$, where $\beta$ is the Bockstein homomorphism.

A twisting is a form of equivariant $K$-theory on a space.  A level
for the loop group, on the other hand, corresponds to a central
extension of $LG$.  One way of relating these two structures is via
the classification~\eqref{eq:43}.  A central extension of $LG$ has a
topological invariant $H^{3}_{G}(G)$, and so give rise to a twisting,
up to isomorphism.  When the group $G$ is simple and simply connected
this invariant determines the central extension up to isomorphism,
the group $H^{1}_{G}(G)$ vanishes, and there is a canonical
isomorphism $H^{3}_{G}(G)\approx \Z$.  In this case an integer can be
used to specify both a level and a twisting.  There is a more refined
version of this correspondence directly relating twistings to central
extension, and the approach to twistings we take in this paper is
designed to make this relationship as transparent as possible.

There is a map from vector bundles to twistings which associates to a
vector bundle $V$ over $X$ the family of $K$-modules $K^{\bar V_{x}}$,
where $\bar V_{x}$ is the one point compactification of the fiber of
$V$ over $x\in X$, and for a space $S$, $K^{S}$ is the $K$-module with
$\pi_{0}K^{S}=K^{0}(S)$.  We denote the twisting associated to $V$ by
$\tau_{V}$, though when no confusion is likely to arise we will just
use the symbol $V$.  The invariants of $\tau_{V}$ in~\eqref{eq:43} are
$\dim V$, $w_{1}(V)$ and $\beta w_{2}(V)$.  These twistings are
described by Donovan and Karoubi in~\cite{donovan70:_graded_brauer_k}
from the point of view of Clifford algebras.

There is also a homomorphism from $KO_{G}^{-1}(X)$ to the group of
twistings of equivariant $K$-theory on $X$.  In topological terms it
corresponds to the map from the stable orthogonal group $O$ to its
third Postnikov section $O\langle 0,\dots,3\rangle$.  It sends an
element of $KO^{-1}$ to the twisting whose components are $(\sigma
w_{1},\sigma w_{2}, x)$, where $\sigma:H^{\ast}(BO)\to H^{\ast-1}(O)$
is the cohomology suspension, and $x\in H^{3}(O;\Z)$ is the unique
element, twice which is the cohomology suspension of $p_{1}$.  In
terms of operator algebras this homomorphism sends a skew-adjoint
Fredholm operator to its graded Pfaffian gerbe.  We will call this map
$\pfaff$.

We now describe two natural twistings on $G$ which are equivariant for
the adjoint action.  The first comes from the adjoint representation
of $G$ regarded as an equivariant vector bundle over a point, and
pulled back to $G$.  We'll write this twisting as $\mathfrak g$.  For
the other, first note that the equivariant cohomology group
$H^{\ast}_{G}(G)$ is just $H^{\ast}(LBG).$ The vector bundle
associated to the adjoint representation gives a class $\ad\in
KO^{0}(BG)$ which we can transgress to $KO^{-1}(LBG)$, and then map to
twistings by the map $\pfaff$.  We'll call this twisting $\dualc$.
When $G$ is simple and simply connected, the integer corresponding to
$\dualc$ is the dual Coxeter number, and $\mathfrak g$ is just a degree
shift.  With these definitions, the formula
\[
\zeta(\tau)=\mathfrak g + \dualc + \tau
\]
in the statement of Theorem~\ref{thm:4} should be clear.  The
twistings $\mathfrak g$ and $\dualc$ are those just described, and $\tau$
is the twisting corresponding to the level.

Theorem~\ref{thm:4} provides a topological description of the group
$\rep^{\tau}(LG)$ and its fusion product when it exists.  But it also
gives more.  The twisted $K$-group $K^{\zeta(\tau)}_{G}(G)$ is
defined for any compact Lie group $G$, and it makes sense for any
level $\tau$.  This points the way to a formulation of an analogue of
the group $\rep^{\tau}(LG)$ for {\em any} compact Lie group $G$ (even
one which is finite).  In Parts II and III we take up these
generalizations and show that the assertions of Theorem~\ref{thm:4}
remain true.

One thing that emerges from our topological considerations is the need
to consider {\em $\Z/2$-graded} central extensions of loop groups.
Such extensions are necessary when working with a group like $SO(3)$
whose adjoint representation is not $\spin$.  Another interesting case
is that of $O(2)$.  When the adjoint representation is not orientable
the dual Coxeter twist makes a non-trivial contribution to
$H^{0}_{G}(G;\Z/2)$.  One sees an extra change in the degree in which
the interesting $K$-group occurs.  In the case of $O(2)$ these degree
shifts are different on the two connected components, again
emphasizing the point that twistings should be regarded as a
generalization of degree.  The ``Verlinde ring'' in this case is
comprised of an even $K$-group on one component and an odd $K$-group
on another.  Such inhomogeneous compositions are not typically
considered when discussing ordinary $K$-groups.

The fusion product on $\rep^{\tau}(LG)$ has been defined for simple
and simply connected $G$, and in a few further special cases.  The
Pontryagin product on $K^{\zeta(\tau)}_{G}(G)$ is defined exactly
when $\tau$ is {\em primitive} in the sense that its pullback along
the multiplication map of $G$ is isomorphic to the sum of its
pullbacks along the two projections.  This explains, for example, why
a fusion product on $\rep^{\tau}(LSO(n))$ exists only at half of the
levels.  Using the Pontryagin product we are able to define a fusion
product on $\rep^{\tau}(LG)$ for any $G$ at any primitive level
$\tau$.  We do not, however, give a construction of this product in
terms of representation theory.

When the fusion product is defined on $\rep^{\tau}(LG)$ it is part of
a much more elaborate structure.  For one thing, there is a trace map
$\rep^{\tau}(LG)\to\Z$ making $\rep^{\tau}(LG)$ into a Frobenius
algebra.  Using twisted $K$-theory we construct this trace map for
general compact Lie groups at primitive levels $\tau$ which are
non-degenerate in the sense that the image of $\tau$ in
\[
H^{3}_{T}(T;\R)\approx H^{1}(T;\R)\otimes H^{1}(T;\R)
\]
is a non-degenerate bilinear form.  Again, in the cases when the
fusion product has been defined, there are operations on
$\rep^{\tau}(LG)$ coming from the moduli spaces of Riemann surfaces
with boundary, making $\rep^{\tau}(LG)$ part of is often called a {\em
topological conformal field theory}.  Using topological methods, we
are able to construct a topological conformal field theory for any
compact Lie group $G$, at levels $\tau$ which are transgressed from
(generalized) cohomology classes on $BG$ and which are non-degenerate.
Some of this work appears in~\cite{math.AT/0206257,freed:_consis}.

Another impact of Theorem~\ref{thm:4} is that it brings the
computational techniques of algebraic topology to bear on the
representations of loop groups.  One very interesting approach, for
connected $G$, is to use the Rothenberg-Steenrod spectral sequence
relating the equivariant $K$-theory of $\Omega G$ to that of $G$.  In
this case one gets a spectral sequence
\begin{equation}\label{eq:44}
\tor^{K^{G}_{\ast}(\Omega G)}(R(G), R(G))\implies K^{G}_{\tau+\ast}(G),
\end{equation}
relating the untwisted equivariant $K$-homology of $\Omega G$, and the
representation ring of $G$ to the Verlinde algebra.  The ring
$K^{G}_{\ast}(\Omega G)$ can be computed using the techniques of
Bott~\cite{bott58:_lie} and Bott-Samelson~\cite{bott58:_applic_morse}
and has also been described by Bezrukavnikov, Finkelberg and
Mirkovi\'c~\cite{bezrukavnikov05:_equiv_k_grass_toda}.  The $K$-groups
in the $E^{2}$-term are untwisted.  The twisting appears in the way
that the representation ring $R(G)$ is made into an algebra over
$K^{G}(\Omega G)$.  The equivariant geometry of $\Omega G$ has been
extensively studied in connection with the representation theory of
$LG$, and the spectral sequence~\eqref{eq:44} seems to express yet
another relationship.  We do not know of a representation theoretic
construction of~\eqref{eq:44}.  An analogue of the spectral
sequence~\eqref{eq:44} has been used by Chris
Douglas~\cite{douglas06:_k_lie} to compute the (non-equivariant)
twisted $K$-groups $K^{\tau}(G)$ for all simple, simply connected $G$.

Using the Lefschetz fixed point formula one can easily conclude 
for connected $G$ that
\[
\Delta^{-1}K^{G}_{\ast}(\Omega G)= \Delta^{-1}\Z[\Lambda\times \Pi],
\]
where $\Delta$ is the square of the Weyl denominator,
$\Pi=\pi_{1}T$ is the co-weight lattice, and $\Lambda$ is the weight
lattice.  When the level $\tau$ is non-degenerate there are no higher
$\tor$ groups, and the spectral sequence degenerates to an isomorphism
\[
\Delta^{-1}K^{G}_{\tau+\ast}(G) \approx \Delta^{-1}R(G)/I^{\tau}
\]
where $I^{\tau}$ is the ideal of representations whose characters
vanish on certain conjugacy classes.  The main computation of this
paper asserts that such an isomorphism holds without inverting
$\Delta$ when $G$ is connected, and $\pi_{1}G$ is torsion free.  The
distinguished conjugacy classes are known as {\em Verlinde conjugacy
classes}, and the ideal $I^{\tau}$ as the {\em Verlinde ideal}.
In~\cite{math.AT/0206257} the ring $K^{\zeta(\tau)}_{G}(G)\otimes\C$
is computed using a fixed point formula, and shown to be isomorphic to
the Verlinde algebra.

The plan of this series of papers is as follows.  In Part I we define
twisted $K$-groups, and compute the groups $K^{\zeta}_{G}(G)$ for
connected $G$ with torsion free fundamental group, at non-degenerate
levels $\zeta$.  Our main result is Theorem~\ref{thm:46}.  In Part II
we introduce a certain family of Dirac operators and our
generalization of $\rep^{\tau}(LG)$ to arbitrary compact Lie groups.
We construct a map from $\rep^{\tau}(LG)$ to $K^{\zeta(\tau)}_{G}(G)$
and show that it is an isomorphism when $G$ is connected with torsion
free fundamental group.  In Part III we show that our map is an
isomorphism for general compact Lie groups $G$, and develop some
applications.

The bulk of this paper is concerned with setting up twisted
equivariant $K$-theory.  There are two things that make this a little
complicated.  For one, when working with twistings it is important to
remember the morphisms between them, and not just the isomorphism
classes.  The twistings on a space form a {\em category} and spelling
out the behavior of this category as the space varies gets a
little elaborate.  The other thing has to do with the kind of
$G$-spaces we use.  We need to define twistings on $G$-equivariant
$K$-theory in such a way as to make clear what happens as the group
$G$ changes, and for the constructions in Part II we need to make the
relationship between twistings and (graded) central extensions as
transparent as possible.  We work in this paper with {\em groupoids}
and define twisted equivariant $K$-theory for groupoids.  Weakly
equivalent groupoids (see Appendix~\ref{sec:groupoids}) have
equivalent categories of twistings and isomorphic twisted $K$-groups.
A group $G$ acting on a space $X$ forms a special kind of groupoid
$X\modmod G$ called a ``global quotient groupoid.''  A central
extension of $G$ by $U(1)$ defines a twisting of $K$-theory for
$X\modmod G$.  If $G$ is a compact connected Lie group acting on
itself by conjugation, and $PG$ denotes the space of paths in $G$
starting at the identity, acted upon by $LG$ by conjugation, then 
\[
PG\modmod LG \to G\modmod G 
\]
is a local equivalence, and a (graded) central extension of $LG$
defines a twisting of $PG\modmod LG$ and hence of $G\modmod G$.  In
general we will define a twisting of a groupoid $X$ to consist of a
local equivalence $P\to X$ and a graded central extension $\tilde P$
of $P$.  While most of the results we prove reduce, ultimately, to
ordinary results about compact Lie groups acting on spaces, not all
do.  In Part III it becomes necessary to work with groupoids which are
not equivalent to a compact Lie group acting on a space.

Nitu Kitchloo~\cite{kitchloo:_domin_k_kac_moody} has pointed out that
the space $PG$ is the universal $LG$ space for proper actions.  Using
this he has described a generalization of our computation to other
Kac-Moody groups.

At the time we began this work, the paper of Atiyah and
Segal~\cite{atiyah:_twist_k} was in preparation, and we benefited a
great deal from early drafts.  Since that time several other
approaches to twisted $K$-theory have appeared.  In addition
to~\cite{atiyah:_twist_k} we refer the reader
to~\cite{bouwknegt02:_twist_k_k,tu04:_twist_k}.  We have chosen to use
``graded central extensions'' because of the close connection with
loop groups and the constructions we wish to make in Part II.  Of
course our results can be presented from any of the points of view
mentioned above, and the choice of which is a matter of personal
preference.

We have attempted to organize this paper so that the issues of
implementation are independent of the issues of computation.
Section~\ref{sec:twisted-k-theory-1} is a kind of field guide to
twisted $K$-theory.  We describe a series of examples intended to give
the reader a working knowledge of twisted $K$-theory sufficient to
follow the main computation in \S\ref{sec:comp-k_gg-twist}.
Section~\ref{sec:twistings} contains our formal discussion of
twistings of $K$-theory for groupoids, and our definition of twisted
$K$-groups appears in \S\ref{sec:twisted-k-groups}.  We have attempted
to axiomatize the theory of twisted $K$-groups in order to facilitate
comparison with other models.  Our main computation appears in
\S\ref{sec:comp-k_gg-twist}.

This paper has been a long time in preparation and we have benefited
from discussion with many people.  The authors would like to thank Sir
Michael Atiyah and Graeme Segal for making available early drafts
of~\cite{atiyah:_twist_k} and~\cite{math/0510674v1}.  As will be
evident to the reader, our approach to twisted $K$-theory relies
heavily on their ideas.  We would like to thank Is Singer for many
useful conversations.  We would also like to thank Ulrich Bunke and
Thomas Schick for their careful study and comments on this work.  The
report of their seminar appears in~\cite{bunke05:_twist_k_tqft}.

We assume throughout this paper that all spaces are locally
contractible, paracompact and completely regular.  These assumptions
implies the existence of partitions of unity~\cite{dold63:_partit} and
locally contractible slices through actions of compact Lie
groups~\cite{mostow57:_equiv_euclid,palais61:_lie}.

\section{Twisted $K$-theory by example}
\label{sec:twisted-k-theory-1}
\numberwithin{equation}{section}

The $K$-theory of a space is assembled from data which is local.  To
give a vector bundle $V$ on $X$ is equivalent to giving vector bundles
$V_{i}$ on the open sets $U_{i}$ of a covering, and isomorphisms
\[
\lambda_{ij}:V_{i}\to V_{j}
\] 
on $U_{i}\cap U_{j}$ satisfying a compatibility (cocycle) condition on
the triple intersections.  In terms of $K$-theory this is expressed by the
Mayer-Vietoris (spectral) sequence relating $K(X)$ and the $K$-groups
of the intersections of the $U_{i}$.  In forming twisted $K$-theory we
modify this descent or gluing datum, by introducing a line bundle
$L_{ij}$ on $U_{i}\cap U_{j}$, and asking for an isomorphism
\[
\lambda_{ij}:L_{ij}\otimes V_{i}\to V_{j}
\] 
satisfying a certain cocycle condition.     In terms of $K$-theory,
this modifies the restriction maps in the Mayer-Vietoris sequence.

In order to formulate the cocycle condition, the $L_{ij}$ must come
equipped with an isomorphism
\[
L_{jk}\otimes L_{ij}\to L_{ik}
\]
on the triple intersections, satisfying an evident compatibility
relation on the quadruple intersections.  In other words, the
$\{L_{ij} \}$ must form a $1$-cocycle with values in the groupoid of
line bundles.  Cocycles differing by a $1$-cochain give isomorphic
twisted $K$-groups, so, up to isomorphism, we can associate a twisted
notion of $K(X)$ to an element 
\[
\tau\in H^{1}(X;\{\text{Line Bundles} \}).
\]

On good spaces there are isomorphisms
\[
H^{1}(X;\{\text{Line Bundles} \}) \approx
H^{2}(X;U(1)) \approx
H^{3}(X;\Z),
\]
and correspondingly, twisted notions of $K$-theory associated to an
integer valued $3$-cocycle.  In this paper we find we need to allow
the $L_{ij}$ to be $\pm$ line bundles, so in fact we consider
twisted notion of $K(X)$ classified by elements%
\footnote{The set of isomorphism classes of twistings has a group
structure induced from the tensor product of graded line bundles.
While there is, as indicated, a set-theoretic factorization of the
isomorphism classes of twistings, the group structure is not, in
general, the product.}%
\[
\tau\in H^{1}(X;\{\text{$\pm$Line Bundles} \})\approx
H^{3}(X;\Z)\times H^{1}(X;\Z/2).
\]
We will write $K^{\tau+n}$ for the version of $K^{n}$, twisted by
$\tau$.

In practice, to compute twisted $K(X)$ one represents the twisting
$\tau$ as a Cech $1$-cocycle on an explicit covering of $X$.  The
twisted $K$-group is then assembled from the Mayer-Vietoris sequence
of this covering, involving the same (untwisted) $K$-groups one would
encounter in computing $K(X)$.  The presence of the $1$-cocycle is
manifest in the restriction maps between the $K$-groups of the open
sets.  They are modified on the two-fold intersections by tensoring
with the ($\pm$) line bundle given by the $1$-cocycle.  This
``operational definition'' suffices to make most computations.  See
\S\ref{sec:twisted-k-groups} for a more careful discussion.  Here
are a few examples.

\begin{eg}
\label{eg:23} Suppose that $X=S^{3}$, and that the isomorphism class
of $\tau$ is $n\in H^{3}(X;\Z)\approx \Z$.  Let $U_{+}= X\setminus
(0,0,0,-1)$ and $U_{-}=X\setminus (0,0,0,1)$.  Then $U_{+}\cap
U_{-}\sim S^{2}$, and $\tau$ is represented by the $1$-cocycle whose
value on $U_{+}\cap U_{-}$ is $L^{n}$, with $L$ the tautological line
bundle.  The Mayer-Vietoris sequence for $K^{\tau}(X)$ takes the form
\begin{multline*}
\dots\to K^{\tau+0}(X)\to 
K^{\tau+0}(U_{+})\oplus
K^{\tau+0}(U_{-})\to 
K^{\tau+0}(U_{+}\cap U_{-}) \\
\to 
K^{\tau+1}(X)\to 
K^{\tau+1}(U_{+})\oplus
K^{\tau+1}(U_{-})\to 
K^{\tau+1}(U_{+}\cap U_{-})\to \cdots .
\end{multline*}
Since the restriction of $\tau$ to $U_{\pm}$ is isomorphic to zero, we
have
\begin{align*}
K^{\tau+0}(U_{\pm})& \approx 
K^{0}(U_{\pm})\approx \Z \\
K^{\tau+1}(U_{\pm})& \approx 
K^{1}(U_{\pm})\approx 0 \\
K^{\tau+1}(U_{+}\cap U_{-})& \approx 
K^{1}(S^{2}) \approx 0,\\
\end{align*}
and
\[
K^{\tau+0}(U_{+}\cap U_{-}) \approx 
K^{0}(S^{2}) \approx \Z\oplus\Z.
\]
with basis the trivial bundle $1$, and the tautological line bundle
$L$.  The Mayer-Vietoris sequence reduces to the exact sequence
\[
0\to K^{\tau+0}(X)\to \Z\oplus\Z\to \Z\oplus\Z \to K^{\tau+1}(X)\to 0.
\]

In ordinary (untwisted) $K$-theory, the middle map is
\[
\begin{pmatrix}
1 & -1 \\
0 & 0
\end{pmatrix}: \Z\oplus\Z \to \Z\oplus\Z.
\]
In twisted $K$-theory, with suitable conventions, the middle map
becomes
\begin{equation}\label{eq:45}
\begin{pmatrix}
1 & n-1 \\
0 & -n
\end{pmatrix}: \Z\oplus\Z \to \Z\oplus\Z,
\end{equation}
and so
\[
K^{\tau+n}(S^{3}) = \begin{cases}
0 &\quad n=0 \\
\Z/n &\quad n=1.
\end{cases}
\]
In the language of twistings, the map~\eqref{eq:45} is accounted for
as follows.  To  identify the twisted $K$-groups with ordinary twisted
$K$-groups we have to choose isomorphisms 
\[
t_{\pm}:\tau\vert_{U_{\pm}}\to 0.
\]
If we use the $t_{+}$ to trivialize $\tau$ on $U_{+}\cap U_{-}$,
then the following diagram commutes
\[
\xymatrix{
K^{\tau+\ast}(U_{+})\ar[r]^{t_{+}}\ar[d]_{\text{restr.}} & K^{0+\ast}(U_{+})\ar[d]^{\text{restr.}} \\
K^{\tau+\ast}(U_{+}\cap U_{-})\ar[r]_{t_{+}} & K^{0+\ast}(U_{+}\cap U_{-}) \\
}
\]
and we can identify the restriction map in twisted $K$-theory from
$U_{+}$ to $U_{+}\cap U_{-}$ with the restriction map in untwisted
$K$-theory.  On $U_{+}\cap U_{-}$ we have
$t_{-}=(t_{-}t_{+}^{-1})\circ t_{+}$ so the restriction map in twisted
$K$-theory is identified with the restriction map in untwisted
$K$-theory, followed by the map $(t_{-}t_{+}^{-1})$.  By definition of
$\tau$, this map is given by multiplication by $L^{n}$.  This 
accounts for the second column of~\eqref{eq:45}.
\end{eg}

\begin{eg}
\label{eg:24} Now consider the twisted $K$-theory of $U(1)$ acting
trivially on itself.    In this case the twistings are classified by 
\[
H^{3}_{U(1)}(U(1);\Z)\times
H^{1}_{U(1)}(U(1);\Z/2)\approx 
\Z\oplus \Z/2.
\]
We consider twisted $K$-theory, twisted by $\tau=(n,\epsilon)$.
Regard $U(1)$ as the unit circle in the complex plane, and set 
\[
U_{+}=U(1)\setminus\{-1 \}\quad
U_{-}=U(1)\setminus\{+1 \}.
\]
The twisting $\tau$ restricts to zero on both $U_{+}$ and $U_{-}$.
Write
\[
K^{0}_{U(1)}=R(^{1})=\Z[L,L^{-1}].
\]
Then the Mayer-Vietoris sequence becomes
\[
0\to K^{\tau+0}_{U(1)}(U(1))\to 
\Z[L^{\pm 1}]\oplus
\Z[L^{\pm 1}] \to 
\Z[L^{\pm 1}]\oplus
\Z[L^{\pm 1}] \to 
K^{\tau+1}_{U(1)}(U(1))\to 0.
\]
The $1$-cocycle representing $\tau$ can be taken to be the equivariant
vector bundle whose fiber over $-i$ is the trivial representation of
$U(1)$ and whose fiber over $+i$ is $(-1)^{\epsilon}L^{n}$.  With
suitable conventions, the middle map becomes
\[
\begin{pmatrix}
1& -(-1)^{\epsilon} L^{n} \\
1 & -1 
\end{pmatrix}: \Z[L^{\pm1}]^{2} \to \Z[L^{\pm1}]^{2}.
\]
It follows that
\[
K^{\tau+k}_{U(1)}(U(1)) =
\begin{cases}
0 &\quad k=0 \\
\Z[L^{\pm1}] / \left((-1)^{\epsilon}L^{n}-1 \right) &\quad k=1.
\end{cases}
\]
When $\epsilon=0$ this coincides the Grothendieck group of
representations of the Heisenberg extension of $\Z\times U(1)$ of
level $n$, and in turn with the Grothendieck group of positive energy
representations the loop group of $U(1)$ at level $n$.
\end{eg}

\begin{eg}
\label{eg:25} Consider the twisted $K$-theory of $SU(2)$ acting on
itself by conjugation.  The group $H^{1}_{SU(2)}(SU(2);\Z/2)$
vanishes, while
\[
H^{3}_{SU(2)}(SU(2);\Z)=\Z,
\]
so a twistings $\tau$ in this case is given by an integer
\[
n\in H^{3}_{SU(2)}(SU(2);\Z)=\Z.
\]
Set
\[
U_{+}=SU(2)\setminus\{-1 \}\quad
U_{-}=SU(2)\setminus\{+1 \}.
\]
The spaces $U_{+}$ and $U_{1}$ are equivariantly contractible, while
$U_{+}\cap U_{-}$ is equivariantly homotopy equivalent to $S^{2}=
SU(2)/T$, where $T=U(1)$ is a maximal torus.  The restrictions of
$\tau$ to $U_{+}$ and $U_{-}$ are isomorphic to zero.  We have
\[
K^{0}_{SU(2)}(U_{\pm})\approx
K^{0}_{SU(2)}(\point)=
R(SU(2))=\Z[L,L^{-1}]^{W}
\]
with the Weyl group $W\approx \Z/2$ acting by exchanging $L$ and
$L^{-1}$, and
\[
K^{0}_{SU(2)}(U_{+}\cap U_{-})\approx
K^{0}_{SU(2)}(SU(2)/T) \approx
K^{0}_{T}(\point) \approx \Z[L,L^{-1}].
\]
The ring $R(SU(2))$ has an additive basis consisting of the
irreducible representations,
\[
\rho_{k}= L^{k}+ L^{k-2}+\dots + L^{-k},\qquad k\ge 0
\]
which multiply according to the Clebsch-Gordon rule
\[
\rho_{l}\,\rho_{k} = \rho_{k+l}+\rho_{k+l-2}+\dots + \rho_{k-l}\qquad k
\ge l.
\]
As in our other example, the Mayer-Vietoris sequence is short exact
\[
0\to K^{\tau+0}_{SU(2)}(SU(2))\to 
R(SU(2))\oplus
R(SU(2))\to 
\Z[L^{\pm 1}] \to 
K^{\tau+1}_{SU(2)}(SU(2))\to 0.
\]
The $1$-cocycle representing the difference between the two
trivializations of the restriction of $\tau$ to $U_{+}\cap U_{-}$ can
be taken to be the element $L^{n}\in K_{SU(2)}(SU(2)/T)\approx R(T)$.
The sequence is a sequence of $R(SU(2))$-modules.  With suitable
conventions, the middle map is
\begin{equation}\label{eq:28}
\begin{pmatrix}
1 & -L^{n}
\end{pmatrix}: R(SU(2))^{2} \to R(T).
\end{equation}
To calculate the kernel and cokernel, note that $R(T)$ is a free
module of rank $2$ over $R(SU(2))$.  We give $R(SU(2))\oplus
R(SU(2)))$ the obvious basis, and $R(T)$ the basis $\{1,L \}$.
It follows from the identity
\[
L^{n} = L\,\rho_{n-1} - \rho_{n-2}
\]
that~\eqref{eq:28} is represented by the matrix
\[
\begin{pmatrix}
1 & \rho_{n-2}\\
0 & -\rho_{n-1}
\end{pmatrix},
\]
and that
\[
K^{\tau+k}_{SU(2)}(SU(2)) =
\begin{cases}
0 &\quad k=0 \\
R(SU(2))/(\rho_{n-1}) &\quad k=1.
\end{cases}
\]
This coincides the Grothendieck group of
positive energy representations the loop group of $SU(2)$ at level
$(n-2)$.
\end{eg}

Examples~\ref{eg:24} and~\ref{eg:25} illustrate the relationship
between twisted $K$-theory and the representations of loop groups.  In
both cases the Grothendieck group of positive energy representations
of the loop group of a compact Lie group $G$ is described by the
twisted equivariant $K$-group of $G$ acting on itself by conjugation.
Two minor discrepancies appear in this relationship.  On one hand, the
interesting $K$-group is in degree $k=1$.  As explained in the
introduction, the representations of the loop group at level $\tau$
correspond to twisted $K$-theory at the twisting $\zeta(\tau)= \mathfrak g
+ \dualc + \tau$.  The shift in $K$-group to $k=1$ corresponds to the term
$\mathfrak g$.  In both examples the adjoint representation is $\spin^{c}$
and so contributes only its dimension to $\zeta(\tau)$.  This term
could be gotten rid of by working with twisted equivariant
$K$-homology rather than $K$-cohomology.  We have chosen to work with
$K$-cohomology in order to make better contact with our geometric
constructions in Parts II and III.  The other discrepancy is the shift
in level in Example~\ref{eg:25}: twisted equivariant $K$-theory at
level $n$ corresponds to the representations of the loop group at
level $(n-2)$.  The shift of $2$ here corresponds to the term $\dualc$ in
our formula for $\zeta(\tau)$.

We now give a series of examples describing other ways in which
twistings of $K$-theory arise.

\begin{eg} 
Let $V$ be a vector bundle of dimension $n$ over a space $X$, and
write $X^{V}$ for the Thom complex of $V$.  Then $\tilde
K^{n+k}(X^{V})$ is a twisted form of $K^{k}(X)$.  To identify the
twisting, choose local $\spinc$ structures $\mu_{i}$ on the
restrictions $V_{i}=V\vert_{U_{i}}$ of $V$ to the sets in  an open cover of $X$.
The $K$-theory Thom classes associated to the $\mu_{i}$ allow one to identify
\[
\tilde
K^{n+k}(U_{i}^{V_{i}})\approx K^{k}(U_{i})
\]
The difference between the two identifications on $U_{i}\cap U_{j}$ is
given by multiplication by the graded line bundle representing the
difference between $\mu_{i}$ and $\mu_{j}$.  Figuratively, the cocycle
representing the twisting is $\mu_{j}\mu_{i}^{-1}$ and the cohomology
class is $(w_{1}(V), W_{3}(V))\in H^{1}(X;\Z/2)\times H^{3}(X;\Z)$,
where $W_{3}=\beta w_{2}$.  This is one of the original examples of
twisted $K$-theory, described by Donovan and
Karoubi~\cite{donovan70:_graded_brauer_k} from the point of view of
Clifford algebras.  We will review their description in
\S\ref{sec:thom-isomorphism}.
\end{eg}

\begin{eg}\label{eg:27}
Let $G$ be a compact Lie group.  The central extensions
\[
\T\to \tilde G\xrightarrow{\tau} G
\]
of $G$ by $\T=U(1)$ are classified by
$H^{3}_{G}(\{\point\};\Z)=H^{3}(BG;\Z)$.  The Grothendieck group
$\rep^{\tau}(G)$, of representations of $\tilde G$ on which $\T$ acts
according to its defining representation, can be thought of as a
twisted form of $\rep(G)$.  In this case, our definition of
equivariant twisted $K$-theory gives
\[
K^{\tau+k}_{G}(\point) = 
\begin{cases}
\rep^{\tau}(G) &\qquad k=0 \\
0 &\qquad k=1
\end{cases}
\]
More generally, if $S$ is a $G$-space, and $\tau\in H^{3}_{G}(S)$ is
pulled back from $H^{3}_{G}(\point)$, then $K^{\tau+k}_{G}(S)$ is the
summand of $K^{k}_{\tilde G}(S)$ corresponding to $\tilde
G$-equivariant vector bundles on which $\T$ acts according to its
defining character.
\end{eg}

\begin{eg}\label{eg:26}
Now suppose that
\[
H\to G\to Q
\]
is an extension of groups, and $V$ is an irreducible representation of
$H$ that is stable, up to isomorphism, under conjugation by elements
of $G$.  Then Grothendieck group of representations of $G$ whose
restriction to $H$ is $V$-isotypical, forms a twisted version of the
Grothendieck group of representations of $Q$.  When $H$ is central,
equal to $\T$, and $V$ is the defining representation, this is the
situation of Example~\ref{eg:27}.  We now describe how to reduce to
this case.

Fix an $H$-invariant Hermitian metric on $V$, and write
$V^{\ast}=\hom(V,\C)$ for the representation dual to $V$.  Let $\tilde
G$ denote the group of pairs
\[
(g,f)\in G\times\hom(V^{\ast},V^{\ast})
\]
for which $f$ is unitary, and satisfies
\[
f(h v) = g h g^{-1} f(v)\qquad h \in H.
\]
Since $V$ is irreducible, and $(\ad g)^{\ast}V\approx V$, the same is
true of $V^{\ast}$, and the map
\begin{align*}
\tilde G &\to G \\
(g,f) &\mapsto g
\end{align*}
is surjective, with kernel $\T$.   The inclusion
\begin{align*}
H &\subset \tilde G \\
h &\mapsto (h,\text{action of }h) 
\end{align*}
is normal, and lifts the inclusion of $H$ into $G$.  We define 
\[
\tilde Q = \tilde G/H.
\]
The group $\tilde Q$ is a central extension of $Q$ by $\T$, which we
denote
\[
\tilde Q\xrightarrow{\tau} Q.
\] 
We now describe an equivalence of categories between $V$-isotypical
$G$-representations, and $\tau$-projective representation of $Q$
(representations of $\tilde Q$ on which $\T$ acts according to its
defining character).

By definition, the representations $V$ and $V^{\ast}$ of $H$ come
equipped with extensions to unitary representation of $\tilde G$.
Given a $V$-isotypical representation $W$ of $G$, we let $M$ denote
the $H$-invariant part of $V^{\ast}\otimes W$:
\[
M = \left(V^{\ast}\otimes W \right)^{H}.
\]
The action of $\tilde G$ on $V^{\ast}\otimes W$ factors through an
action of $\tilde Q$ on $M$.   This defines a functor from
$V$-isotypical representations of $G$ to $\tau$-projective
representations of $Q$.  

Conversely, suppose $M$ is a $\tau$-projective representation of
$\tilde Q$.  Let $\tilde G$ act on $M$ through the projection $\tilde
G\to \tilde Q$, and form
\[
W = V\otimes M.
\]
The central $\T$ of $\tilde G$ acts trivially on $W$, giving $W$ a
$G$-action.  This defines a functor from the category of
$\tau$-projective representations of $Q$ to the category of
$V$-isotypical representations of $G$.  One easily checks these two
functors to form an equivalence of categories.
\end{eg}

\begin{eg}
\label{eg:28}
Continuing with the situation of Example~\ref{eg:26} consider an
extension
\[
H\to G\to Q
\]
and an irreducible representation $V$ of $H$, which this time is not
assumed to be stable under conjugation by $G$.   Write 
\[
G_{0}=\left\{g\in G \mid (\ad g )^{\ast} V\approx V\right\},
\]
and $Q_{0}= G_{0}/H$.  Let $S$ be the set of isomorphism classes of
irreducible representations of $H$ of the form $(\ad g)^{\ast}V$.  The
conjugation action of $G$ on $S$ factors through $Q$, and we have an
identification $S=Q/Q_{0}$.  Let's call a representation of $G$ {\em
$S$-typical} if its restriction to $H$ involves only the
irreducible representations in $S$.  One easily checks that
``induction'' and ``passage to the $V$-isotypical part of the
restriction'' give an equivalence of categories
\[
\left\{S\text{-typical representations of }G \right\}
\leftrightarrow
\left\{V\text{-isotypical representations of }G_{0} \right\},
\]
and therefore an isomorphism of the Grothendieck group $\rep^{S}(G)$
of $S$-typical representations of $G$ with
\[
\rep^{\tau}(Q_{0})\approx K^{\tau}_{Q_{0}}(\point).
\]
We can formulate this isomorphism a little more cleanly in the
language of groupoids.  For each $\alpha\in S$, choose an irreducible
$H$-representation $V_{\alpha}$ representing $\alpha$.  Consider the
groupoid $S\modmod Q$, with set of objects $S$, and in which a
morphism $\alpha\to \beta$ is an element $g\in Q$ for which $(\ad
g)^{\ast}\alpha = \beta$.  We define a new groupoid $P$ with objects
$S$, and with $P(\alpha,\beta)$ the set of equivalences classes of
pairs $(g, \phi)\in G\times\hom(V_{\alpha^{\ast},V_{\beta}}^{\ast})$,
with $\phi$ unitary, and satisfying
\[
\phi(h\,v)= ghg^{-1}\phi(v)
\]
(so that, among other things, $(\ad g)^{\ast} \alpha = \beta$).    The
equivalence relation is generated by 
\[
(g,\phi)\sim (h g, h\phi)\qquad h\in H.
\]
There is an evident functor $\tau:P\to S\modmod Q$, representing $P$
as a central extension of $S\modmod Q$ by $\T$.  The automorphism
group of $V$ in $P$ is the central extension $\tilde Q_{0}$ of
$Q_{0}$.  An easy generalization of the construction of
Example~\ref{eg:26} gives an equivalence of categories
\[
\left\{\tau-\text{projective representations of }P \right\}
\leftrightarrow
\left\{S-\text{typical representations of }G \right\}.
\]
Central extensions of $S\modmod Q$ 
are classified by 
\[
H^{3}(S\modmod Q;\Z) = H^{3}_{Q}(S;\Z) \approx
H^{3}_{Q_{0}}(\point;\Z),
\]
and so represent twistings of $K$-theory.  Our definition of twisted
$K$-theory of groupoids will identify the $\tau$-twisted $K$-groups of
$(S\modmod Q)$ with the summand of the $K$-theory of $P$ on which the
central $\T$ acts according to its defining representation.  We
therefore have an isomorphism
\[
\rep^{S}(G) \approx K^{\tau+0}(S\modmod Q)= K^{\tau+0}_{Q}(S).
\]
\end{eg}

\begin{eg}
\label{eg:29}
Now let $S$ denote the set of isomorphism classes of all irreducible
representations of $H$.  Decomposing $S$ into orbit types, and using
the construction of Example~\ref{eg:28} gives a central extension
$\tau:P\to \left(S\modmod Q \right)$, and an isomorphism
\[
\rep(G)\approx K^{\tau+0}(S\modmod Q)=K^{\tau+0}_{Q}(S).
\] 
More generally, if $X$ is a space with a $Q$-action there is an
isomorphism
\[
K_{G}^{k}(X)\approx K_{Q}^{\tau+k}(X\times S),
\]
in which $\tau$ is the $Q$-equivariant twisting of $X\times S$ pulled
back from the $Q$-equivariant twisting $\tau$ of $S$, given by $P$.  
\end{eg}

\section{Twistings of $K$-theory}\label{sec:twistings}

We now turn to a more careful discussion of twistings of $K$-theory.
Our terminology derives from the situation of Example~\ref{eg:27} in
which a central extension of a group gives rise to a twisted notion of
equivariant $K$-theory.  By working with graded central extensions of
groupoids (rather that groups) we are able to include in a single
point of view both the twistings that come from $1$-cocycles with
values in the group of $\Z/2$-graded line bundles and the twistings
that come from central extensions.  In order to facilitate this, in
the rest of this paper we will use the language of $\T$-bundles and
$\T$-torsors instead of ``line bundles,'' where $\T$ is the group
$U(1)$.  We begin with a formal discussion of ($\Z/2$-)graded
$\T$-bundles.

\subsection{Graded $\T$-bundles}
\label{sec:graded-u1-bundles}

Let $X$ be a topological space.  

\begin{defin}
A {\em graded $\T$-bundle} over $X$ consists of a principal
$\T$-bundle  $P\to X$, and a locally constant function
$\epsilon:X\to\Z/2$.
\end{defin}

We will call a graded $\T$-bundle $(P,\epsilon)$ {\em even}
(resp.~{\em odd}) if $\epsilon$ is the constant function function $1$
(resp.~$-1$).  The collection of graded $\T$-bundles forms a symmetric
monoidal groupoid.  A map of graded $\T$-bundles
$(P_{1},\epsilon_{1})\to(P_{1},\epsilon_{2})$ exists only when
$\epsilon_{1}=\epsilon_{2}$, in which case it is a map of principal
bundles $P_{1}\to P_{2}$.  The tensor structure is given by
\[
(P_{1},\epsilon_{1})\otimes
(P_{2},\epsilon_{2}) =
(P_{1}\otimes P_{2},\epsilon_{1}+\epsilon_{2}), 
\]
in which $P_{1}\otimes P_{2}$ is the usual ``tensor product'' of
principal $\T$-bundles: 
\[
\left(P_{1}\otimes P_{2} \right)_{x}= \left(P_{1}
\right)_{x}\times\left(P_{2}
\right)_{x}/(v\lambda,w)\sim(v,w\lambda).
\]
It is easiest to describe the symmetry transformation
\[
T:(P_{1},\epsilon_{1})\otimes (P_{2},\epsilon_{2})
\to (P_{2},\epsilon_{2})\otimes (P_{1},\epsilon_{1})
\]
fiberwise.  In the fiber over a point $x\in X$ it is 
\[
(v,w)\mapsto (w,v\epsilon_{1}(x)\epsilon_{2}(x)).
\]

We will write $\buonepm$ for the contravariant functor which
associates to a space $X$ the category of graded $\T$-bundles over
$X$, and for $\buone$ the functor ``category of $\T$-bundles''.  We
will also write $H^{1}(X;\tpm)$ for the group of isomorphism classes in
$\buonepm(X)$, and $H^{0}(X;\tpm)$ for the group of automorphisms of
any object.  There is an exact sequence
\begin{equation}\label{eq:16}
\buone\to \buonepm\to Z/2
\end{equation}
in which the rightmost arrow is ``forget everything but the grading.''
In fact this sequence can be split by associating to a locally
constant function $\epsilon:X\to\Z/2$ the ``trivial'' graded
$\T$-bundle
\[
\mathbf{1}^{\epsilon}:=\left(
X\times \T,\epsilon\right).
\]
This splitting is compatible with the monoidal structure, but not with
its symmetry.  It gives an equivalence
\begin{equation}\label{eq:17}
\buonepm\approx \buone\times\Z/2
\end{equation}
of monoidal categories (but not of {\em symmetric} monoidal
categories).  It follows that the group $H^{1}(X;\tpm)$, of
isomorphism classes of graded $\T$-bundles over $X$, is isomorphic to
$H^{0}(X;\Z/2)\times H^{1}(X;\T)$.  Since $X$ is assumed to be
paracompact, this in turn is isomorphic to $H^{0}(X;\Z/2)\times
H^{2}(X;\Z)$.  The automorphism group of any graded $\T$-bundle is the
group of continuous maps from $X$ to $\T$.

\subsection{Graded central extensions} 
\label{sec:grad-centr-extens}

Building on the notion of graded $\T$-bundles we now turn to graded
central extension of groupoids.  The reader is referred to
Appendix~\ref{sec:groupoids} for our conventions on groupoids, and a
recollection of the fundamental notions.  Unless otherwise stated all
groupoids will be assumed to be {\em local quotient groupoids}
(\S\ref{sec:local-glob-quot}), in the sense that they admit a
countable open cover by sub-groupoids, each of which is weakly
equivalent to a compact Lie group acting on a Hausdorff space.

Let $X=(X_{0},X_{1})$ be a groupoid.  Write $\bzt$ for the groupoid
associated to the action of $\Z/2$ on a point.

\begin{defin}\label{def:2}
A {\em graded groupoid} is a groupoid $X$ equipped with
a functor $\epsilon:X\to \bzt$.  The map
$\epsilon$ is called the {\em grading}.
\end{defin}

The collection of gradings on $X$  forms a groupoid, in which a
morphism is a natural transformation.  Spelled out, a grading of $X$ is a
function $\epsilon:X_{1}\to\Z/2$ satisfying $\epsilon(g\circ
f)=\epsilon(g)+\epsilon(f)$, and a morphism from $\epsilon_{0}$ to
$\epsilon_{1}$ is a continuous function $\eta:X_{0}\to \Z/2$
satisfying, for each $(f:x\to y)\in X_{1}$,
\[
\epsilon_{1}(f) = \epsilon_{0}(f)+(\eta(y)-\eta(x)).
\]

\begin{eg}
\label{eg:6}
Suppose that $X=S\modmod{G}$, with $S$ a connected topological space.
Then a grading of $X$ is just a homomorphism $G\to \Z/2$, making $G$
into a graded group.  
\end{eg}

We denote the groupoid of gradings of $X$
\[
\shom(X,\bzt).
\]

\begin{defin}\label{def:3}
A {\em graded central extension} of $X$ is a graded $\T$-bundle $L$
over $X_{1}$, together with an isomorphism of graded $\T$-bundles on
$X_{2}$
\[
\lambda_{g,f}:L_{g}\otimes L_{f}\to L_{g\circ f}
\]
satisfying the cocycle condition, that the diagram
\[
\xymatrix{
\left(L_{h}L_{g} \right)L_{f}\ar[dr] \ar[r]^{\approx} &L_{h}\left(L_{g}L_{f}
\right) \ar[r] &  L_{h}L_{g\circ f} \ar[d]\\
& L_{h\circ g}L_{f}\ar[r] & L_{h\circ g\circ f}
}
\]
of graded $\T$-bundles on $X_{3}$ commutes.
\end{defin}

If $L\to X_{1}$ is a graded central extension of $X$, then the pair
$\tilde X=(X_{0},L)$ is a graded groupoid over $X$, and the functor
$\tilde X\to X$ represents $\tilde X $ as a graded central extension
of $X$ by $\T$ in the evident sense.  Our terminology comes from this
point of view.  Some constructions are simpler to describe in terms of
the graded $\T$-bundles $L$ and others in terms of $\tilde X\to X$.

The collection of graded central extensions of $X$ forms a symmetric
monoidal $2$-category which we denote $\exttwocat(X,\tpm)=\sext{X}$.
The category of morphisms in $\sext{X}$ from $L^{1}\to L^{2}$ is
the groupoid of graded $\T$-torsors $(\eta,\epsilon)$ over $X_{0}$,
equipped with an isomorphism
\[
\eta_{b}L^{1}_{f}\to L^{2}_{f}\eta_{a}
\]
making
\[
\xymatrix{
\eta_{c}L^{1}_{g}L^{1}_{f}\ar[r]\ar[d] & L^{2}_{g}\eta_{b}L^{1}_{f}\ar[r] &
L^{2}_{g}L^{2}_{f}\eta_{a} \ar[d]  \\
\eta_{c}L^{1}_{g\circ f}\ar[rr] && L^{2}_{g\circ f}\eta_{a}
}
\]
commute.  
The {\em tensor product} $L\otimes L'$ is the graded central extension 
\[
L\otimes L' \to X_{1}
\]
with structure map
\[
L_{g}\otimes L'_{g}\otimes L_{f}\otimes
L'_{f}\xrightarrow{1\otimes\text{T}\otimes 1}
L_{g}\otimes L_{f}\otimes L'_{g}\otimes
L'_{f}\xrightarrow{\lambda_{g,f}\otimes\lambda'_{g,f}}
L_{g\circ f}\otimes L'_{g\circ f}.
\]
The {\em symmetry isomorphism} $L\otimes L'\to L'\otimes L$ is derived
from the symmetry of the tensor product of graded $\T$-bundles.

For the purposes of twisted $K$-groups it suffices to work with the
$1$-category quotient of $\sext{X}$.

\begin{defin}
\label{def:18}
The category $\extonecat(X;\tpm)=\hext{X}$ is the category with objects the graded
central extensions of $X$, and with morphisms from $L$ to $L'$ the set
of isomorphism classes in $\sext{X}(L,L').$
\end{defin}

The symmetric monoidal structure on $\sext{X}$ makes $\hext{X}$ into a
symmetric monoidal groupoid in the evident way.

\begin{rem}\label{rem:9}
A $1$-automorphisms of $L$ consists of a graded $\T$ torsor $\eta$
over $X_{0}$, together with an isomorphism $\eta_{a}\to \eta_{b}$ over
$X_{1}$, satisfying the cocycle condition.  In this way the category
of automorphisms of any twisting can be identified with the groupoid
of graded line bundles over $X$.  A graded line bundle on $X$ defines
an element of $K^{0}(X)$ (even line bundles go to line bundles, and
odd line bundles go to their negatives).  The fundamental property
relating twistings and twisted $K$-theory is that the automorphism
$\eta$ acts on twisted $K$-theory as multiplication by the
corresponding element of $K^{0}(X)$ (Proposition~\ref{thm:8}).
\end{rem}

\begin{rem}
\label{rem:6}
\begin{textList}
\item The formation of $\sext{X}$ is functorial in $X$, in the
sense of $2$-categories.  If $F:Y\to X$ is a map of groupoids, and
$L\to X_{1}$ is a graded central extension of $X$, then $F^{\ast}L\to
Y_{1}$ gives a graded central extension, $F^{\ast}L$ of $Y$.  If
$\eta\to X_{0}$ is a morphism from $L^{1}$ to $L^{2}$, then
$F^{\ast}\eta$ defines a morphism from $F^{\ast}L^{1}$ to
$F^{\ast}L^{2}$.  
\item If $T:Y_{0}\to X_{1}$ is a natural transformation from $F$ to $G$, then
the graded line bundle $T^{\ast}L$ determines a morphism from
$F^{\ast}L$ to $G^{\ast}L$.  This is functorial in the sense that
given $\eta:L^{1}\to L^{2}$, there is a $2$-morphism relating the two
ways of going around
\[
\xymatrix@dr{
F^{\ast}L^{1}\ar[r]^{F^{\ast}\eta} \ar[d]_{T^{\ast}L^{1}} &
F^{\ast}L^{2}\ar[d]^{T^{\ast}L^{2}}\\
G^{\ast}L^{1}\ar[r]_{G^{\ast}\eta}\ar@{}[ur]|\implies  &  G^{\ast}L^{2}
}.
\]
The $2$-morphism is the isomorphism 
\[
(G^{\ast}\eta)\circ(T^{\ast}L^{1})\to (T^{\ast}L^{2})\circ(F^{\ast}\eta)
\]
gotten by pulling back map $\eta:L^{1}\to L^{2}$ along $T$.  It is
given pointwise over $y\in Y$, $Ty:Fy\to Gy$ by
\[
(\eta_{Gy}) (L^{1}_{Ty}) \to (L^{2}_{Ty})(\eta_{Fy}). 
\]

\item It follows that the formation of $\hext{X}$ is functorial in
$X$, making $\hext{X}$ a (weak) presheaf of groupoids.  
\end{textList}
\end{rem}

\begin{eg}\label{eg:8}
Suppose that $G$ is a group, and $X=\point\modmod{G}$.  Then a graded
central extension of $X$ is just a graded central extension of $G$ by
$\T$.
\end{eg}

Forgetting the $\T$-bundle gives a functor from $\sext{X}$ to the
groupoid of gradings of $X$, and the decomposition~\eqref{eq:17} gives
a $2$-category equivalence
\begin{equation}\label{eq:46}
\exttwocat(X,\tpm)\approx
\exttwocat(X,\T)\times \shom(X,\bzt)
\end{equation}
which is not, in general, compatible with the monoidal structure.
Here $\exttwocat(X,\T)$ is the $2$-category of evenly graded (ie, ordinary)
central extensions of $X$ by $\T$.

\subsubsection{Classification of graded central extensions}
\label{sec:class-grad-centr}

We now turn to the classification of graded central extensions of a
groupoid $X$.  In view of~\eqref{eq:46}, it suffices to separately
classify $\T$-central extensions (graded central extensions which are
purely even) and gradings.

For the $\T$-central extensions, first recall that the category of
$\T$-torsors on a space $Y$ is equivalent to the category whose
objects are $\T$-valued Cech $1$-cocycles, $\Check{Z}^{1}(Y;\T)$,
and in which a morphism from $z_{0}$ to $z_{1}$ is a Cech $0$-cochain
$c\in \Check{C}^{0}(Y;\T)$ with the property that 
\[
\delta c = z_{1}-z_{0}.
\]
Now consider the double complex for computing the Cech
hyper-cohomology groups of the nerve $X_{\bullet}$, with coefficients
in $\T$:
\begin{equation}\label{eq:3}
\xymatrix@=4ex{
&&& \\
\check{C}^{2}(X_{0};\T) \ar[r]\ar[u] &\check{C}^{2}(X_{1};\T)\ar[r]\ar[u] & 
\check{C}^{2}(X_{2};\T)\ar[r]\ar[u]  & \\
\check{C}^{1}(X_{0};\T) \ar[r]\ar[u] &\check{C}^{1}(X_{1};\T)\ar[r]\ar[u] &
\check{C}^{1}(X_{2};\T)\ar[r]\ar[u]  & \\
\check{C}^{0}(X_{0};\T) \ar[r]\ar[u] &\check{C}^{0}(X_{1};\T)\ar[r]\ar[u] &
\check{C}^{0}(X_{2};\T)\ar[r]\ar[u]  & \\
}
\end{equation}
In terms of the Cech cocycle model for $\T$-bundles, the $2$-category
of $\T$-central extension of $X$ is equivalent to the category whose
objects are cocycles in~\eqref{eq:3}, of total degree $2$, whose
component in $\Check{C}^{2}(X_{0};\T)$ is zero.  The $1$-morphisms
are given by cochains of total degree $1$, whose coboundary has the
property its component in $\Check{C}^{2}(X_{0};\T)$ vanishes.  The
$2$-morphisms are given by cochains of total degree $0$.  Write
\[
\Check{H}^{\ast}(X)\quad\text{and}\quad\Check{H}^{\ast}(X_{0})
\]
for the Cech hyper-cohomology of $X$, and the Cech cohomology of
$X_{0}$ respectively.  Then
the group of isomorphism classes of {\em even} graded $\T$-gerbes is
given by the kernel of the map
\[
\Check{H}^{2}(X;\T)\to 
\Check{H}^{2}(X_{0};\T),
\]
the group of isomorphism classes of $1$-automorphisms of any even
graded $\T$-gerbe is $\Check{H}^{1}(X;\T)$, and the group of
$2$-automorphisms of any $1$-morphism is $\Check{H}^{0}(X;\T)$.

As for the gradings, the group of isomorphism classes of gradings is 
\[
\ker\left\{\Check{H}^{1}(X;\Z/2)\to \Check{H}^{1}(X_{0};\Z/2) \right\}
\]
and the automorphism group of any grading is 
\[
\Check{H}^{0}(X;\Z/2).
\]

For convenience, write
\[
\Check{H}_{\rel}^{t}(X;A) = \ker\left\{\Check{H}^{t}(X;A)\to
\Check{H}^{t}(X_{0};A) \right\}.
\]

\begin{prop}\label{thm:35}
The group $\pi_{0}\sext{X}$ of isomorphism classes of graded central
extension of $X$ is given by the set-theoretically split extension
\begin{equation}\label{eq:26}
\Check{H}_{\rel}^{2}(X;\T)\to \pi_{0}\sext{X} \to
\Check{H}_{\rel}^{1}(X;\Z/2)
\end{equation}
with cocycle 
\[
c(\epsilon,\mu)=\beta(\epsilon\cup\mu),
\]
where $\beta:\Z/2=\{\pm1 \}\subset \T$ is the inclusion.
The group of isomorphism classes of automorphisms of any graded
central extension of $X$ is $\Check{H}^{1}(X;\T)\times
\Check{H}^{0}(X;\Z/2)$, and the group of $2$-automorphisms of any
morphism of graded central extensions is $\Check{H}^{0}(X;\T)$.
\end{prop}

\begin{pf}
Most of this result was proved in the discussion leading up to its
statement.   The decomposition~\eqref{eq:46} gives the exact sequence,
as well as a set-theoretic splitting 
\[
s:\Check{H}_{\rel}^{1}(X;\Z/2)\to \pi_{0}\sext{X}.
\]
It remains to identify the cocycle describing the group structure.
Suppose that 
\[
\epsilon,\mu:X\to\bzt
\]
are two gradings of $X$.  Then by the discussion leading up
to~\eqref{eq:17}, $s(\epsilon)s(\mu)$ is the graded central extension
given by
\[
\left(s(\epsilon)s(\mu) \right)_{f} = 
\mathbf{1}^{\epsilon(f)}\mathbf{1}^{\mu(f)}
\approx \mathbf{1}^{\epsilon(f)+\mu(f)},
\]
and structure map
\begin{equation}\label{eq:27}
\mathbf{1}^{\epsilon(g)}\mathbf{1}^{\mu(g)}
\mathbf{1}^{\epsilon(f)}\mathbf{1}^{\mu(f)}
\to
\mathbf{1}^{\epsilon(g)}\mathbf{1}^{\epsilon(f)}
\mathbf{1}^{\mu(g)}\mathbf{1}^{\mu(f)}
\to 
\mathbf{1}^{\epsilon(g\circ f)}\mathbf{1}^{\mu(g\circ f)}.
\end{equation}
Using the canonical identifications
\[
\mathbf{1}^{a}\mathbf{1}^{b}
=\mathbf{1}^{a+b},
\]
and
\begin{align*}
\epsilon(g\circ f) &= \epsilon(g)+\epsilon(f)\\
\mu(g\circ f) &= \mu(g)+\mu(f)
\end{align*}
one checks that~\eqref{eq:27} can be identified with the automorphism
\[
(-1)^{\epsilon(f)\mu(g)}
\]
of the trivialized graded line
\[
\mathbf{1}^{\epsilon(g)+\mu(g)+\epsilon(f)+\mu(f)}\approx
\mathbf{1}^{\epsilon(g\circ f)+\mu(g\circ f)}.
\]
Similarly, the structure map of $s(\epsilon+\mu)$ can be identified
with the {\em identity map} of the same trivialized graded line.
It follows that $s(\epsilon)s(\mu)= c(\epsilon,\mu)s(\epsilon+\mu)$,
where $c(\epsilon,\mu)$ is graded central extension with 
with $L_{f}=\mathbf{1}$, and 
\[
\lambda_{g,f}=(-1)^{\epsilon(f)\mu(g)}.
\]
Now the $2$-cocycle $\epsilon(f)\mu(g)$ is precisely the
Alexander-Whitney formula for the cup product $\epsilon\cup\mu\in
Z^{2}(X;\Z/2)$.  This completes the proof.
\end{pf}

One easy, but very useful consequence of Proposition~\ref{thm:35} is
the stacky nature of the morphism categories in $\sext{X}$.

\begin{cor}\label{thm:38}
Let $P,Q\in\sext{X}$, and $f:Y\to X$ be a local equivalence.
Then the functor
\[
f^{\ast}:\sext{X}(P,Q)\to \sext{Y}(f^{\ast}P,f^{\ast}Q)
\]
is an equivalence of categories, and so
\[
f^{\ast}:\hext{X}(P,Q)\to \hext{Y}(f^{\ast}P,f^{\ast}Q)
\]
is a bijection of sets. \qed
\end{cor}

Another useful, though somewhat technical consequence of the
classification is

\begin{cor}
\label{thm:36}
Suppose that $X$ is a groupoid with the property that the maps 
\begin{align*}
\Check{H}^{2}(X;\T) &\to \Check{H}^{2}(X_{0};\T)  \\
\Check{H}^{1}(X;\Z/2) &\to \Check{H}^{1}(X_{0};\Z/2)
\end{align*}
are zero.   If $Y\to X$ is a local equivalence, then the maps
\begin{align*}
\Check{H}^{2}(Y;\T) &\to \Check{H}^{2}(Y_{0};\T)  \\
\Check{H}^{1}(Y;\Z/2) &\to \Check{H}^{1}(Y_{0};\Z/2)
\end{align*}
are zero, and 
\[
\sext{X} \to \sext{Y},\text{ and }
\hext X \to \hext Y
\]
are equivalences.
\end{cor}

\begin{pf}
The first assertion is a simple diagram chase.   It has the
consequence that the maps
\begin{align*}
\Check{H}_{\rel}^{2}(X;\T) &\to \Check{H}^{2}(X_{0};\T) \\
\Check{H}_{\rel}^{1}(X;\Z/2) &\to \Check{H}^{1}(X_{0};\Z/2) \\
\Check{H}_{\rel}^{2}(Y;\T) &\to \Check{H}^{2}(Y_{0};\T) \\
\Check{H}_{\rel}^{1}(Y;\Z/2) &\to \Check{H}^{1}(Y_{0};\Z/2)
\end{align*}
are isomorphisms, and the second assertion follows.
\end{pf}

For later reference we note the following additional consequence of
Proposition~\ref{thm:35}
\begin{cor}
\label{thm:30}
Suppose that $X$ is a local quotient groupoid, and that $\tau$ is a
twisting of $X$ represented by a local equivalence $P\to X$ and a
graded central extension $\tilde P\to P$.  Then $\tilde P$ is a local
quotient groupoid.
\end{cor}

\begin{pf}
Since the property of being a local quotient groupoid is an invariant
of local equivalence,  we know that $P$ is a local quotient
groupoid.    The question is also local in $P$, so we may assume $P$
is of the form $S\modmod G$ for a compact Lie group $G$.  By our
assumptions, the action of $G$ on $S$ has locally contractible
slices.   Working still more locally in $S$ we may assume $S$ is
contractible.  But then Proposition~\ref{thm:35} implies that $\tau$ is
given by a (graded) central extension of $\tilde G\to G$, and $\tilde
P=S\modmod{\tilde G}$.  
\end{pf}

\subsection{Twistings}
\label{sec:twistings-1}

We now describe the category $\twist{X}$ of $K$-theory twistings on a
local quotient groupoid $X$ (\S\ref{sec:local-glob-quot}).  The
objects of $\twist{X}$ are pairs $a=(P,L)$ consisting of a local
equivalence $P\to X$, and a graded central extension $L$ of $P$.  The
set of morphisms from $a=(P_{0}, L_{0})$ to $b=(P_{1},L_{1})$ is
defined to be the colimit
\[
\twist X(a,b)=\varinjlim_{p:P\to P_{12}}
\hext{P}(p^{\ast}\pi_{1}^{\ast}a,p^{\ast}\pi_{2}^{\ast}b), 
\]
where $P_{12}=P_{1}\times_{X}P_{2}$, the limit is taken over
$\ocov{P_{12}}$ and our notation refers to the diagram
\[
\xymatrix@=3ex{
& P \ar[d]^{p} & \\
P_{1} & P_{12}\ar[r]_{\pi_{2}}\ar[l]^{\pi_{1}} & P_{2}.
}
\]
We leave to the reader to check that $\hext{P}(a,b)$ does indeed define
a functor on the $1$-category quotient $\ocov{P_{12}}$.

The colimit appearing in the definition of $\twist X(a,b)$ is present
in order that the definition be independent of any extraneous choices.
In fact the colimit is attained at any stage.

\begin{lem}\label{thm:40}
For any local equivalence
\begin{equation}\label{eq:41}
p:P\to P_{12}
\end{equation}
the map
$\hext{P}(p^{\ast}\pi_{1}^{\ast}a,p^{\ast}\pi_{2}^{\ast}b)
\to\twist{X}(a,b)$ is an isomorphism.
\end{lem}

\begin{pf}
By Corollary~\ref{thm:38}, for any 
\[
P'\to P 
\] 
in $\ocov{P_{12}}$, the map 
\[
\hext{P}(a,b)\to \hext{P'}(a,b)
\]
is an isomorphism.  The result now follows from Corollary~\ref{thm:41}.
\end{pf}

For the composition law, suppose we are given three twistings
\[
a=(P_{1},L_{1}), b=(P_{2},L_{2}),\text{ and } c=(P_{3},L_{3}).
\]
Find a $P_{123}\in\cov X$ and maps
\[
\xymatrix{
& P_{123} \ar[dl]_{p_{1}}\ar[d]^{p_{2}}\ar[dr]^{p_{3}} \\
P_{1} & P_{2} & P_{3}
}
\]
(for example one could take $P_{123}$ to be the ($2$-category) fiber
product $P_{1}\times_{X}P_{2}\times_{X}P_{3}$, and $p_{i}$ projection
to the $i^{\text{th}}$ factor).  By Lemma~\ref{thm:40}, the maps
\begin{align*}
\hext{P_{123}}(p_{1}^{\ast}a,p_{2}^{\ast}b) &\to \twist{X}(a,b) \\
\hext{P_{123}}(p_{1}^{\ast}a,p_{2}^{\ast}c) &\to \twist{X}(a,c) \\
\hext{P_{123}}(p_{1}^{\ast}b,p_{2}^{\ast}c) &\to \twist{X}(b,c) 
\end{align*}
are bijections.  With these identifications, the composition law in $\twist{X}$
is formed from that in $\hext{P_{123}}$.  We leave to the reader to
check that this is well-defined.

The formation of $\twist{X}$ is functorial in $X.$  Given $f:Y\to X$,
and $a=(P,L)\in\twist X$ form
\[
\begin{CD}
Y\times_{X}P @>\pi>> P\\
@VVV @VVV \\
Y @>>f> X
\end{CD}
\]
and set
\[
f^{\ast}a=(f^{\ast}P,\pi^{\ast}L).
\]

\begin{prop}\label{thm:43}
The association $X\mapsto \twist{X}$ is a weak presheaf of groupoids.
If $Y\to X$ is a local equivalence then $\twist{Y}\to \twist{X}$ is an
equivalence of categories. \qed
\end{prop}

There is an evident functor
\[
\hext{X}\to \twist{X}.
\] 
\begin{prop}\label{thm:3}
When $X$ satisfies the condition of Corollary~\ref{thm:36}, the
functor
\begin{equation}\label{eq:42}
\hext{X}\to \twist{X}
\end{equation}
is an equivalence of categories.
\end{prop}

\begin{pf}
Lemma~\ref{thm:40} shows that~\eqref{eq:42} is fully faithful.
Essential surjectivity is a consequence of Corollary~\ref{thm:36}.
\end{pf}

\begin{cor}\label{thm:42}
If $Y\to X$ is a local equivalence, and $Y$ satisfies the conditions
of Corollary~\ref{thm:36} then the functors
\[
\twist{X}\to\twist{Y}\leftarrow \hext Y
\]
are equivalences of groupoids. \qed
\end{cor}

Combining this with Proposition~\ref{thm:35} gives
\begin{cor}\label{thm:2}
The group $\pi_{0}\twist X$ of isomorphism classes of twistings on $X$ is
the set-theoretically split extension
\[
\Check{H}^{2}(X;\T)\to \pi_{0}\twist X\to
\Check{H}^{1}(X;\Z/2)
\]
with cocycle 
\[
c(\epsilon,\mu)=\beta(\epsilon\cup\mu).
\]
The group of automorphisms of any twisting is
$\Check{H}^{1}(X;\T)\times \Check{H}^{0}(X;\Z/2)$. \qed
\end{cor}

We now switch to the point of view of ``fibered categories'' in order
to more easily describe the functorial properties of twisted $K$-groups.

Let $\frakext$ denote the category whose objects are pairs $(X,L)$
consisting of a groupoid $X$ and a graded central extension $L$ of
$X$.  A morphism $(X,L)\to (Y,M)$ is a functor $f:X\to Y$, and an
isomorphism $L\to f^{\ast}M$ in $\hext{X}$.  We identify
morphisms
\[
f,g:(X,L)\to (Y,M)
\]
if there is a natural transformation $T:f\to g$ for which the
following diagram commutes:
\[
\xymatrix@R=8pt{
& f^{\ast}M \ar[dd]^{\eta_{T}}\\
L \ar[ur]\ar[dr] & \\
& g^{\ast}M
}
\]
The functor $(X,L)\mapsto X$ from $\frakext$ to groupoids represents
$\frakext$ as a fibered category, fibered over the category of
groupoids.

Similarly, we define a category $\twistcat$ with objects $(X,a)$
consisting of a groupoid $X$, and a $K$-theory twisting $a\in\twist
X$, and morphisms $(X,a)\to (Y,b)$ to be equivalence classes of pairs
consisting of a functor $f:X\to Y$ and an isomorphism $a\to f^{\ast}b$
in $\twist{X}$.

There is an inclusion $\frakext\to\twistcat$ corresponding to the
inclusion $\hext{X}\to \twist{X}$.   Corollary~\ref{thm:42}
immediately implies

\begin{lem}\label{thm:44}
Suppose that $F:\frakext\to \mathcal C$ is a functor 
sending every morphism $(f,t):(X,L)\to (Y,M)$ in which $f$ is a local
equivalence to an isomorphism.  Then there is a factorization 
\[
\xymatrix{
\frakext \ar[d] \ar[r]^{F} & {}\mathcal{C}\\
\twistcat \ar @{-->} [ur]_{F'}   &
}
\]
Moreover any two such factorizations are naturally isomorphic by a
unique natural isomorphism.
\end{lem}

\subsection{Examples of twistings}
\label{sec:examples-twistings}

\begin{eg}\label{eg:12}
Suppose that $X$ is a space, $P\to X$ a principal $G$ bundle, and
\begin{align*}
\tilde G &\to G \\
G &\xrightarrow{\epsilon}{}\Z/2
\end{align*}
is a graded central extension of $G$.  Then $P\modmod G\to X$ is a
local equivalence, $P\modmod{\tilde G}$ is a graded central
extension of $P\modmod G$, and $(P\modmod{\tilde G},P\modmod G)$
represents a twisting of $X$.
\end{eg}

\begin{eg}\label{eg:3}
As a special case, we note that any double cover $P\to X$ defines a
twisting.  In this case $G=\Z/2$, $\tilde G\to G$ is the trivial
bundle, and $\epsilon:G\to \Z/2$ is the identity map.  Any cohomology
theory can be twisted by a double cover, and in fact these are the
only twistings of ordinary cohomology with integer coefficients.
\end{eg}

\begin{eg}
\label{eg:20} Suppose that $X=\point\modmod{G}$, with $G$ a compact
Lie group.  In this case every local equivalence $\tilde X\to X$
admits a section.  Moreover the inclusion
\[
\{\id_{X}\}\to \cov{X}
\]
of the trivial category consisting of the identity map of $X$ into
$\cov{X}$
is an equivalence.  It follows that
\[
\exttwocat(\point\modmod{G},\tpm)\to \twist{X}
\]
is an equivalence of categories, and so twistings of $X$ in this case
are just graded central extensions of $G$.  Using
Corollary~\ref{thm:2}, one can draw the same conclusion for
$S\modmod{G}$ when $S$ is contractible.
\end{eg}

We now describe the main example of twistings used in this paper.

\begin{eg}\label{eg:13}
Suppose that $G$ is a connected compact Lie group, and consider the
path-loop fibration
\begin{equation}\label{eq:2}
\Omega G\to PG\to G.
\end{equation}
We regard $PG$ as a principal bundle over $G$ with structure group
$\Omega G$.    The group $G$ acts on everything by conjugation.  
Write $LG$ for the group $\Map(S^{1},G)$ of smooth maps from
$S^{1}$ to G.   The homomorphism ``evaluation at 1'' $:LG\to G$ is split by the
inclusion of the constant loops.  This exhibits $LG$ as a semidirect
product 
\[
LG \approx \Omega G\rtimes G.
\]
The group $LG$ acts on the fibration~\eqref{eq:2} by conjugation.
The action of $LG$ on $G$ factors through the action of $G$ on itself by
conjugation, through the map $LG\to G$.  This defines a map of
groupoids
\[
PG\modmod LG \to G\modmod{G}
\]
which is easily checked to be a local equivalence.  A graded
central extension $\tilde LG\to LG$ then defines a twisting of $G/\!/G$.
\end{eg}

We will write $\tau$ for a typical twisting of $X$, and write the
monoidal structure  additively: $\tau_{1}+\tau_{2}$.  We will use the
notations $(\tilde P^{\tau},P^{\tau})$ and $(L^{\tau}, P^{\tau})$ for
typical representing graded central extensions.  This is consistent
with writing the monoidal structure additively:
\[
L^{\tau_{1}+\tau_{2}}=L^{\tau_{1}}\otimes L^{\tau_{2}}.
\]

\begin{eg}\label{eg:21}
Suppose that $Y=S\modmod G$ and $H\subset G$ is a normal subgroup.
Write $X=S\modmod H$, and $f:X\to Y$ for the natural map.  If $\tau$
is a twisting of $Y$ then $f^{\ast}\tau$ has a natural action of $G/H$
(in the $2$-category sense), and the map $X\to Y$ is invariant under
this action (again, in the $2$-category sense).  To see this it is
easiest to replace $X$ by the weakly equivalent groupoid
\[
X' = (S\times G/H) \modmod G, 
\]
and factor $X\to Y$ as 
\[
X\xrightarrow{i} X'\xrightarrow{f'} Y.
\]
Since $i^{\ast}:\twist{X'}\to\twist{X}$ is an equivalence of
$2$-categories, it suffices to exhibit an action of $G/H$ on
${f'}^{\ast}\tau$.   The obvious left action of $G/H$ on $S\times G/H$
commutes with the right action of $G$, giving an action of $G/H$ on
$X'$ commuting with $f'$.    The action of $G/H$ on ${f'}^{\ast}\tau$
is then a consequence of naturality.
\end{eg}

\begin{eg}\label{eg:22}
By way of illustration, consider the situation of Example~\ref{eg:21}
in which $H$ is commutative, $S=\{\point \}$, and $\tau$
is given by a central extension
\[
\T\mapsto \tilde G\to G.
\]
Write
\[
\T\mapsto \tilde H\to H,
\]
for the restriction of $\tau$ to $H$ and assume, in addition, that
$\tilde H$ is commutative.  Then the action of $G/H$ on $f^{\ast}\tau$
constructed in Example~\ref{eg:21} works out to be the natural action
of $G/H$ on $\tilde H$ given by conjugation.
\end{eg}

\section{Twisted K-groups}\label{sec:twisted-k-groups}
\subsection{Axioms}\label{sec:axioms}

Before turning to the definition of twisted $K$-theory, we list some
general properties describing it as a cohomology theory on the
category $\twistcat$ of local quotient groupoids equipped with a
twisting.  These properties almost uniquely determine twisted
$K$-theory, and suffice to make our main computation in
Section~\ref{sec:comp-k_gg-twist}.

Twisted $K$-theory is going to be homotopy invariant, so we need
to define the notion of homotopy

\begin{defin}\label{def:1}
A {\em homotopy} between two maps 
\[
f,g:(X,\tau_{X})\to
(Y,\tau_{Y})
\]
is a map
\[
(X\times[0,1], \pi^{\ast}\tau_{X})\to (Y,\tau_{Y})
\]
($\pi:X\times [0,1]\to X$ is the projection) whose restriction to
$X\times\{0 \}$ is $f$, and to $X\times\{1 \}$ is $g$.
\end{defin}

Twisted $K$-theory is also a cohomology theory.  To state this
properly involves defining the relative twisted $K$ theory of a triple
$(X,A,\tau)$ consisting of a local quotient groupoid $X$, a
sub-groupoid $A$, and a twisting of $X$.  We form a category of the
triples $(X,A,\tau)$ in the same way we formed $\twistcat$.  We'll
call this the {\em category of pairs in $\twistcat$}.

We now turn to the axiomatic properties of twisted $K$-theory.  

\begin{prop}\label{thm:7}
The association $(X, A,\tau)\mapsto K^{\tau+n}(X, A)$ to be constructed in
\S\ref{sec:definition-twisted-k} is a contravariant homotopy functor
on the category of pairs $(X, A,\tau)$ in $\twistcat$, taking local
equivalences to isomorphisms.
\end{prop}

\begin{prop}\label{thm:8}
The functors $K^{\tau+n}$ form a cohomology theory:
\begin{thmList}
\item\label{item:les} there is a natural long exact sequence
\begin{multline*}
\dots\to K^{\tau+n}(X,A)\to K^{\tau+n}(X)\to 
K^{\tau+n}(A) \\
\to K^{\tau+n+1}(X,A)\to K^{\tau+n+1}(X)\to K^{\tau+n+1}(A)
\to\cdots.\\
\end{multline*}

\item If $Z\subset A$ is a (full) subgroupoid whose closure is contained in
the interior of $A$, then the restriction (excision) map
\[
K^{\tau+n}(X,A)\to 
K^{\tau+n}(X\setminus Z,A\setminus Z)
\]
is an isomorphism.

\item If $(X,A,\tau)=\coprod_{\alpha} (X_{\alpha},A_{\alpha},\tau_{\alpha})$, then
the natural map
\[
K^{\tau+n}(X,A) \to \prod_{\alpha}K^{\tau_{\alpha}+n}(X_{\alpha},A_{\alpha})
\]
is an isomorphism.

\end{thmList}
\end{prop}

The combination of excision and the long exact sequence of a pair
gives the Mayer-Vietoris sequence
\begin{align*}
\dots &\to K^{\tau+n}(X)\to 
K^{\tau+n}(U)\oplus K^{\tau+n}(V)\to 
K^{\tau+n}(U\cap V) \\
&\to  K^{\tau+n+1}(X)\to \cdots
\end{align*}
when $X$ is written as the union of two sub-groupoids whose interiors
form a covering.   

\begin{prop}
\label{thm:19}
\begin{thmList}
\item\label{item:multiplication} There is a bilinear pairing
\[
K^{\tau+n}(X)\otimes K^{\mu+m}(X)\to
K^{\tau+\mu+n+m}(X)
\]
which is associative and (graded) commutative up to the natural
isomorphisms of twistings coming from Proposition~\ref{thm:7}.

\item\label{item:one-morphism} Suppose that $\eta:\tau\to\tau$ is a
$1$-morphism, corresponding to a graded line bundle $L$ on $X$.  Then
\[
\eta_{\ast}=\text{multiplication by $L$}:K^{\tau+n}(X)\to K^{\tau+n}(X)
\]
where $L$ is regarded as an element of $K^{0}(X)$ and the
multiplication is that of~\thmItemref{item:multiplication}
\end{thmList}
\end{prop}

Twisted $K$-theory also reduces to equivariant $K$-theory in special cases.
\begin{prop}
\label{thm:18} Let $X=S\modmod G$ be a global quotient groupoid, with
$G$ a compact Lie group, and $\tau$ a twisting given by a graded
central extension
\begin{gather*}
\T\to G^{\tau} \to G \\
\epsilon:G \to\Z/2.
\end{gather*}
\begin{thmList}
\item If $\epsilon=0$ then then $K^{\tau+n}(X)$ is the summand
\[
K_{G^{\tau}}^{n}(S)(1)
\subset 
K_{G^{\tau}}^{n}(S)
\]
on which $\T$ acts via its standard (defining) representation.   This
isomorphism is compatible with the product structure.
\item For general $\epsilon$, $K^{\tau+n}(X)$ is isomorphic to 
\[
K_{G^{\tau}}^{n+1}(S\times(\R(\epsilon), \R(\epsilon)\setminus \{0 \}))(1),
\]
in which the symbol $\R(\epsilon)$ denotes the $1$-dimensional
representation $(-1)^{\epsilon}$ of $G^{\tau}$.
\end{thmList}
\end{prop}

In part ii), When $\epsilon=0$, then $\R(\epsilon)$ is the trivial
representation, and the isomorphism can be composed with the
suspension isomorphism to give the isomorphism of i).  

When $\tau=0$, so that $G^{\tau}\approx G\times\T$, Proposition~\ref{thm:18}
reduces to an isomorphism
\[
K^{\tau+n}(X)\approx K_{G}^{n}(S).
\]
In view of this, we'll often write 
\[
K_{G}^{\tau+n}(S) = K^{\tau+n}(X)
\]
in case $X=S\modmod G$ and $\epsilon=0$.  Of course there is also a
relative version of Proposition~\ref{thm:18}.

The reader is referred to Section~4 of~\cite{math.AT/0206257} for a
more in depth discussion of the twistings of equivariant $K$-theory,
and interpretation of ``$\epsilon$'' part the twisting in terms of
graded representations.

Using the Mayer-Vietoris sequence one can easily check that result of
part i) of Proposition~\ref{thm:18} holds for {\em any} local quotient
groupoid $X$.  If the twisting $\tau$ is represented by a central
extension $P\to X$, then the restriction mapping is an isomorphism
\[
K^{\tau+n}(X)\approx K^{n}(P)(1).
\]
In this way, once $K$-theory is defined for groupoids, twisted
$K$-theory is also defined.

\subsection{Twisted Hilbert spaces}
\label{sec:twist-hilb-spac}

Our definition of twisted $K$-theory will be in terms of Fredholm
operators on a twisted bundle of Hilbert spaces.  In this section we
describe how one associates to a graded central extension of a
groupoid, a twisted notion of Hilbert bundle.  We refer the reader to
Appendix~\ref{sec:groupoids}, \S\ref{sec:hilb-space-bundl} for our
notation and conventions on bundles over groupoids, and to
\S\ref{sec:hilb-space-bundl} for a discussion of Hilbert space
bundles.

Let $X$ be a groupoid, and $\tau:\tilde X\to X$ a graded central
extension, whose associated graded $\T$-bundle we denote $L^{\tau}\to
X_{1}$.  As in Appendix~\ref{sec:fiber-bundles-over}, we will use
\[
a \xrightarrow{f}{}b \quad\text{and}\quad
a \xrightarrow{f}{}b\xrightarrow{g}{}c
\]
to refer to generic points of $X_{1}$ and $X_{2}$ respectively, and
so, for example, in a context describing bundles over $X_{2}$, the
symbol $H_{b}$ will refer to the pullback of $X$ along the map 
\begin{gather*}
X_{2} \to X_{0} \\
(a\to b\to c) \mapsto b.
\end{gather*}

\begin{defin}\label{def:7}
A $\tau$-twisted Hilbert bundle on $X$ consists of a 
$\Z/2$-graded Hilbert bundle $H$ on $X_{0}$, together with
an isomorphism (on $X_{1}$)
\[
L^{\tau}_{f}\otimes H_{a}\to H_{b}
\]
satisfying the cocycle condition that
\[
\renewcommand{\otimes}{}
\xymatrix{
L^{\tau}_{g}\otimes \left(L^{\tau}_{f}\otimes H_{a}\right) \ar@{<->}[rr]\ar[d]
&
&
\left(L^{\tau}_{g}\otimes L^{\tau}_{f}\right)\otimes H_{a}
\ar[d]
\\
L^{\tau}_{g}\otimes H_{b} \ar[r]
&
H_{c}
&
L^{\tau}_{g\circ f}\otimes H_{a}\ar[l]
}
\]
commutes on $X_{2}$.
\end{defin}

\begin{rem}
\label{rem:8}
Phrased differently, a twisted Hilbert bundle is just a graded
Hilbert bundle over $\tilde X$, with the property that the map
$H_{a}\to H_{b}$ induced by $(a\xrightarrow{f}{}b)\in \tilde X_{1}$
has degree $\epsilon(f)$, and for which the central $\T$ acts
according to its defining character.
\end{rem}

\begin{eg}\label{eg:17}
Suppose that $X$ is of the form $P\modmod{G}$, and that our twisting
corresponds to a central extension of $G^{\tau}\to G$.  Then a
projective unitary representation of $G$ (meaning a representation of
$G^{\tau}$ on which the central $\T$ acts according to its defining
character) defines a twisted Hilbert bundle over $X$.
\end{eg}

Suppose that $\tau$ and $\mu$ are graded central extensions of $X$,
with associated graded $\T$-torsors $L^{\tau}$ and $L^{\mu}$.  If $H$
is a $\tau$-twisted Hilbert bundle over $X$ and $W$ is a
$\mu$-twisted Hilbert bundle, then the graded tensor product
$H\otimes W$ is a $(\tau+\mu)$-twisted Hilbert bundle, with
structure map
\[
L_{f}^{\tau+\mu}\otimes H_{a}\otimes W_{a} =
L_{f}^{\tau}\otimes L_{f}^{\mu}\otimes H_{a}\otimes W_{a} \to
L_{f}^{\tau}\otimes H_{a}\otimes  L_{f}^{\mu}\otimes W_{a} \to
H_{b}\otimes W_{b}.
\]

Now suppose that $H^{1}$ and $H^{2}$ are $\tau$-twisted, graded
Hilbert bundles over $X$.
\begin{defin}
A {\em linear transformation} $T:H^{1}\to H^{2}$ consists of a
linear transformation of Hilbert bundles%
\ $T:H^{1}\to H^{2}$ on $X_{0}$ for which the following diagram commutes
on $X_{1}$:
\[
\begin{CD}
L_{f}H^{1}_{a} @>>> H^{1}_{b}\\
@V1\otimes TVV @VV TV \\
L_{f}H^{2}_{a} @>>> H^{2}_{b}.
\end{CD}
\]
\end{defin}

If $L^{\tau}$ is a graded central extension of $X$, we'll write
$\mathcal U^{\tau}_{X}$ (or just $\mathcal U^{\tau}$ if $X$ is
understood) for the category in which the objects are $\tau$-twisted
$Z/2$-graded Hilbert bundles, and with morphisms the linear isometric
{\em embeddings.}  If $f:Y\to X$ is a map, there is an evident functor
\[
f^{\ast}:\mathcal U^{\tau}_{X}\to \mathcal U^{f^{\ast}\tau}_{Y}.
\] 
A natural transformation $T:f\to g$ of functors $X\to Y$ gives a
natural transformation $T^{\ast}:f^{\ast}\to g^{\ast}$.  Using
Remark~\ref{rem:8} and descent, one easily checks that $f^{\ast}$ is
an equivalence of categories when $f$ is a local equivalence.

The category $\mathcal U^{\tau}_{X}$ is also functorial in $\tau$.
Indeed, suppose that $\eta:\tau\to\sigma$ is a morphism,
given by a graded $\T$-bundle $\eta$, and an isomorphism 
\[
\eta_{b}\otimes\tau_{f}\to \sigma_{f}\otimes\eta_{a}.
\]
If $H$ is a $\tau$-twisted Hilbert bundle, then $H\otimes\eta$
is a $\sigma$-twisted Hilbert bundle.  One easily checks that
$H\mapsto H\otimes\eta$ gives an equivalence of categories $\mathcal
U^{\tau}_{X}\to \mathcal U^{\sigma}_{X}$, with inverse $H\mapsto
H\otimes \eta^{-1}$.  The $2$-morphisms $\eta_{1}\to\eta_{2}$ give
natural isomorphisms of functors.

The tensor product of Hilbert spaces gives a natural tensor product
\[
\mathcal U^{\tau}_{X}\times \mathcal U^{\mu}_{X}\to \mathcal
U^{\tau+\mu}_{X}.
\]

\subsection{Universal Twisted Hilbert Bundles}

We now turn to the existence of special kinds of $\tau$-twisted
Hilbert bundles, following the discussion of
\S\ref{sec:hilb-space-bundl}.  We keep the notation of \S\ref{sec:twist-hilb-spac}.

\begin{defin}
\label{def:12} A $\tau$-twisted Hilbert bundle $H$ on $X$ is
\begin{thmList}\label{}
\item {\em universal} if for every $\tau$-twisted Hilbert space
bundle $V$ there is a unitary embedding $V\to H$;
\item {\em locally universal} if $H\vert_{U}$ is universal for every open
sub-groupoid $U\subset X$;
\item {\em absorbing} if for every $\tau$-twisted Hilbert space
bundle $V$ there is an isomorphism $H\oplus V\approx
H$; 
\item {\em locally absorbing} if $H\vert_{U}$ is absorbing for every open
sub-groupoid $U\subset X$.
\end{thmList}
\end{defin}

As in Appendix~\ref{sec:hilb-space-bundl}, if $H$ is (locally) universal, then
$H$ is automatically absorbing.

\begin{lem}
\label{thm:29} Suppose that $\tilde X\to X$ is a graded central
extension and $H$ is a graded Hilbert bundle on $\tilde X$.  Let
$H(1)\subset H$ be the eigenbundle on which the central $\T$ acts
according to its defining representation.  Then $H(1)$ is a
$\tau$-twisted Hilbert bundle on $X$ which is (locally) universal if
$H$ is.\qed
\end{lem}

\begin{lem}
\label{thm:20} If $X$ is a local quotient groupoid, and $\tau:\tilde
X\to X$ is a graded central extension then there exists a locally
universal $\tau$-twisted Hilbert bundle $H$ on $X$.  The bundle $H$ is
unique up to unitary equivalence.
\end{lem}

\begin{pf}
By Corollary~\ref{thm:30}, $\tilde X$ is a local quotient groupoid
which, by Corollary~\ref{thm:27}, admits a locally universal Hilbert
bundle.  The result now follows from Lemma~\ref{thm:29}.
\end{pf}

\subsection{Definition of twisted $K$-groups}
\label{sec:definition-twisted-k}

Our task is to define twisted $K$-groups for pairs $(X,A,\tau)$ in
$\twistcat$.  In view of Lemma~\ref{thm:44} it suffices to define
functors $K^{\tau+\ast}(X,A)$ for $(X,A,\tau)$ in $\frakext$, and show
that they take local equivalences to isomorphisms.  We will do this by
using spaces of Fredholm operators to construct a spectrum $\ul
K^{\tau}(X,A)$ and defining $K^{\tau+n}(X,A)=\pi_{-n}\ul
K^{\tau}(X,A)$.  The reader is referred to
\S\ref{sec:fredholm-operators} for some background discussion on
spaces of Fredholm operators.

Suppose then that $(X,\tau)$ is an object of $\frakext$ and $H$ is a
locally universal, $\tau$-twisted Hilbert bundle over $X$.  With
the notation of \S\ref{sec:twist-hilb-spac}, $H$ is given by a Hilbert
bundle $H$ over $X_{0}$, equipped with an isomorphism
\begin{equation}\label{eq:35}
L^{\tau}_{f}\otimes H_{a}\to H_{b}
\end{equation}
over $f:a\to b\in X_{1}$, satisfying the cocycle condition.   
The map 
\[
T\mapsto \id\otimes T
\]
is a homeomorphism between the spaces of Fredholm operators (See
\S\ref{sec:fredholm-operators}) $\fred{n}(H_{a})$ and
$\fred{n}(L^{\tau}_{f}\otimes H_{a})$ compatible with the structure
maps~\eqref{eq:35}.  The spaces $\fred{n}(H_{a})$ therefore form a
fiber bundle over $\fred{n}(H)$ over $X$.

We define spaces $\ul K^{\tau}(X)_{n}$ by
\[
\ul K^{\tau}(X)_{n} = 
\begin{cases}
\Gamma(X;\fred{0}(H)) &\qquad \text{$n$ even} \\
\Gamma(X;\fred{1}(H)) &\qquad \text{$n$ odd},
\end{cases}
\]
By an obvious modification of the arguments of
Atiyah-Singer~\cite{atiyah69:_index_fredh}, the results described in
\S\ref{sec:fredholm-operators} hold for the bundle $\fred{n}(H)$ over
$X$.  In particular, the maps~\eqref{eq:33} and the
homeomorphism~\eqref{eq:34} give weak homotopy equivalences
\[
\ul K^{\tau}(X)_{n} \to 
\Omega\ul K^{\tau}(X)_{n+1},
\]
making the collection of spaces 
\[
\ul K^{\tau}(X) = \{\ul K^{\tau}(X)_{n} \}
\]
into a spectrum. 

\begin{defin}
\label{def:14} Suppose that $(X,\tau)$ is a local quotient groupoid
equipped with a graded central extension $\tau$, and $H$ is a locally
universal, $\tau$-twisted Hilbert bundle over $X$.  The {\em
twisted $K$-theory spectrum} of $X$ is the spectrum $\ul K^{\tau}(X)$
defined above.
\end{defin}

To keep things simple, we do not indicate the choice of Hilbert bundle
$H$ in the notation $\ul K^{\tau}(X)$.  The value of the twisted $K$-group
is, in the end, independent of this choice.  See Remark~\ref{rem:12} below.

We now turn to the functorial properties of $X\mapsto \ul
K^{\tau}(X)$.  Suppose that $f:Y\to X$ is a map of local quotient
groupoids, and $\tau$ is a twisting of $X$.  Let $H_{X}$ be a
$\tau$-twisted, locally universal Hilbert bundle over $X$, and $H_{Y}$
an $f^{\ast}\tau$-twisted, locally universal Hilbert space bundle over
$Y$.  Since $H_{Y}$ is universal, there is a unitary embedding
$f^{\ast}H_{X}\subset H_{Y}$.  Pick one.  There is then an induced map
\begin{align*}
f^{\ast}\fred{n}(H_{X}) & \to
\fred{n}(H_{Y}) \\
T &\mapsto T\oplus\epsilon,
\end{align*}
($\epsilon$ is the base point) and so a map of spectra 
\[
f^{\ast}:\ul K^{\tau}(X)
\to \ul K^{\tau}(Y).
\] 

Suppose that $\eta:\sigma\to\tau$ is a morphism of central extensions
of $X$, given by a graded $T$-bundle $\eta$ over $X_{0}$, and
isomorphism
\[
\eta_{b}\otimes \sigma_{f}\to \tau_{f}\otimes \eta_{a}.
\]
If $H$ is a locally universal $\sigma$-twisted Hilbert bundle,
then $H\otimes\eta$ is a locally universal $\tau$-twisted Hilbert
bundle.  The map
\[
T\mapsto T\otimes \id_{\eta}
\]
then gives a homeomorphism $\fred{n}(H)\to \fred{n}(H\otimes \eta)$,
and so an isomorphism of spectra
\[
\eta_{\ast}:\ul K^{\sigma}(X) \to \ul K^{\tau}(X).
\]
Since automorphisms of $\eta$ commute with the identity map,
$2$-morphisms of twistings have no effect on $\eta_{\ast}$.  In this
way the association $\tau\mapsto K^{\tau}(X)$ can be made into a
functor on $\hext X$.

Now we come to an important point.  Suppose $Y\to X$ is the inclusion
of a (full) subgroupoid of a local quotient groupoid, and $H_{X}$ is locally
universal.  By Corollary~\ref{thm:47}, we may then take $H_{Y}$ to be
$f^{\ast}H_{X}$.  The bundle of spectra $\ul K^{\tau}(Y)$ is then just
the restriction of $\ul K^{\tau}(X)$.  This would not be true for
general groupoids and is the reason for our restriction to local
quotient groupoids.  We use this restriction property in the
definition of the twisted $K$-theory of a pair.  While this could be
avoided, the restriction property plays a key role in the proof of
excision, and does not appear to be easily avoided there.

\begin{defin}
Suppose that $A\subset X$ is a sub-groupoid of a local quotient
groupoid, and that $\tau$ is a graded central extension of $X$.  The
{\em twisted $K$-theory spectrum of $(X,A,\tau)$} is the homotopy
fiber $\ul K^{\tau}(X,A)$ of the restriction map $\ul K^{\tau}(X)\to
\ul K^{\tau}(A)$.
\end{defin}

If we write 
\[
\Gamma(X,A;\fred{n}(H))\subset
\Gamma(X;\fred{n}(H))
\]
for the subspace of sections whose restriction to $A$ is the basepoint
$\epsilon$, then
\[
\ul K^{\tau}(X,A)_{n} = \Gamma(N,A;\fred{n}(H)),
\]
where $N$ is the mapping cylinder of $A\subset X$
\[
N = X\amalg A\times [0,1]/\sim.
\]

\begin{defin}
The twisted $K$-group $K^{\tau+n}(X,A)$ is the group $\pi_{-n}\ul
K^{\tau}(X,A)=\pi_{0}\ul K^{\tau}(X,A)_{n}$.
\end{defin}

\begin{rem}\label{rem:12}
There are several unspecified choices that go into the definition of
the spectra $\ul K^{\tau}(X,A)$, and the induced maps between them as
$X$, $A$ and $\tau$ vary.  It follows from Propositions~\ref{thm:33}
and~\ref{thm:34} that these choices are parameterized by (weakly)
contractible spaces, and so have no effect on the homotopy invariants
(such as twisted $K$-groups, and maps of twisted $K$-groups) derived
from them.  
\end{rem}

\subsection{Verification of the axioms}

\subsubsection{Proof of Proposition~\ref{thm:7}: functoriality}

Most of this result was proved in the process of defining the groups
$K^{\tau+n}(X,A)$.  Functoriality in $\frakext$ follows from the
discussion of \S\ref{sec:definition-twisted-k} and
Remark~\ref{rem:12}.  For homotopy invariance, note that if $H$ is a
locally universal $\tau$-twisted Hilbert bundle over
$(X,A)$, then $\pi^{\ast}H$ is a locally universal Hilbert space
bundle over $(X,A)\times I$, and so
\[
\ul K_{n}^{\tau}((X,A)\times I)=\ul K_{n}^{\tau}(X,A)^{I},
\]
and the two restriction maps to $\ul K_{n}^{\tau}(X)$ correspond to
evaluation of paths at the two endpoints.  The two restriction maps
are thus homotopic, and homotopy invariance follows easily.

The assertion about local equivalences is an immediate consequence of
descent.  As remarked at the beginning of
\S\ref{sec:definition-twisted-k}, this, in turn, gives functoriality
on the category of pairs in $\twistcat$.

\subsubsection{Proof of Proposition~\ref{thm:8}: cohomological
properties}\label{sec:proof-prop-refthm:8}

The long exact sequence of a pair (assertion i)is just the long exact
sequence in homotopy groups associated to the fibration of spectra
\[
\ul K^{\tau}(X,A)\to
\ul K^{\tau}(X) \to 
\ul K^{\tau}(A), 
\]
The ``wedge axiom'' (part iii) is immediate from the definition.
More significant is excision (part ii).  In describing the proof, we
will freely use, in the context of groupoids, the basic constructions
of homotopy theory as described in~\ref{sec:point-set-topology}.
Write $U=X\setminus Z$ and let $N$ be the double mapping cylinder of 
\[
U \leftarrow  U\cap A \rightarrow A.
\]
Then the map $N\to X$ is a homotopy equivalence of groupoids, and so
by the diagram of fibrations
\[
\begin{CD}
\ul K^{\tau}(X,A) @>>> \ul K^{\tau}(X) @>>> \ul K^{\tau}(A) \\
@VVV @VVV @VVV \\
\ul K^{\tau}(N,A) @>>> \ul K^{\tau}(N) @>>> \ul K^{\tau}(A) \\
\end{CD}
\]
the map 
\[
\ul K^{\tau}(X,A)\to
\ul K^{\tau}(N,A)
\]
is a weak equivalence.  Similarly, if $N'$ denote the mapping
cylinder of $U\cap A\to U$, then 
\[
\ul K^{\tau}(X\setminus Z,A\setminus Z)=
\ul K^{\tau}(U,U\cap A)\to
\ul K^{\tau}(N',U\cap A)
\]
is a weak equivalence.  We therefore need to show that for each $n$,
the map
\[
\ul K^{\tau}(N,A)_{n}\to
\ul K^{\tau}(N',U\cap A)_{n}
\]
is a weak equivalence.  Let $H$ be a locally universal $\tau$-twisted
Hilbert bundle over $X$.  Then the pullback of $H$ to each of the
(local quotient) groupoids $N$, $N'$ $A$, $U$, $U'$, $U\cap A$ is also
locally universal.  It follows that the twisted $K$-theory spectra of
each of these groupoids is defined in terms of sections of the bundle
pulled back from $\fred{n}(H\otimes C_{1})$.  To simplify the notation
a little, let's denote all of these pulled back bundles $\fred{n}$.
Now consider the diagram
\[
\begin{CD}
\Gamma(N,A;\fred{n}) @>>> \Gamma(N',U\cap A;\fred{n})\\
@VVV @VVV \\
\ul K^{\tau}(N,A)_{n} @>>> \ul K^{\tau}(N',U\cap A)_{n}.
\end{CD}
\]
We are to show that the bottom row is a weak equivalence.  But the top row
is a homeomorphism, and the vertical arrows are
weak equivalences since the maps
\begin{align*}
(N\cup \cyl(A),A) &\to (N,A) \\
(N'\cup \cyl(U\cap A),U\cap A) &\to (N',U\cap A)
\end{align*}
are relative homotopy equivalences.

\subsubsection{Proof of Proposition~\ref{thm:19}: Multiplication}

The multiplication is derived from the pairing 
\begin{align*}
\fred{n}(H_{1})\times \fred{m}(H_{2}) &\to
\fred{n+m}(H_{1}\otimes H_{2}) \\
(S,T) &\mapsto S\ast T = S\otimes \id + \id\otimes T,
\end{align*}
the tensor structure on the category of twisted Hilbert space
bundles described in \S\ref{sec:twist-hilb-spac}, and the natural
identification of the $\Z/2$-graded tensor product $\cl(\R^{n})\otimes
\cl(\R^{m})\approx \cl(R^{n+m}$
Verification of part i) is left to the reader.

Even if one of $S$, $T$ is not acting on a locally universal Hilbert
bundle the product $S\ast T$ will.  This is particularly useful when
describing the product of an element of untwisted $K$-theory, with one
of twisted $K$-theory.  For example if $V$ is a vector bundle over
$X$, we can choose a Hermitian metric on $V$, regard $V$ as a bundle
of finite dimensional graded Hilbert spaces, with odd component $0$,
and take $S=0$.  Then $S\ast T$ is just the identity map of $V$
tensored with $T$.  More generally, a virtual difference $V-W$ of
$K^{0}(X)$ can be represented by the odd, skew-adjoint Fredholm
operator $S=0$ on the graded Hilbert space whose even part is $V$ and
whose odd part is $W$, and $S\ast T$ represents the product of $V-W$
with the class represented by $T$.  The assertion of Part ii) is the
special case in which $V$ is a graded line bundle.

\subsubsection{Proof of Proposition~\ref{thm:18}: Equivariant
$K$-theory}
\label{sec:proof-prop-refthm:18}

Let $X=S\modmod G$ be a global quotient, and $\tau$ a twisting given
by a graded central extension $G^{\tau}$ of $G$, and a homomorphism
$\epsilon:G\to\Z/2$.  Replacing $X$ with
$X\times(\R(\epsilon),\R(\epsilon)\setminus\{0 \})$ and
using~\eqref{eq:38}, if necessary, we may reduce to the case
$\epsilon=0$.  
Write  $V(1)$  for the summand of
\[
V=C_{1}\otimes L^{2}(G^{\tau})\otimes \ell^{2}
\]
on which the central $\T$ acts according to its defining character.
Then $H=S\times V(1)$ is a locally universal Hilbert bundle.
Our definition of $K^{\tau}(X)$ becomes
\[
K^{\tau-n}(X) = [S, \fred{n}(C_{n}\otimes V(1))]^{G}
\] 
which is the summand of 
\[
[S, \fred{n}(C_{n}\otimes V)]^{G^{\tau}}
\]
corresponding to the defining representation of $\T$.  So the result
follows from the fact that $\fred{n}(C_{n}\otimes V)$ is a classifying
space for $K_{G^{\tau}}^{-n}$.  While this is certainly well-known, we
were unable to find an explicit statement in the literature.  It
follows easily from the case in which $G$ is trivial.  Indeed, the
universal index bundle is classified by a map to any classifying space
for equivariant $K$-theory, and it suffices to show that this map is a
weak equivalence on the fixed point spaces for the closed subgroups
$H$ of $G$.  The assertion for the fixed point spaces easily reduces
to the main result of~\cite{atiyah69:_index_fredh}.

\subsection{The Thom isomorphism, pushforward, and the Pontryagin product}
\label{sec:thom-isomorphism}

We begin with a general discussion.  Let $E=\{E_{n}\xrightarrow{t_{n}}
\Omega E_{n+1}\}_{n=0}^{\infty}$ be a spectrum.  For a real vector space $V$,
equipped with a positive definite metric let $\Omega^{V}(E_{n})$
denote the space of maps from the unit ball $B(V)$ to $E_{n}$, sending
the unit sphere $S(V)$ to the base point.  The collection of spaces
$\Omega^{V}E_{n}$ forms a spectrum $\Omega^{V}E$.  An isomorphism
$V\approx \R^{k}$ gives an identification $\Omega^{V}E_{n}\approx
E_{n-k}$, and of $\Omega^{V}E$ with the spectrum derived from $E$ by
simply shifting the indices.  Such a spectrum is called a ``shift
desuspension'' of $E$ (see \cite{LMayS}).  Some careful organization
is required to avoid encountering signs by moving loop coordinates
past each other.  The reader is referred to~\cite{LMayS} for more
details.  Of course, for a space $X$ one has
\[
\left(\Omega^{V}E \right)^{n}(X)\approx
E^{n}(X\times(V,V\setminus\{0\})) \approx
E^{-k+n}(X).
\]

Now suppose that $V$ is a vector bundle of dimension $k$ over a space
$X$.  The construction described above can be formed fiberwise to form
a bundle
\[
\Omega^{V}E = \{\Omega^{V}E_{n} \}
\]
of spectra over $X$.  The group of vertical homotopy classes of
sections 
\begin{equation}\label{eq:36}
\pi_{0}\Gamma(X,\Omega^{V}E_{n})
\end{equation}
can then be thought of as a twisted form of $[X,E_{-k+n}]=E^{-k+n}(X)$.
We denote this twisted (generalized) cohomology group 
\[
E^{-\tau_{V}+n}(X).  
\]
Now the group~\eqref{eq:36} is the group of pointed homotopy classes
of maps $[X^{V},E_{n}]$ from the Thom complex of $V$ to $E_{n}$.  This
gives a tautological Thom isomorphism
\[
\tilde E^{n}(X^{V})=E^{n}(B(V),S(V))\approx \tilde E^{-\tau_{V}+n}(X).  
\]
The more usual Thom isomorphisms arise when a geometric construction
is used to trivialize the bundle $\Omega^{V}E$.  Such a trivialization
is usually called an ``$E$-orientation of $V$.''

We now return to the case $E=K$, with the aim of identifying the
twisting $\tau_{V}$ with type defined in \S\ref{sec:twistings}.  The
main point is that the action of the orthogonal group $O(k)$ on
$\Omega^{k}\fred{n}$ lifts through the Atiyah-Singer map
$\fred{k+n}\to \Omega^{k}\fred{n}$.  Our discussion of this matter is
inspired by the Stoltz-Teichner~\cite{stolz04:_what} description of
$\spin$-structures, and, of course
Donovan-Karoubi~\cite{donovan70:_graded_brauer_k}.

Let $X$ be a local quotient groupoid, and $V$ a real vector bundle
over $X$ of dimension $k$, and $\cl(V)$ the associated bundle of
Clifford algebras.  The bundle $\cl(V)\otimes H$ is 
a locally universal $\cl(V)$-module.  The Atiyah-Singer
construction~\cite{atiyah69:_index_fredh}  gives a map
\[
\fredholm_{\cl(V)}(\cl(V)\otimes H) \to
\Omega^{V}\fred{0}(\cl(V)\otimes H)
\]
which is a weak equivalence on global sections.  We can therefore
trivialize the bundle of spectra $\Omega^{V} K$ by trivializing the
bundle of Clifford algebras $\cl(V)$.  

Of course something weaker will also trivialize $\Omega^{V}K$.  We
don't really need a bundle isomorphism $\cl(V)\approx X\times C_{k}$.
We just need a way of going back and forth between $\cl(V)$-modules
and $C_{k}$ modules.  It is enough to have a bundle of irreducible
$\cl(V)\!-\!C_{k}$ bimodules giving a Morita equivalence.

Let $M=C_{k}$, regarded as a $\cl(\R^{k})\!-\!C_{k}$-bimodule.  We
equip $M$ with the Hermitian metric in which the monomials in the
$\epsilon_{i}$ are orthonormal.  Consider the group $\pinc(k)$ of
pairs $(t,f)$ in which $t:\R^{k}\to \R^{k}$ is an orthogonal map, and
\[
f:t^{\ast}M\to M
\]
is a unitary bimodule isomorphism.  The group $\pinc(k)$ is a graded
central extension of $O(k)$, graded by the sign of the determinant.

We now identify the twisting $\tau_{V}$ in the terms of
\S\ref{sec:twistings-1}.  Let $E\to X$ be the bundle of orthonormal
frames in $V$.  Thus $E\to X$ is a principal bundle with structure
group $O(k)$.  Write $P=E\modmod O(k)$, $\tilde P=E\modmod\pinc(k)$.
Then
\[
P\to X
\]
is a local equivalence. and $\tilde P\to P$ is a graded central
extension, defining a twisting $\tau$.  Over $\tilde P$ we can form the
bundle of bimodules 
\[
\tilde M = (E\times M) \modmod \pinc(k),
\]
giving a Morita equivalence between bundles of $\cl(V)$-modules and
bundles of $C_{k}$-modules.  In particular, 
\[
H' = \hom_{\cl(V)}(\tilde M,\cl(V)\otimes H)
\]
is a locally universal $\tau$-twisted
$C_{k}$-module, and the map
\begin{align*}
\Gamma(\fredholm_{\cl(V)}(\cl(V)\otimes H)) & \to 
\Gamma(\fred{k}(H')) \\
T & \mapsto T\circ(\slot)
\end{align*}
is a homeomorphism.  Thus the group $K^{-\tau_{V}+n}(X)$ is isomorphic
to the twisted $K$-group $K^{-\tau+n}(X)$, and, as in
Donovan-Karoubi~\cite{donovan70:_graded_brauer_k} we have a
tautological Thom isomorphism
\[
K^{n}(X^{V})\approx K^{-\tau+n}(X).
\]
More generally, the same construction leads to a tautological Thom
isomorphism 
\begin{equation}\label{eq:38}
K^{\sigma+n}(B(V),S(V)) \approx
K^{-\tau+\sigma+n}(X),
\end{equation}
when $V$ is a vector bundle over a groupoid $X$. 

With the Thom isomorphism in hand, one can define the pushforward, or
umkehr map in the usual way.    Let $f:X\to Y$ be a map of smooth
manifolds, or a map of groupoids forming a bundle of smooth manifolds,
$T=T_{X/Y}$ the corresponding relative (stable) tangent bundle, and
$\tau_{0}$ the twisting on $X$ corresponding to $T$.    Given a
twisting $\tau$ on $Y$, and an isomorphism $f^{\ast}\tau\approx
\tau_{0}$ one can combine the Pontryagin-Thom collapse with the
Thom-isomorphism to form a pushforward map
\[
f_{!}:K^{f^{\ast}\sigma+n}(X)\to K^{-\tau+\sigma+n}(Y)
\]
where $\sigma$ is any twisting on $Y$.  We leave the details to the
reader.

We apply this to the situation in which $X=(G\times G)\modmod G$,
$Y=G\modmod G$ (both with the adjoint action) and $X\to Y$ is the
multiplication map $\mu$.  In this case the twisting $\tau_{0}$ can be
taken to be the twisting we denoted $\mathfrak g$ in the introduction.
Since $\mathfrak g$ is pulled back from $\text{pt}\modmod G$, 
there are canonical isomorphisms 
\[
\mu^{\ast}\mathfrak g\approx p_{1}^{\ast}\mathfrak g\approx
p_{2}^{\ast}\mathfrak g.
\]
We'll just write $\mathfrak g$ for any of these twistings.  Suppose
$\sigma$ is any twisting of $G\modmod G$ which is ``primitive'' in the
sense that it comes equipped with an associative isomorphism
$\mu^{\ast}\sigma\approx p_{1}^{\ast}\sigma+p_{2}^{\ast}\sigma$.  Then
the group $K^{\sigma+\mathfrak g}_{G}(G)$ acquires a Pontryagin
product
\[
K^{\sigma+\mathfrak g}_{G}(G)\otimes K^{\sigma+\mathfrak g}_{G}(G)\to
K^{\mu^{\ast}\sigma+2\mathfrak g}_{G}(G\times G)
\xrightarrow{\mu_{!}}{}
K^{\sigma+\mathfrak g}_{G}(G),
\] 
making it into an algebra over $K_{G}^{0}(\text{pt})=R(G)$.  

\subsection{The fundamental spectral sequence}\label{sec:fund-spectr-sequ}

Our basic technique of computation will be based on 
a variation of the Atiyah-Hirzebruch spectral sequence, which
is constructed using the technique of
Segal~\cite{Segal:CSandSS}.  The identification of the $E^{2}$-term
depends only on the properties listed in \S\ref{sec:axioms}. 

Suppose that $X$ is a local quotient groupoid, and 
write $\prek^{\tau+t}$ for the presheaf on $\cms{X}$ given by
\[
\prek^{\tau+t}(U)= K^{\tau+t}(X_{U}).
\]
Write
\[
\sk^{\tau+t}=\ass\prek^{\tau+t}
\]
for the associated sheaf.  The limit of the Mayer-Vietoris spectral
sequences associated to the (hyper-)covers of $\cms{X}$ is a spectral
sequence 
\[
H^{s}(\cms{X};\sk^{\tau+t})\implies K^{\tau+s+t}(X).
\]
Since $X$ admits locally contractible slices the stalk
of $\sk^{\tau+t}$ at a point $c\in \cms{X}$ is 
\[
K^{\tau+t}(X_{c})\approx 
\begin{cases}
0&\qquad t \text{ odd}\\
\rep^{\tau}(G_{x})&\qquad t \text{ even},
\end{cases}
\]
where $x\in X_{0}$ is a representative of $c$, and $G_{x}=X(x,x)$.

There is also a relative version.  Suppose that $A\subset X$ is a pair
of groupoids, and write $\sk^{\tau+t}_{\text{rel}}$ for the sheaf on
$\cms{X}$ associated to the presheaf
\[
U\mapsto K^{\tau+t}(X_{U}, X_{\cms{A}\cap U}).
\]
Then the limit of the Mayer-Vietoris spectral sequences associated to
the hyper-covers of $\cms{X}$ gives 
\[
H^{s}(\cms{X};\sk^{\tau+t}_\text{{rel}})\implies K^{\tau+s+t}(X,A).
\]
This spectral sequence is most useful when $A\subset X$ is closed, and
has the property that for all sufficiently small $U\subset \cms{X}$,
the map $K^{\tau+t}(X_{U})\to K^{\tau+t}(X_{\cms{A}\cap U})$ is surjective.
In that case there is (for sufficiently small $U$) a short exact
sequence 
\[
K^{\tau+t}(X_{U}, X_{\cms{A}\cap U})\to K^{\tau+t}(X_{U})\to
K^{\tau+t}(X_{\cms{A}\cap U})
\]
and the sheaf $\sk^{\tau+t}_{\text{rel}}$ can be identified with the
extension of $i^{\ast}\sk^{\tau+t}$ by zero
\[
\sk^{\tau+t}_{\text{rel}}=i_{!}(\sk^{\tau+t}),
\]
where $i:V\subset \cms{X}$ is the inclusion of the complement of $A$.
We will make use of this situation in the proof of
Proposition~\ref{thm:10}.

\section{Computation of $K_G^{\tau}(G)$}
\label{sec:comp-k_gg-twist}

The aim of this section is to compute the groups
$K_{G}^{\tau+\ast}(G)$ for non-degenerate $\tau$.  We'll begin by
considering general twistings, and adopt the non-degeneracy hypothesis
as necessary.  Our main results are Theorem~\ref{thm:46},
Corollary~\ref{thm:17} and Corollary~\ref{thm:50}.  

\subsection{Notation and assumptions}
\label{sec:notation-assumptions}

We first fix some notation.  Let
\begin{itemize}
\item $G$ be a compact connected Lie group;
\item $\mathfrak{g}$ the Lie algebra of $G$;
\item $T$ a fixed maximal torus of $G$;
\item $\mathfrak{t}$ the Lie algebra of $T$;
\item $N$ the normalizer of  $T$;
\item $W=N/T$ the Weyl group;
\item $\Pi=\ker \exp:\mathfrak{t}\to T$;
\item $\Lambda=\hom(\Pi,Z)$, the character group of $T$;
\item $\na = \Pi\rtimes NT$
\item $\wa = \Pi\rtimes W=\na/T$, the extended affine Weyl group;
\end{itemize}
The group $\wa$ can be identified the group of symmetries of
$\mathfrak{t}$ generated by translations in $\Pi$ and the reflections
in $W$.  When $G$ is connected, the exponential map, from the orbit
space $\mathfrak{t}/\wa$ to the space of conjugacy classes in $G$, is
a homeomorphism.

We will make our computation for groups satisfying the
equivalent conditions of the following lemma.

\begin{lem}\label{thm:1}
For a Lie group $G$ the following are equivalent
\begin{thmList}
\item For each $g\in G$ the centralizer $Z(g)$ is connected; 
\item $G$ is connected and $\pi_{1}G$ is torsion free;
\item $G$ is connected and any central extension 
\[
\T\to G^{\tau}\to G
\]
splits.
\end{thmList}
\end{lem}

\begin{pf}
The equivalence of~(ii) and~(iii) is elementary: Since $G$~is
connected, its classifying space~$BG$ is simply connected, and from
the Hurewicz theorem and the universal coefficient theorem the torsion
subgroup of~$\pi _1G$ is isomorphic to the torsion subgroup
of~$H_2(BG)$, and so to the torsion subgroup of~$H^3(BG;\Z)$.  But for
any compact Lie group the odd dimensional cohomology of the
classifying space is torsion---the real cohomology of the classifying
space is in even degrees (and is given by invariant polynomials on the
Lie algebra).
 
The implication (ii)$\implies$(i) is~\cite[(3.5)]{bourbaki82:_elemen}.
For the converse (i)$\implies$(ii) we note first that $G=Z(e)$ is
connected by hypothesis.  Let $G'\subset G$ denote the derived
subgroup of $G$, the connected Lie subgroup generated by commutators
in $G$, and $Z_1\subset G$ the connected component of the center of
$G$.  Set $A=Z(G')\cap Z_1$.  Then from the principal fiber bundle
$G'\to G\to Z_1/A$ we deduce that the torsion subgroup of $\pi _1G$ is
$\pi _1G'$.  We must show the latter vanishes.  Now the inclusion $\pi
_1Z(g)\to \pi _1G$ is surjective for any $g\in G$, since any
centralizer contains a maximal torus $T$ of $G$ and the inclusion $\pi
_1T\to\pi _1G$ is surjective---the flag manifold $G/T$ is simply
connected.  It follows that $Z(g)$ is connected if and only if the
conjugacy class $G/Z(g)$ is simply connected.  Furthermore, the
conjugacy class in $G$ of an element of $G'$ equals its conjugacy
class in $G'$, from which we deduce that all conjugacy classes in $G'$
are connected and simply connected.  Let $\widetilde{G'}$ denote the
simply connected (finite) cover of $G'$.  Then the set of conjugacy
classes in $\widetilde{G'}$ may be identified as a bounded convex
polytope in the Lie algebra of a maximal torus, and furthermore $\pi
_1G'$ acts on it by affine transformations with quotient $\widetilde{G'}/G'$;
see~\cite[\S3.9]{duistermaat00:_lie}.  The center of mass of the vertices of the
polytope is a fixed point of the action.  The finite group $\pi _1G'$
acts freely on the corresponding conjugacy class of $\widetilde{G'}$
with quotient a conjugacy class in $G'$.  Since the former is
connected and the latter simply connected, it follows that $\pi_1G'$
is trivial, as desired.
\end{pf}

\subsection{The main computation}
\label{sec:main-computation}

Let $X=G\modmod{G}$ be the groupoid formed from $G$ acting on itself
by conjugation.  We will compute $K^{\tau+\ast}(X)=K^{\tau+\ast}_{G}(G)$
using the spectral sequence described in
\S\ref{sec:fund-spectr-sequ}.  In this case it takes the form 
\begin{equation}\label{eq:6}
H^{s}(G/G;\sk^{\tau+t})\implies K^{\tau+s+t}_{G}(G).  
\end{equation}

The orbit space $G/G$ is the space of conjugacy classes in $G$,
which is homeomorphic via the exponential map to $\mathfrak{t}/\wa$.
Our first task is to identify the sheaf $\sk^{\tau+t}$ on
$G/G\approx\mathfrak{t}/\wa$.

Since $G\modmod{G}$ admits locally contractible slices, the stalk of
$\sk^{\tau+t}$ at a conjugacy class $c\in G/G$ is the twisted
equivariant $K$-group
\[
K^{\tau+t}_{G}(c).
\]
A choice of point $g\in c$ gives an identification $c=G/Z(g)$, and an
isomorphism
\begin{equation}\label{eq:21}
K^{\tau+t}_{G}(c)\approx 
K^{\tau_{g}+t}_{Z(g)}(\{g \})\approx 
\begin{cases}
\rep^{\tau_{g}}(Z(g)) &\quad t \text{ even} \\
0 &\quad t \text{ odd} \\
\end{cases}
\end{equation}
We have denoted by $\tau_{g}$ the restriction of $\tau$ to $\{g
\}\modmod Z(g)$, in order to emphasize the dependence on the choice of
$g$.  Among other things, this proves that
\[
\sk^{\tau+\text{odd}}=0.
\]

The twisting $\tau_{g}$ corresponds to a graded central extension
\begin{equation}\label{eq:9}
\T\to \tilde Z(g)\to Z(g).
\end{equation}
The group $Z(g)$ has $T$ for a maximal torus, and is connected when
$G$ satisfies the equivalent conditions of Lemma~\ref{thm:1}.  Denote
\begin{equation}\label{eq:11}
\T\to\tilde T\to T
\end{equation}
the restriction of~\eqref{eq:9} to $T$.  Then $\tilde T$ is a maximal
torus in $\tilde Z(g)$.  The map from the Weyl group of $\tilde Z(g)$
to the Weyl group $W_{g}$ of $Z(g)$ is an isomorphism, and
\[
\rep^{\tau_{g}}(Z(g))\to \rep^{\tau_{g}}(T)^{W_{g}}
\]
is an isomorphism.  We can therefore re-write~\eqref{eq:21} as
\begin{equation}\label{eq:22}
K^{\tau+t}_{G}(c)\approx 
\begin{cases}
\rep^{\tau_{g}}(T)^{W_{g}} &\quad t \text{ even} \\
0 &\quad t \text{ odd} \\
\end{cases}
\end{equation}

We now reformulate these remarks in order to eliminate the explicit
choice of $g\in c$.  We can cut down the size of $c$ by requiring that
$g$ lie in $T$.  That helps, but it doesn't eliminate the dependence
of $\tau_{g}$ on $g$.  We can get rid of the reference to $g$ by
choosing a geodesic in $T$ from each $g$ to the identity element, and
using it to identify the twisting $\tau_{g}$ with $\tau_{0}$.  This
amounts to considering the set of elements of $\mathfrak{t}$ which
exponentiate into $c$.  This set admits a transitive action of $\wa$,
and the stabilizer of an element $v$ is canonically isomorphic to
$W_{g}$ where $g=\exp(v)$.

We are thus led to look at the groupoid $\mathfrak{t}\modmod T$, and
the action of $\wa$.  Writing it this way, however, does not
conveniently display the action of $\wa$ on the twisting $\tau$.
Following Example~\ref{eg:21}, we work instead with the weakly
equivalent groupoid $(\wa\times\mathfrak{t})\modmod \na$.

Consider the map
\[
K^{\tau+t}_{G}(G)\to
K^{\tau+t}_{\na}(\wa\times\mathfrak{t})
\]
induced by 
\[
(\wa\times\mathfrak{t})\modmod\na\xrightarrow{\text{projection}}
\mathfrak{t}\modmod\na\xrightarrow{\text{exp}}
G\modmod G
\]
Since the right action of $\wa=\na/T$ commutes with the diagonal
left action of $\na$ on $\wa\times\mathfrak{t}$, the group $\wa$
acts on the groupoid $\wa\times \mathfrak{t}\modmod\na$.  The twisting $\tau$
is fixed by this action since it is pulled back from $G\modmod G$.  The left
action of $\wa$ on $(\wa\times\mathfrak{t})\modmod\na$ therefore induces
a right action of $\wa$ on $K^{\tau+t}_{\na}(\mathfrak{t})$, and the
image of $K^{\tau+t}_{G}(G)$ is invariant:
\begin{equation}\label{eq:25}
K^{\tau+t}_{G}(G)\to
K^{\tau+t}_{\na}(\wa\times \mathfrak{t})^{\wa}.
\end{equation}
Since $\wa=\na/T$, the map
\[
\mathfrak t\modmod T \xrightarrow{} (\wa\times\mathfrak{t})\modmod \na
\]
is a local equivalence, and so gives an isomorphism
\[
K^{\tau+t}_{\na}(\wa\times \mathfrak{t})\approx K^{\tau+t}_{T}(\mathfrak{t}).
\]
There is therefore an action of $\wa$ on $K^{\tau+t}_{T}(\mathfrak
t)$, and we may re-write~\eqref{eq:25} as
\[
K^{\tau+t}_{G}(G)\to K^{\tau+t}_{T}(\mathfrak t)^{\wa}.
\]

For $c\in G/G\approx \mathfrak{t}/\wa$, let 
\[
S_{c}=\left\{s\in\mathfrak{t}\mid \exp(s)\in c \right\}
\]
be the corresponding $\wa$-orbit in $\mathfrak{t}$.   A similar
discussion gives a map 
\begin{equation}\label{eq:23}
K^{\tau+t}_{G}(c)\to K^{\tau+t}_{N}(\wa\times S_{c})^{\wa} \approx
K^{\tau+t}_{T}(S_{c})^{\wa}.
\end{equation}

\begin{prop}
\label{thm:6}
If $G$ satisfies the conditions of Lemma~\ref{thm:1}
then the map 
\[
K^{\tau+t}_{G}(c)\to K^{\tau+t}_{T}(S_{c})^{\wa}
\]
constructed above is an isomorphism.
\end{prop}

\begin{pf}
Choose $v\in S_{c}$, and let $W_{v}\subset\wa$ be the stabilizer of
$v$.  We then have an identification
$S_{c}\approx \wa/W_{v}$, and so an isomorphism
\[
K^{\tau+t}_{T}(S_{c})^{\wa}\approx
K^{\tau+t}_{T}(\{v \})^{W_{v}}.
\]
Write $g=\exp(v)$.  The restriction of $\na\to G$ identifies $W_{v}$
with the Weyl group of $Z(g)$, $\{v \}\modmod T$ with $\{g \}\modmod
T$, and the restriction of $\tau$ to $\{v \}\modmod T$ with
$\tau_{g}$.  By Example~\ref{eg:22} action of $W_{v}$ on
$K^{\tau+t}_{T}(\{v \})$ coincides with the action $W_{v}$ by
conjugation.  The result then follows from~\eqref{eq:22}.
\end{pf}

We now identify the sheaf $\sk^{\tau+t}$.  Since $\{0
\}\to\mathfrak{t}$ is an equivariant homotopy equivalence, the
restriction map
\[
K^{\tau+t}_{T}(S_{c}\times \mathfrak{t})
\to K^{\tau+t}_{T}(S_{c}\times \{0 \})
\]
is an isomorphism.    Next note that the aggregate of the restriction
maps to the points of $S_{c}$ gives a map from 
\[
K^{\tau+t}_{T}(S_{c}\times\mathfrak{t})^{\wa}
\]
to the set of $\wa$-equivariant maps
\[
S_{c}\to K^{\tau+t}_{T}(\mathfrak{t}),
\]
which, using the fact that $\wa$ acts transitively on $S_{c}$, is
easily checked to be an isomorphism.  Write
$p:\mathfrak{t}\to\mathfrak{t}/\wa$ for the projection, and for an
open $U\subset G/G=\mathfrak{t}/\wa$ set
\[
S_{U} = p^{-1}(U).
\]
Let $\st^{\tau+t}$ be the presheaf which associates to $U\subset G/G$
the set of locally constant $\wa$-equivariant maps
\[
S_{U}\to K_{T}^{\tau+t}(\mathfrak{t}).
\]
There is then a map of presheaves
\[
\prek^{\tau+t}\to\st^{\tau+t},
\]
hence a map of sheaves
\begin{equation}\label{eq:8}
\sk^{\tau+t}\to\st^{\tau+t}.
\end{equation}

\begin{cor}
The map~\eqref{eq:8} is an isomorphism.  
\end{cor}

\begin{pf}
Proposition~\ref{thm:6} implies that~\eqref{eq:8} is an isomorphism of
stalks, hence an isomorphism.
\end{pf}

We now re-interpret the sheaf $\st$ in a form more suitable to
describing its cohomology.  Since
\[
\twist{\mathfrak{t}\modmod T}\to
\twist{\{0 \}\modmod T}
\]
is an equivalence of categories, the restriction of $\tau$ to
$\mathfrak{t}\modmod T$ corresponds to a graded central extension
\begin{equation}\label{eq:24}
\T\to T^{\tau}\to T
\end{equation}
equipped with an action of $\wa$.  The $\wa$-action fixes $\T$ and
acts on $T$ through its quotient $W$, the Weyl group.  Write
$\LambdaTorsor$ for the set of splittings of~\eqref{eq:24}.  Note that
$\LambdaTorsor$ is a torsor for $\Lambda$ and inherits a compatible
$\wa$ action from~\eqref{eq:24}.

By Proposition~\ref{thm:18} the group
\[
K^{\tau+0}_{T}(\mathfrak{t})\approx
K^{\tau+0}_{T}(\{0 \})
\]
may be identified with with the set of compactly supported functions
on $\LambdaTorsor$ with values in $\Z$.  We will see shortly that the
action of $\wa$ is the combination of its natural action on
$\LambdaTorsor$ and an action on $\Z$ given by a homomorphism
$\epsilon:\wa\to\Z/2$.  Writing $\Z(\epsilon)$ for the sign
representation associated to $\epsilon$, we then have an isomorphism
of $\wa$-modules
\begin{equation}\label{eq:4}
K^{\tau+0}_{T}(\{0 \})\approx \Hom_{c}(\LambdaTorsor,\Z(\epsilon)).
\end{equation}

To verify the claim about the action first note that the automorphism
group of the restriction of $\tau$ to $\mathfrak{t}\modmod T$ is
\[
H^{2}(BT;\Z)\times H^{0}(BT;\Z/2) \approx \Lambda\times \Z/2\approx
R(T)^{\times} \approx K_{T}^{0}(\text{pt})^{\times}
\]
By Part \thmItemref{item:one-morphism} of Proposition~\ref{thm:19},
the factor $\Lambda$ acts on $K^{\tau+0}_{T}(\{0\})$ through its
natural action on $\LambdaTorsor$, while the $\Z/2$ acts by its sign
representation.

Since $\wa=\Pi\rtimes W$, and the action of $\wa$
on~$K^{\tau+0}_{T}(\{0 \})$ is determined by its restriction to $\Pi$
and $W$.  The group $\Pi$ acts trivially on $T$ and so it acts
on~$K^{\tau+0}_{T}(\{0 \})$ through a homomorphism
\[
\Pi\xrightarrow{(b,\epsilon_{\Pi})}{} \Lambda\times\Z/2.  
\]
The group $W$ does act on $T$, and so on 
\[
H^{2}(BT;\Z)\times H^{0}(BT;\Z/2) \approx \Lambda\times \Z/2, 
\]
by the product of the natural (reflection) action on
$\Lambda$ and the trivial action on $\Z/2$.    The restriction of the
action of $\wa$ to $W$ is therefore determined by a crossed
homomorphism 
\[
W\to \Lambda
\]
compatible with $b$, and an ordinary homomorphism $\epsilon_{W}:W\to
\Z/2$.   The maps $\epsilon_{\Pi}$ and $\epsilon_{W}$ combine to give
the desired map $\epsilon:\wa\to \Z/2$, while the map
$b:\Pi\to\Lambda$ and the crossed homomorphism $W\to\Lambda$
correspond to the natural action of $\wa$ on $\LambdaTorsor$.  This
verifies the isomorphism~\eqref{eq:4} of $\wa$-modules.

We can now give a useful description of $\st$.  First recall a
construction.  Suppose $X$ is a space equipped with an action of a
group $\Gamma$, and that $\mathcal G$ is an equivariant sheaf on $X$.
Write $p:X\to X/\Gamma$ for the projection to the orbit space.  There
is then a sheaf, $\left(p_{\ast}\mathcal G \right)^{\Gamma}$ on
$X/\Gamma$ whose value on an open set $V$ is the set of
$\Gamma$-invariant elements of $\mathcal G(p^{-1}V)$.  A very simple
situation is when $\mathcal G$ is the constant sheaf $\Z$.  In that
case $\left(p_{\ast}\mathcal G \right)^{\Gamma}$ is again the constant
sheaf $\Z$.  This will be useful in the proof of
Proposition~\ref{thm:56} below.

\begin{cor}
\label{thm:55}
Write
\[
\tilde{\mathfrak{t}}= \mathfrak{t}\times_{\wa}\LambdaTorsor
\]
and let 
\begin{align*}
p &:\mathfrak{t}\times\LambdaTorsor\to\tilde{\mathfrak{t}}
\quad\text{and} \\
f &:\tilde{\mathfrak{t}}\to \mathfrak{t}/\wa
\end{align*}
denote the projections.  There is a canonical isomorphism
\[
\st^{\tau+0}\approx
f_{\ast}^{c}\left(p_{\ast}\Z(\epsilon)^{\wa}\right),
\]
where $f_{\ast}^{c}$ denotes pushforward with proper supports. 
\qed
\end{cor}

To go further we need to make an assumption.
\begin{assumption}\label{ass:non-degeneracy}
The twisting $\tau$ is non-degenerate in the sense that $b$ is a monomorphism.
\end{assumption}

In terms of the classification of twistings, this is equivalent to
requiring that the image of the isomorphism class of $\tau$ in
\[
H^{3}_{T}(T;\R)\approx H^{1}(T;\R)\otimes H^{1}(T;\R)
\]
is a non-degenerate bilinear form.

Next note
\begin{lem}\label{thm:53}
The map $\epsilon_{W}:W\to\Z/2$ is trivial.
\end{lem}

\begin{pf}
The homomorphism in question corresponds to the element in
$H^{1}_{W}(\text{pt})=H^{1}(BW;\Z/2)$ given by restricting the
isomorphism class of the twisting along
\begin{multline*}
H^{1}_{G}(G;\Z/2)\times H^{3}_{G}(G;\Z) \to 
H^{1}_{G}(G;\Z/2) \to 
H^{1}_{G}(\{e \};\Z/2)  \\
\to 
H^{1}_{N}(\{e \};\Z/2) \approx 
H^{1}_{W}(\{e \};\Z/2).
\end{multline*}
Since $G$ is assumed to be connected $H^{1}_{G}(\text{pt})=0$ and the
result follows.
\end{pf}

\begin{cor}\label{thm:52}
There is an isomorphism $\Z\approx \Z(\epsilon)$ of equivariant sheaves
on $\mathfrak{t}\times\LambdaTorsor$.
\end{cor}

\begin{pf}
The sheaf $\Z(\epsilon)$ is classified by the element
\[
\tilde\epsilon\in H^{1}_{\wa}(\mathfrak{t}\times\LambdaTorsor;\Z/2)
\]
pulled back from the $\epsilon\in H^{1}_{\wa}(\text{pt};\Z/2)$.  By
Corollary~\ref{thm:53} the restriction of $\tilde\epsilon$ to $W$ is
trivial.  By Assumption~\ref{ass:non-degeneracy}, the group $\Pi$
acts freely on $\LambdaTorsor$, so the restriction of $\tilde\epsilon$
to $H^{1}_{\wa}(\mathfrak{t}\times\LambdaTorsor;\Z/2)$ is also
trivial.  This proves that $\tilde\epsilon=0$.
\end{pf}

\begin{prop}\label{thm:56}
There is an isomorphism
\[
\st^{\tau+0}\approx f_{\ast}^{c}(\Z),
\]
where
$f:\tilde{\mathfrak{t}}=\mathfrak{t}\times_{\wa}\LambdaTorsor\to\mathfrak{t}/W$
is the projection and $f_{\ast}^{c}$ denotes pushforward with proper
supports.
\end{prop}

\begin{pf}
By Corollaries~\ref{thm:52} and~\ref{thm:55} there is an isomorphism
\[
\st^{\tau+0}\approx f_{\ast}^{c}\left(p_{\ast}\Z^{\wa}\right),
\]
so the result follows from the fact that $p_{\ast}\Z^{\wa}\approx \Z$.
\end{pf}

Finally, using the existence of contractible local slices we can
describe the cohomology of $\sk^{\tau+0}\approx\st^{\tau+0}$.
\begin{lem}\label{thm:11}
The edge homomorphism of the Leray spectral sequence for $f$ is an
isomorphism
\[
H^{\ast}(\mathfrak{t}/\wa;\sk^{\tau})\approx
H^{\ast}_{c}(\tilde{\mathfrak{t}};\Z),
\]
where $H^{\ast}_{c}$ denotes cohomology with compact supports. \qed
\end{lem}

To calculate $H^{\ast}_{c}(\tilde{\mathfrak{t}};\Z)$ we need to understand
the structure of $\tilde{\mathfrak{t}}$.  This amounts to describing
more carefully the action of $\wa$ on $\LambdaTorsor$.

\begin{lem}\label{thm:9}
There exists an element $\efix\in \LambdaTorsor$ fixed by $W$.
\end{lem}

\begin{pf}
The inclusion $\{0 \}\modmod T\subset \mathfrak{t}\modmod T$ is
equivariant for the action of $W$.  We can therefore study the action
of $W$ on the central extension of $T$ defined by the restriction of
$\tau$ to $\{0\}\modmod T$.  Now $\tau$ started out as a twisting of
$G\modmod G$, so our twisting of $\{0\}\modmod T$ is the restriction
of a twisting $\tau_{G}$ of $\{e \}\modmod G$.  Moreover, the action
of $W$ is derived from the action of inner automorphisms of $G$ on
$\tau_{G}$.  Now the twisting $\tau_{G}$ corresponds to a central
extension
\[
\T\to G^{\tau}\to G.
\]
By our assumptions on $G$ (Lemma~\ref{thm:1}), this central extension
splits.  Choose a splitting 
\begin{equation}\label{eq:13}
G^{\tau}\to \T
\end{equation}
and let $\efix\in \LambdaTorsor$ be the composition
\[
\tilde T\to G^{\tau}\to \T.
\]
Since $\T$ is abelian, the splitting~\eqref{eq:13} is preserved by
inner automorphisms of $G$.   It follows that splitting $e$ is fixed
by the inner automorphisms of $\tilde G$ which normalize $\tilde T$.
The claim follows.
\end{pf}

\begin{rem}
Any two choices of $\efix$ differ by a character of $G$, so the
element $\efix$ is unique if an only if the character group of $G$ is
trivial. Since we've assumed that $G$ is connected and $\pi_{1}G$ is
torsion free, this is in turn equivalent to requiring that $G$ be
simply connected.
\end{rem}

\begin{rem}
A {\em primitive} twisting $\tau$ comes equipped with a trivialization
of its restriction to $\{e \}\modmod G$, or in other words a splitting
of the graded central extension $G^{\tau}\to G$.   A primitive
twisting therefore comes equipped with a canonical choice of $\efix$.
\end{rem}

Using a fixed choice of $\lambda_{0}$, we can identify $\LambdaTorsor$
with $\Lambda$ as a $W$-space.  To sum up, we can make an
identification $\LambdaTorsor\approx \Lambda$, the action of $\Pi$ is
given by a $W$-equivariant homomorphism $\Pi\to\Lambda$, and the
$W$-action is the natural one on $\Lambda$.

\begin{lem}\label{thm:14}
When $\tau$ is non-degenerate the $\wa$-set $\LambdaTorsor$ admits an
(equivariant) embedding in $\mathfrak{t}$.  There are finitely many
$\wa$-orbit in $\LambdaTorsor$, and each orbits is of the form
$\wa/W_{c}$, with $\left(W_{c},\mathfrak{t} \right)$ a finite (affine)
reflection group.
\end{lem}

\begin{pf}
Since $b$ is a monomorphism the map $\mathfrak{t}=\Pi\otimes\R\to
\Lambda\otimes\R$ is an isomorphism.  The first assertion now follows
from our identification of $\LambdaTorsor$ with $\Lambda$.  As for the
finiteness of the number of orbits, since $b$ is a monomorphism, the
group $\Lambda/\Pi$ is finite, and there are already only finitely
many $\Pi$-orbits in $\LambdaTorsor$.  The remaining assertions follow
from standard facts about the action of $\wa$ on $\mathfrak{t}$
(Propositions~\ref{thm:12} and~\ref{thm:13} below).
\end{pf}

\begin{cor}
When $\tau$ is non-degenerate, there is a homeomorphism
\[
\tilde{\mathfrak{t}}\equiv\coprod_{s\in S}\mathfrak{t}/W_{s}
\]
with $S$ finite, and $W_{s}$ a finite reflection group of isometries
of $\mathfrak{t}$.  Moreover 
\[
\mathfrak{t}/W_{s}\equiv \R^{n_{1}}\times[0,\infty)^{n_{2}}
\]
with $n_{2}=0$ if and only if $W_{s}$ is trivial.
\end{cor}

\begin{pf}
This is immediate from Lemma~\ref{thm:14} above and
Proposition~\ref{thm:13} below.
\end{pf}

Since $H_{c}^{\ast}([0,\infty);\Z)=0$ and
\[
H_{c}^{\ast}(\R;\Z) = \begin{cases}\Z\quad&\ast=1\\0\quad&\text{otherwise,}\end{cases}
\]
the Kunneth formula gives
\[
H_{c}^{\ast}(\R^{n_{1}}\times[0,\infty)^{n_{2}};\Z)=
\begin{cases}
\Z &\quad \text{$n_{2}=0$ and $\ast=n_{1}$} \\
0 &\quad\text{otherwise}.
\end{cases}
\]

In summary, we have
\begin{prop}\label{thm:15}
The cohomology group $H^{\ast}_{c}(\tilde{\mathfrak{t}};\Z)$ is zero
unless $\ast=n$, and $H^{n}_{c}(\tilde{\mathfrak{t}};\Z)$ is isomorphic
to the free abelian group on the set of free $\wa$-orbits in
$\LambdaTorsor$.   More functorially, 
\[
H^{n}_{c}(\tilde{\mathfrak{t}};\Z)\approx
\Hom_{\wa}\left(\LambdaTorsor, H^{n}_{c}(\mathfrak{t})\otimes\Z(\epsilon) \right).
\]
\qed
\end{prop}

Proposition~\ref{thm:15} implies that the spectral
sequence~\eqref{eq:6}  collapses, giving
\begin{thm}\label{thm:46}
Suppose that $G$ is a Lie group of rank $n$ satisfying the conditions
of Lemma~\ref{thm:1}, and that $\tau$ is a non-degenerate twisting of
$G\modmod G$, classified by 
\[
[\tau]\in H^{3}_{G}(G;\Z)\times H^{1}_{G}(G;\Z/2).
\]
The restriction of $\tau$ to $\point\modmod T$
determines a central extension
\begin{equation}\label{eq:14}
\T\to T^{\tau}\to T
\end{equation}
with an action of $\wa$.  Write $\LambdaTorsor$ for the set of
splittings of~\eqref{eq:14},  $\epsilon:\wa\to\Z/2$  for the map
corresponding to the restriction of
$[\tau]$ to
\[
H^{1}_{N}(T;\Z/2)\approx H^{1}(\wa;\Z/2),
\]
and $\Z(\epsilon)$ for the associated sign representation.  Then
$K^{\tau+n+1}_{G}(G)=0$, and the twisted $K$-group $K^{\tau+n}_{G}(G)$
is given by
\[
K^{\tau+n}_{G}(G)\approx \Hom_{\wa}\left(\LambdaTorsor,
H^{n}_{c}(\mathfrak{t})\otimes\Z(\epsilon)\right),
\]
which can be identified with the free abelian group on the set of free
$\wa$-orbits in $\LambdaTorsor$, after choosing a point in each free
orbit.  This isomorphism is natural in the sense that if $i:H\subset
G$ is a subgroup of rank $n$ also satisfying the conditions of
Lemma~\ref{thm:1}, then the restriction map
\[
K^{\tau+n}_{G}(G)\to K^{\tau+n}_{H}(H)
\]
is given by the inclusion
\[
\Hom_{\wa(G)}\left(\LambdaTorsor, H^{n}_{c}(\mathfrak{t})\otimes\Z(\epsilon) \right) \subset
\Hom_{\wa(H)}\left(\LambdaTorsor, H^{n}_{c}(\mathfrak{t})\otimes\Z(\epsilon) \right).
\]
\qed
\end{thm}

\begin{rem}\label{rem:3}
In Theorem~\ref{thm:46} the group $\wa$ acts on
$H^{n}_{c}(\mathfrak{t})$ through the action of $W$ on $\mathfrak{t}$.
The reflections thus act by $(-1)$ and a choice of orientation on
$\mathfrak{t}$ identifies $H^{n}_{c}(\mathfrak{t})$ with the usual
sign representation of $W$ on $\Z$.
\end{rem}

The group $K^{\tau+n}_{G}(G)$ is a module over $R(G)$.  Our next goal
is to identify this module structure.  Because $G$ is connected we can
identify $R(G)$ with the ring of $W$-invariant elements of
$\Z[\Lambda]$ or with the convolution algebra of compactly supported
functions $\Hom_{c}(\Lambda,\Z)$.  The algebra $\Hom_{c}(\Lambda,\Z)$
acts on
$\Hom(\LambdaTorsor,H^{n}_{c}(\mathfrak{t})\otimes\Z(\epsilon))$ by
convolution, and one easily checks that the $W$-invariant elements
preserve the $\wa$-equivariant functions.

\begin{prop}\label{thm:28}
Under the identification 
\[
K_{G}^{\tau+n}(G)\approx 
\Hom_{\wa}\left(\LambdaTorsor, H^{n}_{c}(\mathfrak{t})\otimes\Z(\epsilon) \right) 
\]
the action of $R(G)\approx \Hom_{c}(\Lambda,\Z)^{W}$ corresponds to
convolution of functions.
\end{prop}

\begin{pf}
This is straightforward to check in case $G$ is a torus.  The case of
general $G$ is reduced to this case by looking at the restriction map
to a maximal torus and using Theorem~\ref{thm:46}.
\end{pf}

Proposition~\ref{thm:28} leads to a very useful description of
$K^{\tau+n}_{G}(G)$.  Choose an orientation of $\mathfrak{t}$ and
hence an identification $H^{n}_{c}(\mathfrak{t})\approx\Z$ of abelian
groups.  By definition, the elements of $\LambdaTorsor$ are characters
of $T^{\tau}$, all of which restrict to the defining character of
$\T$.  To a function
\[
f\in\Hom_{\Pi}\left(\LambdaTorsor, H^{n}_{c}(\mathfrak{t})\otimes\Z(\epsilon) \right) 
\approx \Hom_{\Pi}\left(\LambdaTorsor, \Z(\epsilon) \right) 
\]
we associate the series
\begin{equation}\label{eq:1}
\delta_{f}=\sum_{\lambda\in\LambdaTorsor} f(\lambda)\,\lambda^{-1},
\end{equation}
which is the Fourier expansion of the distribution on $T^{\tau}$
satisfying $\delta_{f}(\lambda)= f(\lambda)$.  Our next aim is to work
out more explicitly which distribution it is, especially when $f$
comes from an element of
\[
\Hom_{\wa}\left(\LambdaTorsor, H^{n}_{c}(\mathfrak{t})\otimes\Z(\epsilon) \right) 
\approx \Hom_{\wa}\left(\LambdaTorsor, \Z(\epsilon) \right). 
\]

Since all of the characters $\lambda$ restrict to the defining
character of the central $\T$, we'll think of the distribution
$\delta_{f}$ as acting on the space of functions $g:T^{\tau}\to\C$
satisfying $g(\zeta\,v)=\zeta g(v)$ for $\zeta\in\T$.  This
space is the space of sections of a suitable complex line bundle
$L_{\tau}$ over $T$.  The character $\chi$ of a representation of $G$
is a function on $T$, and action of $\chi$ on $\delta_{f}$ is given by
\[
\chi\cdot \delta_{f}(g)=\delta_{f}(g\cdot \chi).
\]

Since $\Pi$ and $\Lambda$ are duals, we have
$\hom(\Pi,\Z/2)=\Lambda\otimes\Z/2$, and we may regard
$\epsilon_{\Pi}$ as an element of $\Lambda/2\Lambda$.  This determines
an element 
\[
\lambda_{\epsilon}=\tfrac12\epsilon_{\Pi}\in
\tfrac12\Lambda/\Lambda\subset \Lambda\otimes\R/\Z.
\]
Thinking of $\Pi$ as the character group of $\Lambda\otimes\R/\Z$, 
the function $\epsilon_{\Pi}:\Pi\to\Z/2$ is given by evaluation of
characters on $\lambda_{\epsilon}$:
\[
\epsilon_{\Pi}(\pi)= \pi(\lambda_{\epsilon}).
\]
Since $\epsilon_{\Pi}$ is $W$-invariant, so is $\lambda_{\epsilon}$.

From the embedding $b:\Pi\subset \Lambda$ we get a map
\[
b:T=\Pi\otimes\R/\Z \to \Lambda\otimes\R/\Z.
\]
We'll write $F=\Lambda/\Pi$ for the kernel of this map, and
$F_{\epsilon}$ for the inverse image of $\lambda_{\epsilon}$.  Set
\[
F^{\tau}=\LambdaTorsor/\Pi.
\]
The elements of $F^{\tau}$ can be interpreted as sections of the
restriction of $L_{\tau}$ to $F$.  Finally, let
$F_{\epsilon,\text{reg}} \subset F_{\epsilon}$ and
$F^{\tau}_{\text{reg}}\subset F^{\tau}$ be the subsets
consisting of elements on which the Weyl group $W$ acts freely.

\begin{prop}\label{thm:37}
For 
\[
f\in \Hom_{\Pi}(\LambdaTorsor,H^{n}_{c}(\mathfrak{t})\otimes\Z(\epsilon))\approx
\Hom_{\Pi}(\LambdaTorsor,\Z(\epsilon))
\]
The value of the distribution $\delta_{f}$
on a section $g$ of $L_{\tau}$ is given by
\begin{equation}\label{eq:10}
\delta_{f}(g)= \frac1{|F|}\sum_{(\lambda,x)\in F^{\tau}\times F_{\epsilon}}
f(\lambda) \lambda^{-1}(x) g(x)
\end{equation}
When $f$ is $\wa$-invariant, then
\[
\delta_{f}(g)= \frac1{|F|}\sum_{(\lambda,x)\in F^{\tau}_{\text{reg}}
\times F_{\epsilon,\text{reg}}} f(\lambda) \lambda^{-1}(x) g(x)
\]
\end{prop}

\begin{pf}
Let's first check that~\eqref{eq:10} is well-defined.  Under 
\[
\lambda \mapsto \lambda\pi\qquad \pi\in\Pi,
\]
the term $f(\lambda) \lambda^{-1}(x) g(x)$ gets sent to
\begin{align*}
f(\lambda\pi) (\lambda\pi)^{-1}(x) g(x) &=
\epsilon(\pi)f(\lambda)\lambda^{-1}(x)\pi^{-1}(x) g(x) \\ &= 
\epsilon(\pi)f(\lambda)\lambda^{-1}(x)\epsilon(\pi^{-1}) g(x) \\ &= 
f(\lambda) \lambda^{-1}(x) g(x)
\end{align*}
so $f(\lambda) \lambda^{-1}(x) g(x)$ does indeed depend only
on the $\Pi$-coset of $\lambda$.  
 
To establish~\eqref{eq:10}, it suffices by linearity to consider the
case in which $\delta_{f}$ is of the form
\[
\delta_{f} = \lambda^{-1} \cdot \delta_{\Pi},
\]
with $\lambda\in\LambdaTorsor$,
\[
\delta_{\Pi} = \sum_{\pi\in\Pi}\epsilon(\pi)\pi^{-1},
\]
and $g$ is an element of $\LambdaTorsor$.  In this case $f$ vanishes
off of the $\Pi$-orbit through $\lambda$, and $f(\lambda)=1$.  The
sum~\eqref{eq:10} is then
\begin{equation}\label{eq:20}
\frac1{|F|}\sum_{x\in F_{\epsilon}}\lambda^{-1}(x) g(x).
\end{equation}
By definition, for $\eta\in \Pi$, 
\[
\delta_{\Pi}(\eta)= \eta(\lambda_{\epsilon}) = \frac1{|F|}\sum_{x\in
F_{\epsilon}} \eta(x).
\]
For $\eta\in \Lambda\setminus \Pi$, there is an $a\in F$ with
$\eta(a)\ne 1$.  In that case 
\[
\frac1{|F|}\sum_{x\in F_{\epsilon}} \eta(x) = 
\frac1{|F|}\sum_{x\in F_{\epsilon}} \eta(a x) = 
\eta(a)\frac1{|F|}\sum_{x\in F_{\epsilon}} \eta(x),
\]
so 
\[
\frac1{|F|}\sum_{x\in F_{\epsilon}} \eta(x) = 0.
\]
It follows that for every $\eta\in\Lambda$, 
\[
\delta_{\Pi}(\eta)= \frac1{|F|}\sum_{x\in F_{\epsilon}} \eta(x).
\]
Now suppose $g\in\LambdaTorsor$.  Then 
\[
\delta_{f}(g)= \delta_{\Pi}(\lambda^{-1}g)= \frac1{|F|}\sum_{x\in
F_{\epsilon}} \lambda^{-1}(x) g(x),
\]
which is~\eqref{eq:20}.  This proves the first assertion of
Proposition~\ref{thm:37}.

For the second assertion, note that if $\lambda\in
F^{\tau}_{\text{reg}}$ is fixed by an element of $W$ it is fixed by an
element $w\in W$ which is a reflection.  By Weyl-equivariance, we have
$f(\lambda)=f(w\,\lambda)=w\cdot f(\lambda)=-f(\lambda)$, and so
$f(\lambda) = 0$.  This gives
\[
\delta_{f}(g)=\frac1{|F|}\sum_{(\lambda,x)\in
F^{\tau}_{\text{reg}} \times F} f(\lambda) \lambda^{-1}(x) g(x).
\]
If $x\in F$ is an element fixed by a reflection $w\in W$ then
\begin{align*}
\sum_{\lambda\in F^{\tau}_{\text{reg}}} f(\lambda)\lambda^{-1}(x)g(x) &=
\sum_{\lambda\in F^{\tau}_{\text{reg}}} f(\lambda)\lambda^{-1}(w\cdot
x)g(x) \\
&=
\sum_{\lambda\in F^{\tau}_{\text{reg}}} f(\lambda)(\lambda^{w})^{-1}(x)g(x)
\\
&=\sum_{\lambda\in F^{\tau}_{\text{reg}}} f(\lambda^{w})\lambda^{-1}(x)g(x)
\\
&= -\sum_{\lambda\in F^{\tau}_{\text{reg}}}f(\lambda)\lambda^{-1}(x)g(x)
\end{align*}
so the terms involving such an $x$ sum to zero, and
\[
\delta_{f}(g)=\frac1{|F|}\sum_{(\lambda,x)\in F^{\tau}_{\text{reg}}
\times F_{\epsilon,\text{reg}}} f(\lambda) \lambda^{-1}(x) g(x)
\]
\end{pf}

Let $I^{\tau}\subset R(G)$ be the ideal consisting of virtual
representations whose character vanishes on the elements of
$F_{\epsilon,\text{reg}}$.  

\begin{cor}
\label{thm:48}
The ideal $I^{\tau}$ annihilates $K^{\tau+\ast}_{G}(G)$. 
\end{cor}

\begin{pf}
Write $\chi$ for the character of an element of $I^{\tau}$.  For $f\in
\Hom_{\wa}(\LambdaTorsor,H^{n}_{c}(\mathfrak{t})\otimes\Z(\epsilon))$ we have, by
Proposition~\ref{thm:37}
\[
\chi\,\delta_{f}(g)= \delta_{f}(g\cdot \chi)=
\sum_{(\lambda,x)\in F^{\tau}_{\text{reg}}\times F_{\text{reg}}}
f(\lambda) \lambda^{-1}(x) g(x)\chi(x)=0.
\]
\end{pf}

\begin{rem}
\label{rem:2}
The conjugacy classes in $G$ of the elements in $F_{\epsilon,\text{reg}}$ are known
as the {\em Verlinde conjugacy classes}, and the ideal $I^{\tau}$ as
the {\em Verlinde ideal}.
\end{rem}

\begin{prop}\label{thm:49}
The $R(G)$-module $K^{\tau+n}_{G}(G)$ is cyclic.
\end{prop}

\begin{pf}
Using Lemma~\ref{thm:9} choose a $\wa$-equivariant isomorphism
$\LambdaTorsor\approx \Lambda$, and an orientation of $\mathfrak{t}$
giving an isomorphism $H_{c}^{n}(\mathfrak{t})\approx \Z$.  
We can then identify $K^{\tau+n}_{G}(G)$ with 
\[
\Hom_{\wa}(\Lambda,\Z(\epsilon)),
\]
though we remind the reader that $\wa$ acts on $\Z$ through its sign
representation.  We'll continue the convention of writing elements
\[
f\in\Hom_{\wa}(\Lambda,\Z)
\]
as Fourier series 
\[
\sum f(\lambda)\lambda^{-1}.
\]
Set
\[
\delta_{\Pi}=\sum_{\pi\in\Pi}\epsilon(\pi)\pi^{-1},
\]
and for $\lambda\in \Lambda$ write
\[
a(\lambda)=\sum_{w\in W}(-1)^{w} w\cdot\lambda.
\]
Then the elements 
\[
a(\lambda)*\delta_{\Pi}
\]
span $\Hom_{\wa}(\Lambda,\Z(\epsilon))$.  Since $\pi_{1}G$ is
torsion-free, there is an exact sequence
\[
G'\to G\to J
\]
where $J$ is a torus, and $G'$ is simply connected.  The character
group of $J$ is the subgroup $\Lambda^{W}$ of Weyl-invariant elements
of $\Lambda$, and the weight lattice for $G'$ is the quotient
$\Lambda/\Lambda^{W}$.  Choose a Weyl chamber for $G$ and let $\rho\in
\Lambda\otimes\Q$ be $1/2$ the sum of the positive roots of $G$ (which
we will write as a product of square roots of elements in our Fourier
series notation).  Since $J$ is a torus, the image $\rho'$ of
$\rho$ in $\Lambda/\Lambda^{W}\otimes\Q$ is $1/2$ the sum of the
positive roots of $G'$, which since $G'$ is simply connected, lies in
$\Lambda/\Lambda^{W}$.  Let $\tilde \rho\in\Lambda$ be any element
congruent to $\rho'$ modulo $\Lambda^{W}$.  Claim: for any
$\lambda\in\Lambda$, the ratio
\[
\frac{a(\lambda)}{a(\tilde \rho)}
\]
is the character of a (virtual) representation.  The claim shows that
the class corresponding to $a(\tilde\rho)\cdot \delta_{\Pi}$ is an
$R(G)$-module generator of $K^{\tau+n}_{G}(G)$.  For the claim, first
note that the element $\mu=\rho/\tilde \rho$ is $W$-invariant (and is
in fact the square root of a character of $J$).  It follows from the
Weyl character formula that
\[
\frac{a((\lambda \tilde\rho^{-1})\cdot\rho)}{a(\rho)}
\]
is, up to sign, the character of an irreducible representation.   But
then 
\begin{align*}
\frac{a(\lambda\, \tilde\rho^{-1}\rho)}{a(\rho)} & = 
\frac{a(\lambda\, \mu)}{a(\rho)} 
= \mu\frac{a(\lambda)}{a(\rho)} \\
&= \frac{a(\lambda)}{a(\mu^{-1}\rho)}
= \frac{a(\lambda)}{a(\tilde\rho)}.
\end{align*}
\end{pf}

\begin{cor}
\label{thm:17} Let $U\in K^{\tau+n}_{G}(G)$ be the class corresponding
to $a(\tilde\rho)\cdot \delta_{\Pi}$.  The map ``multiplication by
$U$'' is an isomorphism
\[
R(G)/I^{\tau}\to K_{G}^{\tau+n}(G)
\]
of $R(G)$-modules.
\end{cor}

\begin{pf}
That the map factors through the quotient by $I^{\tau}$ is
Corollary~\ref{thm:48}, and that it is surjective is
Proposition~\ref{thm:49}.  The result now follows from the fact that
both sides are free of rank equal to the number of free $W$-orbits in
$A_{\epsilon,\text{reg}}$ (ie, the number of Verlinde conjugacy classes).
\end{pf}

As described at the end of \S\ref{sec:thom-isomorphism}, when $\tau$
is primitive the $R(G)$-module $K_{G}^{\tau+n}(G)$ acquires the
structure of an $R(G)$-algebra

\begin{cor}\label{thm:50}
When $\tau$ is primitive, there is a canonical algebra isomorphism
\[
K^{\tau+n}_{G}(G)\approx R(G)/I^{\tau}.
\]
\qed
\end{cor}

\begin{rem}
When $\tau$ is primitive, the pushforward map $K_{G}^{\tau}(e)\to
K_{G}^{{\tau+n}}$ is a ring homomorphism.  Being primitive, the
restriction of $\tau$ to $\{e \}\modmod G$ comes equipped with a
trivialization and so $K_{G}^{\tau}(e)\approx R(G)$.
\end{rem}

The isomorphisms of Corollaries~\ref{thm:17} and~\ref{thm:50} are
proved after tensoring with the complex numbers
in~\cite{math.AT/0206257}, where the distributions $\delta_{f}$ are
related to the Kac numerator at $q=1$.  We refer the reader to \S6 and
\S7 of~\cite{math.AT/0206257} for further discussion.

We conclude with a further computation which will be used in Part II.
We consider the situation of this section in which $G=T$ is a torus of
dimension $n$.   The group $K_{T}^{\tau}(\{e\})$ is the free abelian
group on $\LambdaTorsor$, and the pushforward map $K_{T}^{\tau}(\{e
\})\to K_{T}^{\tau+n}(T)$ is defined.
\begin{prop}\label{thm:10}
The pushforward map
\[
i_{!}:K_{T}^{\tau}(\{e
\})\to K_{T}^{\tau+n}(T)
\]
sends the class corresponding to $\lambda\in\LambdaTorsor$ to the class
in $K_{T}^{\tau+n}(T)\approx \Hom_{\Pi}(\LambdaTorsor ,\Z(\epsilon))$
corresponding to the distribution with Fourier expansion
\[
\lambda^{-1}\sum_{\pi\in\Pi}\epsilon(\pi)\pi^{-1}.
\]
\end{prop}

\begin{pf}
The pushforward map is the composition of the Thom isomorphism
\[
K^{\tau+0}_{T}(\{e \})\to K^{\tau+n}_{T}(\mathfrak{t},\mathfrak{t}\setminus\{0 \}) 
\approx K^{\tau+n}_{T}(T,T\setminus\{e \})
\]
with the restriction map
\begin{equation}\label{eq:19}
K^{\tau+n}_{T}(T,T\setminus\{e \})\to K^{\tau+n}_{T}(T).
\end{equation}
We wish to compute these maps using the spectral sequence for relative
twisted $K$-theory described at the end of
\S\ref{sec:fund-spectr-sequ}.  In order to do so, however, we need to
replace $T\setminus\{e \}$ and $\mathfrak{t}\setminus \{0 \}$ by the
smaller $T\setminus B_{e}$ and $\mathfrak{t}\setminus\{B_{0} \}$,
where $B_{e}$ and $B_{0}$ are small open balls containing $e$ and $0$
respectively.  This puts us in the situation described at the end of of
\S\ref{sec:fund-spectr-sequ}, where the sheaf
$\sk^{\tau+t}_{\text{rel}}$ works out to be the extension by zero of
the restriction of $\sk^{\tau+t}$ to $B_{e}$.

Applying the spectral sequence argument of this section to the pairs
$(\mathfrak{t},\mathfrak{t}\setminus B_{0})$ and $(T,T\setminus
B_{e})$ gives isomorphisms
\begin{equation}\label{eq:18}
\begin{aligned}
K^{\tau+n}_{T}(T,T\setminus B_{e}) &\approx 
\hom_{\Pi}(\LambdaTorsor, H^{n}_{c}(\mathfrak{t},\mathfrak{t}\setminus
B_{\Pi})\otimes\Z(\epsilon)) \\
K_{T}^{\tau+n}(\mathfrak{t},\mathfrak{t}\setminus B_{0}) &\approx 
\hom_{c}(\LambdaTorsor, H^{n}_{c}(\mathfrak{t},\mathfrak{t}\setminus
B_{0})\otimes \Z(\epsilon)), 
\end{aligned}
\end{equation}
where $B_{\Pi}\subset\mathfrak{t}$ is the inverse image of $B_{e}$
under the exponential map.  The same argument identifies the
restriction mapping~\eqref{eq:19} with the map induced by
\[
H^{n}_{c}(\mathfrak{t},\mathfrak{t}\setminus
B_{\Pi})\to H^{n}_{c}(\mathfrak{t}),
\]
and the isomorphism
\[
K^{\tau+n}_{T}(T,T\setminus B_{e})
\approx K^{\tau+n}_{T}(\mathfrak{t},\mathfrak{t}\setminus B_{0}) 
\]
with the map
\[
\Hom_{\Pi}(\LambdaTorsor, H^{n}_{c}(\mathfrak{t},\mathfrak{t}\setminus
B_{\Pi})\otimes\Z(\epsilon))\to 
\Hom_{c}(\LambdaTorsor, H^{n}_{c}(\mathfrak{t},\mathfrak{t}\setminus
B_{0})\otimes\Z(\epsilon))
\]
which first forgets the $\Pi$-action and then uses 
\[
H^{n}_{c}(\mathfrak{t},\mathfrak{t}\setminus B_{\Pi}) \to 
H^{n}_{c}(\bar{B}_{0}, B_{0}) \approx
H^{n}_{c}(\mathfrak{t},\mathfrak{t}\setminus B_{0}). 
\]
Finally, the isomorphism
\[
K^{\tau+n}_{T}(\mathfrak{t},\mathfrak{t}\setminus B_{0})\approx 
K^{\tau}_{T}(\text{pt})\otimes K^{n}(\mathfrak{t},\mathfrak{t}\setminus B_{0})
\]
shows that the Thom isomorphism is simply the tensor product of the
identity map with suspension isomorphism
$K^{0}(\text{pt})\to K^{n}(\mathfrak{t},\setminus B_{0})$ (which uses
the orientation of $\mathfrak{t}$).  In terms of~\eqref{eq:18} this
means that the Thom isomorphism 
\[
K^{\tau}_{T}(\text{pt}) \approx \Hom_{c}(\LambdaTorsor,\Z(\epsilon))\to 
K^{\tau}_{T}(\mathfrak{t},\mathfrak{t}\setminus B_{0}) 
\approx 
\hom_{c}(\LambdaTorsor, H^{n}_{c}(\mathfrak{t},\mathfrak{t}\setminus
B_{0})\otimes \Z(\epsilon)),
\]
is simply the map derived from the suspension isomorphism
\[
H^{0}(\text{pt})\approx 
H^{n}_{c}(\mathfrak{t},\mathfrak{t}\setminus B_{0}).
\]
The result follows easily from this.
\end{pf}

\subsection{The action of $\wa$ on $\mathfrak{t}$}
\label{sec:acti-waff-mathfr}

We summarize here some standard facts about about affine Weyl groups
and conjugacy classes in $G$.  Our basic references
are~\cite{bourbaki82:_elemen,helgason01:_differ_lie}.   Recall that
we have fixed a maximal torus $T$ of $G$.  We write $\Lambda$ for the
character group of $T$ and $R$ for the set of roots.  Following Bourbaki,
write $N(T,R)$ for the subgroup of $\mathfrak{t}$ consisting of
elements on which the roots vanish modulo $2\pi\,\Z$.  
There is a short exact sequence 
\[
N(T,R)\rightarrowtail \Pi\twoheadrightarrow \pi_{1}G.
\]

Let $\mathcal{H}$ be the set of hyperplanes forming the {\em diagram}
of $G$.  Thus
\[
\mathcal{H}=\left\{H_{k,\alpha}\mid k\in \Z ,\alpha\in R\right\},
\]
where 
\[
H_{k,\alpha}=\{x\in\mathfrak{t}\mid \alpha(x)=2\pi \,k
\}.
\]
The collection $\mathcal{H}$ is locally finite in the sense that
each $s\in\mathfrak{t}$ has a neighborhood meeting only finitely many
hyperplanes in $\mathcal{H}$.  The {\em affine Weyl group} is the
group $\waff$ be the group generated by reflections in the hyperplanes
$H_{k,\alpha}\in\mathcal{H}$.  It has the structure
\[
N(T,R)\rtimes W.
\]

\begin{prop}\label{thm:16}
Let $x\in\mathfrak{t}$.  The stabilizer of $x$ in $\waff$ is the
finite reflection group generated by reflections through the
hyperplanes $H_{k,\alpha}$ containing $x$.\qed
\end{prop}

Write $\wa=L\rtimes W$.  There is a short exact sequence 
\begin{equation}\label{eq:15}
\waff\rightarrowtail\wa\twoheadrightarrow\pi_{1}G.
\end{equation}

\begin{prop}\label{thm:12}
Let $x\in\mathfrak{t}$.  If $\pi_{1}G$ is torsion free, then the
stabilizer of $x$ in $\wa$ coincides with the stabilizer of $x$ in
$\waff$.  It is therefore the finite reflection group generated by
reflections through the hyperplanes $H_{k,\alpha}$ containing $x$.
\end{prop}

\begin{pf}
Write $W_{x}$ for the stabilizer of $x$ in $\wa$.  The image of
$W_{x}$ in $\mathfrak{t}\rtimes W$ is conjugate to a subgroup of $W$,
and so $W_{x}$ is finite.   By assumption $\pi_{1}G$ has no
non-trivial finite subgroups.  The exact
sequence~\eqref{eq:15} then shows that $W_{x}\subset\waff$.  The
result then follows from Proposition~\ref{thm:16}.
\end{pf}

Write $\R_{\ge0}=[0,\infty)$.
\begin{prop}\label{thm:13}
Suppose that $(W,V)$ is a finite reflection group.  The orbit space
$V/W$ is homeomorphic to $\R^{n_{1}}\times \R_{\ge0}^{n_{2}}$.  The
group $W$ is generated by $n_{2}$ reflections.  In particular, if $W$
is non-trivial, then $n_{2}\ne 0$.
\end{prop}

\begin{pf}
This follows immediately from the Theorems on pages 20 and 24
of~\cite{brown89:_build}.
\end{pf}

\appendix

\section{Groupoids}\label{sec:groupoids}

We remind the reader that we are assuming throughout this paper that,
unless otherwise specified, all spaces are locally contractible,
paracompact and completely regular.  These assumptions implies the
existence of partitions of unity~\cite{dold63:_partit} and locally
contractible slices through actions of compact Lie
groups~\cite{mostow57:_equiv_euclid,palais61:_lie}.

\subsection{Definition and First Properties}\label{sec:defin-first-prop}

A {\em groupoid} is a category in which all morphisms are isomorphism.
We will consider groupoids in the category of topological spaces.
Thus a groupoid $X=(X_{0},X_{1})$ consists of a space $X_{0}$ of
objects, a space $X_{1}$ of morphisms, and map ``identity map''
$X_{0}\to X_{1}$, a pair of maps ``domain'' and ``range'' $X_{1}\to
X_{0}$, an associative composition law $X_{1}\times_{X_{0}}X_{1}\to
X_{1}$, and an ``inverse'' map $X_{1}\to X_{1}$.  Write
$X_{n}=X_{1}\times_{X_{0}}\dots\times_{X_{0}}X_{1}$ for the space of
$n$-tuples of composeable maps.  Then the collection $\{X_{n} \}$ is a
simplicial space.   The $i^{\text{th}}$ face map
\[
d_{i}:X_{n}\to X_{n-1}
\]
is given by
\[
d_{i}(f_{1},\dots f_{n})=
\begin{cases}
(f_{2},\dots,f_{n}) \qquad& i=0 \\
(f_{1},\dots,f_{i}\circ f_{i+1},\dots,f_{n}) \qquad& 0<i<n \\
(f_{1},\dots,f_{n-1}) \qquad& i=n \\
\end{cases}
\]

Even though a groupoid is a special kind of simplicial space, we'll
refer to the simplicial space as the {\em nerve} of $X=(X_{0},X_{1})$
and write $X_{\bullet}$.   Finally, we let 
\[
\left|X \right|=\coprod_{n}X_{n}\times\Delta^{n}/\sim
\]
denote the geometric realization of $X_{\bullet}$.

\begin{eg}\label{eg:1}
(cf Segal~\cite{Segal:CSandSS}) Suppose that $G$ is a topological
group acting on a space $X$.  Then the pair $(G,X)$ forms a groupoid
with space of objects $X$ and in which a morphism from $x$ to $y$ is
an element of $g$ for which $g\cdot x=y$.  In this case $X_{0}=X$ and
$X_{1}= G\times X$.  The composition law is given by the
multiplication in $G$.  We will write $X\modmod{G}$ for this groupoid.
\end{eg}

\begin{eg}\label{eg:2} (Segal~\cite{Segal:CSandSS}) Suppose that $X$ is a space
and $\mathcal{U}=\{U_{i} \}$ is a covering of $X$.  The nerve of the
covering $\mathcal{U}$ is the nerve of a groupoid.  Indeed, let
$N_{\mathcal{U}}$ be the category whose objects are pairs $(U_{i},x)$
with $U_{i}\in\mathcal{U}$ and $x\in U_{i}$ and in which a morphism
from $(U_{i},x)$ to $(U_{j},y)$ is an element $w\in
U_{i}\times_{X}U_{j}$ whose projection to $U_{i}$ is $x$ and whose
projection to $U_{j}$ is $y$.  Then $N_{\mathcal{U}}$ is a groupoid.
If $U_{i}$ and $U_{j}$ are open {\em subsets} of $X$ then such a map
exists if and only if $x=y$, in which case it is unique.  Writing
$X_{0}=\coprod U_{i}$, then $X_{n}=
X_{0}\times_{X}\dots\times_{X}X_{0}$, and the nerve of this groupoid
is just the nerve of the covering $\mathcal{U}$.
\end{eg}

\begin{defin}\label{def:4}
A map of groupoids $F:X\to Y$ is an {\em
equivalence} if it is fully faithful and essentially surjective;
that is, if every object $y\in Y_{0}$ is isomorphic to one of the
form $Fx$, and if for every $a,b\in X_{0}$ the map 
\[
F:X(a,b)\to Y(Fa,Fb)
\]
is a homeomorphism.
\end{defin}

An equivalence of ordinary categories automatically admits an inverse
(up to natural isomorphism).  Examples~\ref{eg:10} and~\ref{eg:11}
below show that the same is not necessarily true of an equivalence of
topological categories, or groupoids.  Requiring the existence of a
globally defined inverse is, on the other hand too restrictive.  

\begin{defin}\label{def:15}
A {\em local equivalence} $X\to Y$ is an equivalence of groupoids
with the additional property that each $y\in Y_{0}$ has a
neighborhood $U$ admitting a lift in the diagram 
\[
\xymatrix{
&\tilde X_{0}\ar[r]\ar[d] & X_{0} \ar[d]\\
&Y_{1}\ar[r]^{\text{range}}\ar[d]^{\text{domain}} & Y_{0}\\
U \ar@{.>}[uur]\ar[r] & Y_{0} & 
}
\]
in which the square is Cartesian.
\end{defin}

\begin{rem}
The term {\em local equivalence} derives from thinking of a groupoid
$X$ as defining a sheaf $U\mapsto X(U)$ on the category of topological
spaces.  An equivalence $X\to Y$ has a globally defined inverse if and
only if for every space $U$ the map $X(U)\to Y(U)$ is an equivalence.
As one easily checks, a map $X\to Y$ is a local equivalence if and
only if it is an equivalence on stalks.  We will say that two
groupoids $X$ and $Y$ as being {\em weakly equivalent}, if there is a
diagram of local equivalences
\[
X\leftarrow Z \rightarrow Y.
\]
\end{rem}

\begin{eg}\label{eg:10}
If $\mathcal{U}$ is an open covering of a space $X$, then the map
\[
N_{\mathcal{U}}\to X
\]
is a local equivalence.  More generally, if $\mathcal{U}\to
\mathcal{V}$ is a map of coverings of $X$, then 
\[
N_{\mathcal{U}}\to N_{\mathcal{V}}
\]
is a local equivalence.
\end{eg}

\begin{eg}\label{eg:14}
Given groupoids $X$ and $A$, write $X(A)$ for the groupoid of maps $A
\to X$.  Then a map $X\to Y$ is a local equivalence if and only if
for spaces $S$, the map
\[
\varinjlim_{\mathcal U} X(U)\to
\varinjlim_{\mathcal U} Y(U)
\]
is an equivalence of groupoids, where $\mathcal U$ ranges over all
coverings of $S$.  Stated more succinctly, a map of groupoids is a
local equivalence if and only if the corresponding map of presheaves
of groupoids is a stalkwise equivalence.
\end{eg}

\begin{eg}\label{eg:11}
If $P\to X$ is a principal $G$-bundle over $X$, then 
\[
P\modmod{G} \to X
\]
is a local equivalence.  
\end{eg}

\begin{eg}\label{eg:15}
If $H\subset G$ is a subgroup, the map of groupoids
\[
\point\modmod{H} \to 
(G/H)\modmod{G}
\]
is a local equivalence.
\end{eg}

The fiber product of functors
\[
\xymatrix@C=7pt{
P\ar[dr]_{i} &  & Q\ar[dl]^{j} \\
& X &
}
\]
is the groupoid $P\times_{X}Q$ whose objects consist of 
\[
p\in P,\quad q\in Q,\quad x\in X
\]
and isomorphisms 
\[
ip\to x \leftarrow jq.
\]
The morphisms are the evident commutative diagrams.  To give a functor
$S\to P\times_{X}Q$ is to give a functors 
\[
S\xrightarrow{p} P,\quad S\xrightarrow{q} Q,\quad S\xrightarrow{x} X
\]
and natural isomorphisms 
\[
i\circ p\to x \leftarrow j\circ q.
\]

The groupoid $P\times_{X} Q$ is usually called {\em fiber product} of
$P$ and $Q$ over $X$, even though strictly speaking it is a kind of
homotopy fiber product and not the categorical fiber product.  We will
also say that the morphism $P\times_{X}Q\to Q$ is obtained from $P\to
X$ by change of base along $j:Q\to X$.  A natural transformation
$T:j_{1}\to j_{2}$ gives a natural isomorphism between the groupoids
obtained by change of base along $j_{1}$ and $j_{2}$.

One easily checks that the class of local equivalences is stable
under composition and change of base.  Consequently, if $P\to X$ and
$Q\to X$ are both local equivalences, so is $P\times_{X}Q\to X$.
Using Example~\ref{eg:14} one easily checks that if two of three maps
in a composition are local equivalences so is the third.

\begin{defin}\label{def:16} 
The $2$-category $\cov X$ is the category whose objects are local
equivalences $p:P\to X$ and in which a $1$-morphism from $p_{1}:P_{1}\to
X$ to $p_{2}:P_{2}\to X$ consists of a functor $F:P_{1}\to P_{2}$ and
a natural transformation $T:p_{1}\to p_{2}\circ F$ making
\[
\xymatrix@R=8pt@C=7pt{
P_{1}\ar[rr]^{F}\ar[ddr]_{p_{1}}&& P_{2}\ar[ddl]^{p_{2}} \\
\ar @{}@<1.5ex> [rr] | {\overset{T}{\implies}} && \\
& X &
}
\]
commute.  A $2$-morphism $(F_{1},T_{1})\to (F_{2},T_{2})$ is a natural
transformation $\eta:F_{1}\to F_{2}$  for which $T_{2}=p_{2}\eta\circ T_{1}$.
\end{defin}

We will denote by $\ocov X$ the $1$-category quotient of $\cov X$.
The objects of $\ocov X$ are those of $\cov X$, and $\ocov X(a,b)$
is the set of isomorphism classes in $\cov X(a,b)$.  We will see that
$\cov X$ and $\ocov X$ are not that different from each other.

\begin{lem}
\label{thm:39} For every $a,b\in \cov X$, the category $\cov X(a,b)$
is a codiscrete groupoid: there is a unique morphism between any two
objects.
\end{lem}

\begin{pf}
Write $a=(F_{1},T_{1})$ and $b=(F_{2},T_{2})$
\[
\xymatrix@R=8pt@C=8pt{
P_{1}\ar[rr]^{F_{1},F_{2}}\ar[ddr]_{p_{1}}&& P_{2}\ar[ddl]^{p_{2}} \\
\ar @{}[rr]^{{T_{1},T_{2}}} |{\implies} && \\
& X &
}.
\]
A morphism (natural transformation) $\eta\in\cov
X(a,b)$  associates to $x\in P_{1}$ a map 
\[
\eta_{x}:F_{1}x\to F_{2}x 
\]
whose image 
\[
p_{2}\eta_{x}:p_{2} F_{1}x\to p_{2}F_{2}x
\] 
is prescribed to fit into the diagram of isomorphisms
\[
\xymatrix@C=10pt{
 & p_{1}x \ar[dl]_{T_{1}x}\ar[dr]^{T_{2}x} & \\
p_{2}F_{1}x \ar[rr]_{p_{2}\eta_{x}} && p_{2}F_{2}x.
}
\]
The map $p_{2}\eta_{x}$ is therefore forced to be $(T_{2}x)\circ
(T_{1}x)^{-1}$, and so $\eta_{x}$ is uniquely determined since $P_{2}\to
X$ is an equivalence.
\end{pf}

\begin{cor}\label{thm:41}
The $1$-category quotient $\ocov X$ is a (co-)directed class.
\end{cor}

\begin{pf}
Suppose that $p_{i}:P_{i}\to X$, $i=1,2$ are two objects of $\cov X$.
The groupoid $P_{12}=P_{1}\times_{X}P_{2}$ comes equipped with maps
$P_{12}\to P_{1}$ and $P_{12}\to P_{2}$.   If $f,g:P\to Q$ are two
morphisms in $\cov X$ there is, by Lemma~\ref{thm:39}, a unique
$2$-morphism relating them, and so in fact $f=g$ in $\ocov X$
\end{pf}

\subsection{Further Properties of Groupoids}
\label{sec:furth-prop-group}

We now turn to several constructions which are invariants of local
equivalences.  A groupoid is a presentation of a {\em stack}, and the
invariants of local equivalences are in fact the invariants of the
underlying stack.  

\subsubsection{Point set topology of
groupoids}\label{sec:point-set-topology}

The {\em orbit space} or {\em coarse moduli space} of a groupoid $X$ is
the space of isomorphism classes of objects, topologized as a quotient
space of $X_{0}$.  We denote the coarse moduli space of $X$ by
$\cms{X}$.  At this level of generality, the space $\cms{X}$ can be
somewhat pathological, and without some further assumptions might not
be in our class of locally contractible, paracompact, and completely
regular spaces.  When $X=S\modmod{G}$, then $\cms{X}$ is the orbit
space $S/G$, and in that case we will revert to the more standard
notation $S/G$.  A local equivalence $Y\to X$ gives a homeomorphism
$\cms{Y}\to \cms{X}$.

For a subspace $S\subset \cms{X}$ we denote $X_{S}$ the full
sub-groupoid of $X$ consisting of objects in the isomorphism class of
$S$.  There is a one to one correspondence between full sub-groupoids
$A\subset X$ containing every object in their $X$-isomorphism class
and subspaces $\cms A$ of $\cms{X}$.  With this we transport many
notions from the point set topology of spaces to the context of
groupoids.  When $S$ is closed (resp.~open) we will say that $X_{S}$
is a closed (resp.~open) subgroupoid of $X$.  We can speak of the {\em
interior} and {\em closure} of a full subgroupoid.  By an {\em open
covering} of a groupoid, we mean an open covering $\{S_{\alpha} \}$ of
$\cms{X}$, in which case the collection $\{X_{S_{\alpha}} \}$ forms a
covering of $X$ by open sub-groupoids.

More generally, if $f:S\to \cms{X}$ is a map, we can form a groupoid $X_{S}$
with objects the pairs $(s,x)\in S\times X_{0}$ for which $x$ is in
the isomorphism class of $f(s)$.  A map $(s,x)\to (t,y)$ is just a map
from $x$ to $y$ in $X$.  Phrased differently, the groupoid $X_{S}$ is
the groupoid whose nerve fits into a pull-back square
\[
\begin{CD}
\left(X_{S} \right)_{\bullet} @>>> X_{\bullet}\\
@VVV @VVV \\
S @>>> \cms{X}.
\end{CD}
\]
We will say that $X_{S}$ is {\em defined by pullback} from the map
$S\to\cms{X}$.  With this we can transport many of the maneuvers of
homotopy theory to the context of groupoids.  For instance if $S=\cms
X\times I$, and $f$ is the projection, then $X_{S}$ is the groupoid
$X\times I$.  One can then form mapping cylinders and other similar
constructions.

For example, suppose that $\cms{X}$ is paracompact, and written as the
union of two sets $S_{0}$, $S_{0}$ whose interiors cover.  Write
$U_{i}=X_{S_{i}}$.  The $U_{i}$ are (full) subgroupoids whose interiors cover
$X$.  Let $N$ denote the groupoid constructed from the map
\[
\cyl \left(S_{0}\leftarrow S_{0}\cap S_{1}
\rightarrow S_{1} \right) \to \cms{X}.
\]
It is the double mapping cylinder of
\[
U_{0}\leftarrow U_{0}\cap U_{1}\rightarrow U_{1}.
\]
Following Segal~\cite{Segal:CSandSS}, a partition of unity
$\{\phi_{0},\phi_{1} \}$ subordinate to the covering $\{S_{0},S_{1}
\}$, defines a map
\[
\cms{X} \to \cyl \left(S_{0}\leftarrow S_{0}\cap S_{1}
\rightarrow S_{1} \right).
\]
The composite 
\[
\cms{X} \cyl \left(S_{0}\leftarrow S_{0}\cap S_{1}
\rightarrow S_{1} \right)\to \cms{X}
\]
is the identity, and so we get functors 
\[
X\to N \to X
\]
whose composite is the identity.   On the other hand, the composite 
\[
\cyl \left(S_{0}\leftarrow S_{0}\cap S_{1} \rightarrow S_{1} \right)
\cms{X}\to 
\cyl \left(S_{0}\leftarrow S_{0}\cap S_{1} \rightarrow S_{1} \right)
\]
is homotopic to the identity, by a homotopy (the obvious linear
homotopy) that preserves the map to $\cms{X}$ (it is a homotopy in the
category of spaces over $\cms{X}$).  This homotopy then defines by
pullback, a homotopy
\[
N\times\Delta^{1}\to N
\]
from the composite $N\to X\to N$ to the identity map of $N$, fixing
the map to $X$.  In this way $X$ becomes a strong deformation retract
of $N$, and $N$ is decomposed in a way especially well-suited for
constructing sequences of Mayer-Vietoris type.

A groupoid $X$ has {\em proper diagonal} if the map 
\begin{equation}\label{eq:7}
X_{1}\xrightarrow{(\text{domain},\text{range})}{X_{0}\times X_{0}} 
\end{equation}
is proper, and $\cms{X}$ is Hausdorff.  If~\eqref{eq:7} is proper and
$X_{0}$ is Hausdorff then $X$ is proper.  If $Y\to X$ is a local
equivalence, then $X$ has proper diagonal if and only if $Y$ has
proper diagonal.

\subsubsection{Local and global quotients}
\label{sec:local-glob-quot}

A groupoid which is related by a chain of local equivalences to one of
the form $S\modmod{G}$, obtained from a group $G$ acting on a space
$S$, is said to be a {\em global quotient}.  A {\em local quotient
groupoid} is a groupoid $X$ admitting a countable open cover
$\{U_{\alpha} \}$ with the property that each $X_{U_{\alpha}}$ is
weakly equivalent to a groupoid of the form $S\modmod{G}$ with $G$ a
compact Lie group, and $S$ a Hausdorff space.  If $Y\to X$ is a local
equivalence, then $Y$ is a local quotient groupoid if and only if $X$
is, so the property of being a local quotient is intrinsic to the
underlying stack.

If $X$ is a local quotient groupoid, then $\cms{X}$ is paracompact,
locally contractible and completely regular.  If $X$ is a local
quotient groupoid with the property that there is at most one map
between any two objects (ie $X_{1}\to X_{0}\times X_{0}$ is an
inclusion), the map $X\to \cms{X}$ is a local equivalence, and so
$X$ is just a space.

The following lemma is straightforward.

\begin{lem}\label{thm:45}
Any groupoid constructed by pullback from a local quotient groupoid is
a local quotient groupoid.  In particular, any (full) subgroupoid of a local
quotient groupoid is a local quotient groupoid, and the mapping
cylinder of a map $X_{S}\to X$ constructed by pullback along a map $S\to \cms{X}$
to the orbit space of a local quotient groupoid is a local quotient groupoid.
\end{lem}

\subsection{Fiber bundles over groupoids and descent}
\label{sec:fiber-bundles-over}

In this section we define the category of fiber bundles over a
groupoid, and show that a local equivalence gives an equivalence of
categories of fiber bundles (Proposition~\ref{thm:5}).  Thus the
category of fiber bundles over a groupoid is intrinsic to the
underlying stack.  

A fiber bundle over a groupoid $X=(X_{0},X_{1})$ consists of a fiber
bundle $P$ on $X_{0}$ together with identifications of certain
pullbacks to $X_{n}$ for various $n$.  We introduce some convenient
notation for describing these pulled back bundles.

Let's denote a typical point of $X_{n}$ by 
\[
x_{0}\xrightarrow{f_{1}}{}\dots
\xrightarrow{f_{n}}{}
x_{n}.
\]
Given a bundle $P\to X_{0}$ we'll write $P_{x_{i}}$ for the pullback of $P$
along the map 
\begin{align*}
X_{n} &\to X_{0}\\
(x_{0}\to\dots \to x_{n}) &\mapsto x_{i}. 
\end{align*}
Similarly, if $P\to X_{1}$ is given, we'll write $P_{f_{i}}$ 
for the pullback of $P$ along the map
\begin{align*}
X_{n} &\to X_{1}\\
(x_{0}\xrightarrow{f_{1}}{} \dots \xrightarrow{f_{n}} x_{n}) &\mapsto (x_{i-1}\xrightarrow{f_{i}}{}x_{i}),
\end{align*}
and $P_{f_{i}\circ f_{i+1}}$ for the pullback along
\begin{align*}
X_{n} &\to X_{1}\\
(x_{0}\xrightarrow{f_{1}}{}\dots \xrightarrow{f_{n}} x_{n}) &\mapsto
(x_{i-1}\xrightarrow{f_{i}\circ f_{i+1}}{}x_{i+1}),
\end{align*}
etc.   For small values of $n$ we'll use symbols like 
\begin{gather*}
(a\xrightarrow{f}{}b)\in X_{1} \\
(a\xrightarrow{f}{}b\xrightarrow{g}{}c)\in X_{2}.
\end{gather*}
to denote typical points.

\begin{defin}\label{def:5}
A {\em fiber bundle} on $X$ consists of a fiber bundle $P\to X_{0}$,
together with an bundle isomorphism
\begin{equation}\label{eq:5}
t_{f}:P_{a}\to P_{b}
\end{equation}
on $X_{1}$, for which $t_{\id}=\id$, and satisfying the cocycle
condition that
\[
\xymatrix{
P_{a} \ar[rr]^{t_{f}}\ar[dr]_{t_{g\circ f}} && P_{b}\ar[dl]^{t_{g}}\\
& P_{c} &
}
\]
commutes on $X_{2}$.
\end{defin}

This way of describing a fiber bundle is convenient when thinking of
$X$ as a category.  The association $a\to P_{a}$ is a functor from $X$
to spaces, that is continuous in an appropriate sense.  There is a
more succinct way of describing a fiber bundle on a groupoid.  Namely,
a fiber bundle on a groupoid $X=(X_{0},X_{1})$ is a groupoid
$P=(P_{0},P_{1})$ and a functor $P\to X$ making $P_{i}\to X_{i}$ into
fiber bundles, and all of the structure maps into maps of fiber
bundles (i.e., pullbacks squares).

A functor $F:Y\to X$ between groupoids defines, in the evident way, a
pullback functor $F^{\ast}$ from the category of fiber bundles over
$X$ to the category of fiber bundles over $Y$.  A natural
transformation $T:F\to G$ defines a natural transformation
$T^{\ast}:F^{\ast}\to G^{\ast}$.  

\begin{eg}\label{eg:9}
Let $\mathcal{U}=\{U_{i} \}$ be a covering of a space $X$.  To give a
fiber bundle over $N_{\mathcal{U}}$ is to give a fiber bundle $P_{i}$
on each $U_{i}$ and the clutching (descent) data needed to assemble
the $P_{i}$ into a fiber bundle over $X$.  Indeed, pullback along the
map $N_{\mathcal{U}}\to X$ gives an equivalence between the category
of fiber bundles over $X$ and the category of fiber bundles over
$N_{\mathcal{U}}$.
\end{eg}

The following generalization of Example~\ref{eg:9} will be referred to
as {\em descent} for fiber bundles over groupoids.

\begin{prop}\label{thm:5}
Suppose that $F:X\to Y$ is a local equivalence.  Then the pullback functor 
\[
F^{\ast}:\left\{\text{Fiber bundles on } Y \right\}
\to
\left\{\text{Fiber bundles on } X \right\}
\]
is an equivalence of categories.
\end{prop}

\begin{pf}
Suppose that $P$ is a fiber bundle over $Y$, which we think of as a
functor from $Y$ to the category of topological spaces.  Since $Y\to
X$ is an equivalence of categories, the functor $F^{\ast}$ has a left
adjoint $F_{\ast}$, given by
\[
F_{\ast}P(x) = \varinjlim_{Y/x} P,
\]
where $Y/x$ is the category of objects in $y\in Y$ equipped with a
morphism $Fy\to x$.  Since $Y\to X$ is an equivalence of groupoids,
there is a unique map between any two objects of $Y/x$, and so
$F_{\ast}P(x)$ is isomorphic to $P_{y}$ for any $y\in Y/x$.  For each
$x\in X$, choose a neighborhood $x\in U\subset X_{0}$, a map $t:U\to
Y_{0}$, and a family of morphisms $U\to X_{1}$ connecting $F\circ t$
to the inclusion $U\to X_{0}$.   We topologize
\[
\bigcup_{x\in X}F_{\ast}P_{x}
\]
by requiring that the canonical map
\[
t^{\ast}P \to F_{\ast}P\vert_{U}
\]
be a homeomorphism.  This gives $F_{\ast}P$ the structure of a fiber
bundle over $X_{0}$.  Naturality provides $F_{\ast}P$ with the
additional structure required to make it into a fiber bundle over $X$.
One easily checks that the pair $(F_{\ast}, F^{\ast})$ is an adjoint
equivalence of the category of fiber bundles over $X$ with the
category of fiber bundles over $Y$.
\end{pf}

For a fiber bundle $p:P\to X$ write $\Gamma(P)$ for the space of
sections
\[
\Gamma(P) = \Gamma(X;P) = \left\{ s:X\to P \mid p\circ s = \id_{X}\right\}, 
\]
topologized as a subspace of 
\[
X_{0}^{P_{0}}\times X_{1}^{P_{1}}.
\]
If $f:Y\to X$ is a
local equivalence, and $P\to X$ is a fiber bundle, then the evident
map
\[
\Gamma(X;P)\to \Gamma(Y;f^{\ast}P)
\]
is a homeomorphism.

If $P$ is a pointed fiber bundle, with $s:X\to P$ as a basepoint, and
$A\subset X$ is a (full) subgroupoid, write $\Gamma(X,A;P)$ for the space of
section $x$ of $P$ for which $x\vert_{A}=s$.  

Now suppose that $P\to Q$ is a map of fiber bundles over $X$, and
$\{U_{\alpha} \}$ is a covering of $X$ by open sub-groupoids.  Write
$P_{\alpha}\to U_{\alpha}$ for the restriction of $P$ to $U_{\alpha}$,
and $P_{\alpha_{1},\dots,\alpha_{n}}$ for the restriction of $P$ to
$U_{\alpha_{1}}\cap\dots\cap U_{\alpha_{n}}$, and similarly for $Q$.

\begin{prop}\label{thm:31}
If for each non-empty finite collection
$\{\alpha_{1}\,\dots\alpha_{n}\}$ the map
\[
\Gamma(P_{\alpha_{1},\dots,\alpha_{n}})
\to \Gamma(Q_{\alpha_{1},\dots,\alpha_{n}})
\]
is a weak homotopy equivalence, then so is 
\[
\Gamma(P)\to\Gamma(Q).
\]
\end{prop}

\begin{pf}
This is a straightforward application of the techniques of
Segal~\cite{Segal:CSandSS}.  Let's first consider the case in which
$X$ is covered by just two open sub-groupoids $U$ and $V$.  We form
the ``double mapping cylinder''
\[
C = U \amalg U\cap V\times [0,1] \amalg V/\sim,
\]
and consider the functor $g:C\to X$.  A choice of partition unity on
$\cms{X}$ subordinate to the covering $\{\cms U,\cms V \}$ gives a
section of $g$ making $\Gamma(X;P)\to \Gamma(X;Q)$ a retract of
$\Gamma(C;g^{\ast}P)\to \Gamma(C;g^{\ast}Q)$.  It therefore suffices,
in this case, to show that $\Gamma(C;g^{\ast}P)\to
\Gamma(C;g^{\ast}Q)$ is a weak equivalence.   But
$\Gamma(C;g^{\ast}P)$ fits into a homotopy pullback square
\[
\begin{CD}
\Gamma(C;g^{\ast}P) @>>> \Gamma(U;P) \\
@VVV @VVV \\
\Gamma(V;g^{\ast}P) @>>> \Gamma(U \cap V;P),
\end{CD}
\]
and similarly for $\Gamma(C;g^{\ast}Q)$ (to simplify the diagram we
have not distinguished in notation between $P$ and its restriction to
$U$, $V$, and $U\cap V$).  The result then follows from the long exact
(Mayer-Vietoris) sequence of homotopy groups.  An easy induction then
gives that the map on spaces of sections of $P\to Q$ restricted to any
finite union $U_{\alpha_{1}}\cup\dots\cup U_{\alpha_{n}}$ is a weak
equivalence.  For the case the collection $\{U_{\alpha} \}$ is
countable (to which we are reduced when $\cms{X}$ is second countable),
order the $U_{\alpha}$ and write 
\[
V_{n} = U_{1}\cup\dots\cup U_{n}.
\]
Form the ``infinite mapping cylinder''
\[
C = \coprod V_{i}\times [i,i+1]/\sim,
\]
and consider $g:C\to X$.  As before, a partition of unity on $\cms{X}$
subordinate to the covering $\cms V_{i}$ defines a section of
$\cms C\to \cms{X}$ and hence of $C\to X$, making $\Gamma(X;P)\to
\Gamma(X;Q)$ a retract of 
\begin{equation}\label{eq:29}
\Gamma(C;g^{\ast}P)\to
\Gamma(C;g^{\ast}Q).
\end{equation}
It therefore suffices to show that~\eqref{eq:29} is a weak
equivalence.  But~\eqref{eq:29} is the homotopy inverse limit of the
tower 
\begin{equation}\label{eq:30}
\Gamma(V_{n};P)\to \Gamma(V_{n};Q)
\end{equation}
and so its homotopy groups (or sets, in the case of $\pi_{0}$) are
related to those of~\eqref{eq:30} by a Milnor sequence, and the result follows.

Alternatively, following Segal~\cite{Segal:CSandSS}, one can avoid the
countability hypothesis and the induction by using for $C$ the nerve
of the covering $\{U_{\alpha} \}$ and the homotopy spectral sequences
of Bousfield-Kan and Bousfield~\cite{BousKan,bousfield89:_homot}.
\end{pf}

\subsection{Hilbert bundles}\label{sec:hilb-space-bundl}

A {\em Hilbert bundle} over a groupoid $X$ is a fiber bundle
whose fibers have the structure of a separable $\Z/2$-graded Hilbert
space.

\begin{rem}\label{rem:11}
There is a tricky issue in the point set topology here.  In defining
Hilbert bundles as special kinds of fiber bundles, we're implicitly
using the compact open topology on $U(H)$ and not the norm topology.
This causes trouble when we form the associated bundle of Fredholm
operators (\S\ref{sec:fredholm-operators}), since we cannot then use
the norm topology on the space of Fredholm operators.  This issue is
raised and resolved by Atiyah-Segal~\cite{atiyah:_twist_k}, and we are
following their discussion in this paper.
\end{rem}

\begin{defin}\label{def:10}
A Hilbert bundle $H$ is {\em universal} if for each Hilbert bundle $V$
there exists a unitary embedding $V\subset H$.  The bundle $H$ is said
to have the {\em absorption property} if for any $V$, there is a
unitary equivalence $H\oplus V\approx H$.
\end{defin}

\begin{lem}\label{thm:51}
A universal Hilbert bundle has the absorption property.
\end{lem}

\begin{pf}
First note that if $H$ is universal, then 
\[
H\otimes \ell^{2}\approx H\oplus H\oplus\cdots
\]
has the absorption property.  Indeed, given $V$ write $H=W\oplus V$,
and use the ``Eilenberg swindle''
\begin{multline*}
V\oplus H\oplus H\oplus\cdots \approx
V\oplus (W\oplus V)\oplus (W\oplus V)\oplus\cdots \\
\approx
(V\oplus W)\oplus (V\oplus W)\oplus (V\oplus W)\cdots  \approx
H\oplus H\oplus H\oplus\cdots.
\end{multline*}
We can then write $H\approx H\otimes\ell^{2}\oplus V\approx
H\otimes\ell^{2}$, to conclude that $H$ is absorbing.
\end{pf}

\begin{defin}\label{def:11}
A Hilbert bundle $H$ over $X$ is {\em locally universal}
if for every open sub-groupoid $X_{U}\subset X$ the restriction of $H$
to $X_{U}$ is universal.  
\end{defin}

\begin{lem}\label{thm:21}
If $H$ and $H'$ are universal Hilbert bundles on $X$, then there
is a unitary equivalence $H \approx H'$. \qed
\end{lem}

\begin{rem}
\label{rem:10} Since the category of Hilbert bundles on $X$ depends
only on $X$ up to local equivalence, if $f:Y\to X$ is a local
equivalence and $H$ is a (locally) universal Hilbert bundle on $X$,
then $f^{\ast}H$ is a (locally) universal Hilbert bundle on $Y$.
Similarly, if $H'$ is a (locally) universal Hilbert bundle on $Y$ there is
a (locally) universal Hilbert bundle $H$ on $X$, and a unitary equivalence
$f^{\ast}H\approx H'$.
\end{rem}

We now show that the existence of a locally universal Hilbert bundle
is a local issue.

\begin{lem}
\label{thm:22}
Suppose that $X$ is a groupoid, and that 
\[
\{U_{i}\mid i=1\dots\infty \}
\]
is a covering of $X$ by open sub-groupoids.  If $H$ is a Hilbert
bundle with the property that
\[
H_{i}=H\vert_{U_{i}}
\]
is universal, then $H\otimes \ell^{2}$ is universal.
\end{lem}

\begin{pf}
Let $V$ be a Hilbert-space bundle on $X$.  Choose a partition of unity
$\{\lambda_{i} \}$ on $\cms{X}$ subordinate to the open cover
$\cms U_{i}$.  For each $i$ choose an embedding
\[
r_{i}:V\vert_{U_{i}}\hookrightarrow H_{i}.
\]
The map 
\[
V\to H\oplus H\oplus\dots = H\otimes \ell^{2}
\]
with components $\lambda_{i}\, r_{i}$ is then an embedding of $V$ in
$H\otimes \ell^{2}$.
\end{pf}

\begin{cor}
\label{thm:24}
In the situation of Lemma~\ref{thm:22}, if $H\vert_{U_{i}}$ is locally
universal, then $H\otimes \ell^{2}$ is locally universal.\qed
\end{cor}

\begin{lem}
\label{thm:23} Suppose that $X$ is a groupoid, and that 
\[
\{U_{i}\mid i=1\dots\infty \}
\]
is a covering of $X$ by open sub-groupoids.  If $H_{i}$ is a locally
universal Hilbert bundle on $\{U_{i} \}$, then there exists a Hilbert
bundle $H$ on $X$ with $H\vert_{U_{i}}\approx H_{i}$.
\end{lem}

\begin{pf}
This is an easy induction, using Lemma~\ref{thm:21}.
\end{pf}

\begin{cor}
\label{thm:25}
Suppose that $X$ is a groupoid, and that 
\[
\{U_{i}\mid
i=1\dots\infty \}
\]
is a covering of $X$ by open sub-groupoids.  If
$H_{i}$ is a locally universal Hilbert bundle on $\{U_{i} \}$,
then there exists a locally universal Hilbert bundle $H$ on
$X$.  \qed
\end{cor}

\begin{lem}
\label{thm:26} Suppose that $X=S\modmod G$ is a global quotient of a
space $S$ by a compact Lie group $G$.  Then the equivariant Hilbert
bundle $S\times L^{2}(G)\otimes C_{1}\otimes \ell^{2}$ is a locally
universal Hilbert bundle on $X$.
\end{lem}

Here $C_{1}$ is the complex Clifford algebra on one (odd) generator.
It is there simply to make the odd component of our Hilbert bundle large enough.

\begin{pf}
Since the open (full) subgroupoids of $S\modmod G$ correspond to the
$G$-stable open subsets of $S$ it suffices to show that that
$L^{2}(G)\otimes C_{1}\otimes \ell^{2}$ is universal.  Let $V$ be any
Hilbert bundle on $S\modmod G$, ie an equivariant Hilbert bundle on
$S$.  By Kuiper's theorem, $V$ is trivial as a (non-equivariant)
Hilbert bundle on $S$.  Choose an orthonormal homogeneous basis
$\{e_{i} \}$, and let $e^{i}=\langle e_{i},\slot\rangle:V\to C_{1}$ be
the corresponding projection operator.  By the universal property of
$L^{2}(G)$, each $e^{i}$ lifts uniquely to an equivariant map
\[
V\to L^{2}(G)\otimes C_{1}.
\]
Taking the sum of these maps gives an embedding of $V$ in
$L^{2}(G)\otimes C_{1}\otimes\ell^{2}$.
\end{pf}

Combining Lemma~\ref{thm:26} with Lemma~\ref{thm:23} gives:

\begin{cor}
\label{thm:27}
If $X$ is a local quotient groupoid, then there exists a locally universal
Hilbert bundle on $X$. \qed
\end{cor}

\begin{cor}
\label{thm:47} Suppose that $X$ is a local quotient groupoid, $f:Y\to
X$ is a map constructed by pullback from $\cms{Y}\to \cms{X}$.  If $H$
is locally universal on $X$, then $f^{\ast}H$ is locally universal on
$Y$.
\end{cor}

\begin{pf}
This is an easy consequence of Lemma~\ref{thm:26} and
Corollary~\ref{thm:24}.
\end{pf}

Corollary~\ref{thm:47} is needed is the proof of excision in twisted
$K$-theory, and is the reason for our restriction to the class of
local quotient groupoids.

The following result is well-known, but we could not quite find a
reference.  Our proof is taken from~\cite[Theorem~1.5]{LMayS} which
gives the analogous result for equivariant embeddings of countably
infinite dimensional inner product spaces (and not Hilbert spaces).
Of course the result also follows from Kuiper's theorem, since the
space of embeddings is $U(H\otimes\ell^{2})/U(V^{\perp})$.  But the
contractibility of the space of embeddings is more elementary than the
contractibility of the unitary group, so it seemed better to have
proof that doesn't make use of Kuiper's theorem.

\begin{lem}
\label{thm:33} Suppose that $V$ and $H$ are Hilbert bundles over
a groupoid $X$, and that there is a unitary embedding $V\subset
H\otimes \ell^{2}$.  Then the space of embeddings $V\hookrightarrow
H\otimes\ell^{2}$ is contractible.
\end{lem}

\begin{pf}
Let $f:V\subset H\otimes\ell^{2}$ be a fixed embedding, and write
\begin{align*}
H\otimes\ell^{2} &= H\oplus H\oplus\cdots \\
f &= (f_{1},f_{2},\dots).
\end{align*}
The contraction is a concatenation of two homotopies.  The first takes
an embedding 
\[
g=(g_{1},g_{2},\dots)
\] 
to 
\[
(0,g_{1},0,g_{2},\dots)
\]
and then the second is
\[
\cos(\pi t/2)\cdot(0,g_{1},0,g_{2},\dots) + 
\sin(\pi t/2)\cdot(f_{1},0,f_{2},0\dots).
\]
It is easier to write down the reverse of the first homotopy.  It, in
turn, is the concatenation of an infinite sequence of $2$-dimensional
rotations
\begin{align*}
(0,g_{1},0,g_{2},0,g_{3},\dots) &\mapsto (g_{1},0,0,g_{2},0,g_{3},\dots)
\qquad 0\le t\le 1/2 \\
(g_{1},0,0,g_{2},0,g_{3},\dots) &\mapsto (g_{1},g_{2},0,0,0,g_{3},\dots)
\qquad 1/2\le t\le 3/4 \\
&\cdots .
\end{align*}
One must check that the limit as $t\mapsto 1$ is
$(g_{1},g_{2},\dots)$, and the that path is continuous in the
compact-open topology.  Both facts are easy and left to the reader.
\end{pf}

\begin{lem}
\label{thm:34} Suppose that $H$ is a locally universal Hilbert bundle
over a local quotient groupoid $X$.  Then the space of sections
$\Gamma(X;U(H))$ of the associated bundle of unitary groups is weakly
contractible.
\end{lem}

\begin{pf}
This follows easily from Kuiper's theorem (see Appendix~3
of~\cite{atiyah:_twist_k}), and Proposition~\ref{thm:31}.
\end{pf}

We conclude this section with a useful criterion for a local-quotient
stack to be equivalent to a global quotient by a compact
Lie group.  

\begin{prop}
The (locally) universal Hilbert bundle over a compact, local-quotient
groupoid, splits into a sum of finite-dimensional bundles iff the
groupoid is equivalent (in the sense of local equivalence) to one of
the form $X\modmod G$, with $X$ compact, and $G$ a compact group.
\end{prop}

\begin{rem}\label{rem:4}
\begin{trivlist}\itemsep0ex
\item (i) This implies right away that the extensions of groupoids
corresponding to twistings whose invariant in $H^3$ has infinite order
are \textit{not} quotient stacks: indeed, any $1$-eigenbundle for the
central $\T$ is a projective bundle representative for the twisting,
and hence must be infinite-dimensional.
\item (ii) There are simple obstructions to a groupoid being related
by a chain of local equivalences to a global quotient by a compact
group; for instance, such quotients admit continuous choices of
Ad-invariant metrics on the Lie algebra stabilizers which are integral
on the co-weight lattices.  The stack obtained by gluing the
boundaries of $B(\T\times\T)\times [0,1]$ via the shearing
automorphism of $\T\times\T$ does not carry such metrics. The same is
true for the quotient stack $\mathcal A\modmod \T\ltimes LT$, where
$T$ is a torus, and $\mathcal A$ is the space of connections on the
trivial $T$-bundle over the circle.  In this case, the stack is
fibered over $T$ in $B(\T\times T) $-stacks with the tautological
shearing holonomies.  Hence, the larger stacks $\mathcal A\modmod
\T\ltimes LG$ where $G$ is a compact Lie group and $\mathcal A$ is the
space of connections on the trivial $G$-bundle over the circle are not
global quotients either.
\item (iii) The result is curiously similar to Totaro's characterisation
of smooth quotient stacks as the Artin stacks where coherent sheaves
admit resolutions by vector bundles~\cite{totaro04}.

\end{trivlist}
\end{rem}
\begin{proof}
The `if' part follows from our construction of the universal Hilbert
bundle.  For the `only if part,' first note that a local quotient
groupoid $X$ is weakly equivalent to a groupoid of the form $S\modmod
G$, if and only if there is a principal $G$-bundle $P\to X$ with the
property that there is at most one map between any two objects in $P$.
In that case $P$ is equivalent to $\cms{P}$
(\S\ref{sec:local-glob-quot}), and $X$ is weakly equivalent to
$\cms{P}\modmod G$.  This latter condition holds if and only if for
each $x\in X_{0}$ the map $\Aut(x)\to G$ associated to $P$ is a
monomorphism.  Suppose that $H$ is the (locally) universal Hilbert
bundle on $X$, and that we can find an orthogonal decomposition
$H=\oplus H_{\alpha}$ with each $H_{\alpha}$ of dimension
$n_{\alpha}<\infty$.  Take $P$ to be the product of the frame bundles
of the $H_{\alpha}$, and $G$ to be the product of the unitary groups
$U(H_{\alpha})\approx U(n_{\alpha})$.  To check that $\Aut(x)\to G$ is
a monomorphism in this case it suffices to check locally near $x$.
The assertion is thus reduced to the case of a global quotient by a
compact group, where it follows from our explicit construction.
\end{proof}

\subsection{Fredholm operators and
$K$-theory}\label{sec:fredholm-operators}

We will build our model of twisted $K$-theory using the ``skew-adjoint
Fredholm'' model of Atiyah-Singer~\cite{atiyah69:_index_fredh}.  In
this section we recall this theory, and the modifications 
described in Atiyah-Segal~\cite{atiyah:_twist_k}

Let $H$ be a $\Z/2$-graded Hilbert bundle over a groupoid $X$.  We
wish to associate to $H$ a bundle of Fredholm operators over $X$.  As
mentioned in Remark~\ref{rem:11}, we cannot just use the norm topology
on the space of Fredholm operators here.  We have used the
compact-open topology on $U(H)$, and the (conjugation) action of
$U(H)$ in the compact-open topology on Fredholm in the norm topology
is not continuous.  Following Atiyah-Segal~\cite[Definition
3.2]{atiyah:_twist_k}, we make the following definition.

\begin{defin}(\cite{atiyah:_twist_k}) The space $\fred{0}(H)$ is the
space of odd skew-adjoint Fredholm operators $A$, for which $A^{2}+1$
is compact, topologized as a subspace of
$\mathcal{B}(H)\times\mathcal{K}(H)$, with $\mathcal{B}(H)$ given the
compact-open topology and $\mathcal{K}(H)$ the norm topology.
\end{defin}

Let $C_{n}=T\{\C^{n}\}/(z^{2}+q(z)=0)$ denote the complex Clifford
algebra associated to the quadratic form $q(z)=\sum z_{i}^{2}$.  We
write $\epsilon_{i}$ for the $i^{\text{th}}$ standard basis element of
$\C^{n}$, regarded as an element of $C_{n}$.  Following
Atiyah-Singer~\cite{atiyah69:_index_fredh}, for an operator
$A\in\fred(C_{n}\otimes H)$, with $n$ odd, let
\[
w(A) = \begin{cases}
\epsilon_{1}\dots\epsilon_{n}\, A &\qquad n \equiv -1 \mod 4 \\
i^{-1}\epsilon_{1}\dots\epsilon_{n}\, A &\qquad n \equiv 1 \mod 4.
\end{cases}
\]
The operator $A$ is then even and self-adjoint.

\begin{defin}(\cite{atiyah:_twist_k})\label{def:13}
The space $\fred{n}(H)$ is the subspace of $\fred{0}(C_{n}\otimes H)$
consisting of odd operators $A$ which commute (in the graded sense)
with the action of $C_{n}$, and for which the essential spectrum of
$w(A)$, in case $n$ is odd, contains both positive and negative
eigenvalues.
\end{defin}

Atiyah and Segal~\cite{atiyah:_twist_k} show that the ``identity'' map
from $\fred{n}(\ell^{2})$ in the norm topology to $\fred{n}(\ell^{2})$
in the above topology is a weak homotopy equivalence.  It then follows
from Atiyah-Singer~\cite[Theorem~B(k)]{atiyah69:_index_fredh} that the
map
\begin{align*}
\fred{n}(\ell^{2}) &\to 
\Omega'\fred{n-1}(C_{1}\otimes \ell^{2}) \\
A &\mapsto \epsilon_{k}\cos(\pi t)+A\,\sin(\pi t)
\end{align*}
is weak homotopy equivalence, where we are making the evident
identification $C_{n}\approx C_{n-1}\otimes C_{1}$, and $\Omega'$
denotes the space of paths from $\epsilon_{k}$ to $-\epsilon_{k}$.
Combining these leads to the following simple consequence.

\begin{prop}
\label{thm:32} If $X$ is a local quotient groupoid, and $H$ a
$\Z/2$-graded, locally universal Hilbert bundle over $X$, the
map
\[
\Gamma(X;\fred{n+1}(H))\to \Omega'\Gamma(X;\fred{n}(H))
\]
is a weak homotopy equivalence.
\end{prop}

\begin{pf}
By Proposition~\ref{thm:31}, the question is local in $X$, so we may
assume $X=S\modmod G$, with $G$ a compact Lie group.  By our
assumption on the existence of locally contractible slices, we may
reduce to the case in which $S$ is equivariantly contractible to a
fixed point $s\in S$.  Finally, since the question is homotopy
invariant in $X$, we reduce to the case $S=\point$.  We are therefore
reduced to showing that if $H$ is a universal $G$-Hilbert space, then
the map of $G$-fixed points
\begin{equation}\label{eq:31}
\fred{n}(H)^{G}\to
\Omega'\fred{n+1}(H)^{G}
\end{equation}
is a weak equivalence.  For each irreducible representation $V$ of
$G$, let $H_{V}$ denote the $V$-isotypical component of $H$.
Then~\eqref{eq:31} is the product over the irreducible representations
$V$ of $G$, of 
\begin{equation}\label{eq:32}
\fred{n}(H_{V})^{G}\to
\Omega'\fred{n+1}(H_{V})^{G}
\end{equation}
Since $H$ is universal, the Hilbert space $H_{V}$ is isomorphic to
$V\otimes\ell^{2}$, and the map
\begin{align*}
\fred{n}(\ell^{2}) &\to 
\fred{n}(V\otimes\ell^{2})^{G} \\
T &\mapsto \id\otimes T
\end{align*}
is a homeomorphism.  The Proposition is thus reduced to the result of
Atiyah-Singer quoted above.
\end{pf}

We now assemble the spaces $\Gamma(X;\fred{n}(H))$ into a spectrum in
the sense of algebraic topology.  To do this requires specifying
basepoints in $\fred{n}(H)$.  Since our operators are odd, we can't
take the identity map as a basepoint and a different choice must be
made.  There are some technical difficulties that arise in trying to
specify consistent choices and we have just chosen to be unspecific on
this point.  The difficulties don't amount to a serious problem since
any invertible operator can be taken as a basepoint, and the space of
invertible operators is contractible.  The reader is referred
to~\cite{joachim01:_k_o} for further discussion.

We will use the symbol $\epsilon$ to refer to a chosen basepoint in
$\fred{n}(H)$, as well as to the constant section with value $\epsilon$
in $\Gamma(X;\fred{n}(H))$.

Proposition~\ref{thm:32} gives a homotopy equivalence
\begin{equation}\label{eq:33}
\Gamma(X;\fred{n+1}(W)) \to\Omega
\Gamma(X;\fred{n}(W))
\end{equation}
As described in \cite{atiyah69:_index_fredh}, the fact
that $C_{2}$ is a matrix algebra gives a homeomorphism
\begin{equation}\label{eq:34}
\Gamma(X\fred{m}(W))\approx\Gamma(X\fred{m+2}(W)).
\end{equation}
We the spectrum $\ul K(X)$ by taking 
\[
\ul K(X)_{n} = 
\begin{cases}
\Gamma(X,\fred{0}(W)) &\qquad \text{$n$ even} \\
\Gamma(X,\fred{1}(W)) &\qquad \text{$n$ odd} \\
\end{cases}
\]
with structure map $\ul K(X)_{n} \to \Omega \ul K(X)_{n+1}$ to be the
map~\eqref{eq:33} when $n$ is odd, and the composite of~\eqref{eq:34}
and~\eqref{eq:33} when $n$ is even.  The group $K^{n}(X)$ is then
defined by
\[
K^{n}(X)=\pi_{0}\ul K(X)_{n}\approx \pi_{k}\ul K(X)_{n+k}.
\]
Because $H$ is locally universal, when $X=S\modmod G$, we have 
\[
K^{n}(X)\approx [S,\fred{n}(L^{2}(G)\otimes\ell^{2})]_{G}.
\]
Since, as remarked in~\S\ref{sec:proof-prop-refthm:18},
$\fred{n}(L^{2}(G)\otimes\ell^{2})$ is a classifying space for
equivariant $K$-theory, this latter group can be identified
with
\[
K^{n}(G)(S).
\]


\begin{thebibliography}{10}

\bibitem{atiyah69:_index_fredh}
M.~F. Atiyah and I.~M. Singer, \emph{Index theory for skew-adjoint {F}redholm
  operators}, Inst. Hautes \'Etudes Sci. Publ. Math. (1969), no.~37, 5--26.
  \MR{44 \#2257}

\bibitem{math/0510674v1}
Michael Atiyah and Graeme Segal, \emph{{Twisted K-theory and cohomology}}.

\bibitem{atiyah:_twist_k}
\bysame, \emph{Twisted {$K$}-theory}, Ukr. Mat. Visn. \textbf{1} (2004), no.~3,
  287--330. \MR{MR2172633}

\bibitem{bezrukavnikov05:_equiv_k_grass_toda}
Roman Bezrukavnikov, Michael Finkelberg, and Ivan Mirkovi{\'c},
  \emph{Equivariant homology and {$K$}-theory of affine {G}rassmannians and
  {T}oda lattices}, Compos. Math. \textbf{141} (2005), no.~3, 746--768.
  \MR{MR2135527 (2006e:19005)}

\bibitem{bott58:_lie}
Raoul Bott, \emph{The space of loops on a {L}ie group}, Michigan Math. J.
  \textbf{5} (1958), 35--61. \MR{MR0102803 (21 \#1589)}

\bibitem{bott58:_applic_morse}
Raoul Bott and Hans Samelson, \emph{Applications of the theory of {M}orse to
  symmetric spaces}, Amer. J. Math. \textbf{80} (1958), 964--1029.
  \MR{MR0105694 (21 \#4430)}

\bibitem{bourbaki82:_elemen}
Nicolas Bourbaki, \emph{\'{E}l\'ements de math\'ematique: groupes et alg\`ebres
  de {L}ie}, Masson, Paris, 1982, Chapitre 9. Groupes de Lie r\'eels compacts.
  [Chapter 9. Compact real Lie groups]. \MR{84i:22001}

\bibitem{bousfield89:_homot}
A.~K. Bousfield, \emph{Homotopy spectral sequences and obstructions}, Israel J.
  Math. \textbf{66} (1989), no.~1-3, 54--104. \MR{MR1017155 (91a:55027)}

\bibitem{BousKan}
A.~K. Bousfield and D.~M. Kan, \emph{Homotopy limits, completions and
  localizations}, Lecture Notes in Mathematics, no. 304, Springer--Verlag,
  Berlin, 1972.

\bibitem{bouwknegt02:_twist_k_k}
Peter Bouwknegt, Alan~L. Carey, Varghese Mathai, Michael~K. Murray, and Danny
  Stevenson, \emph{Twisted {$K$}-theory and {$K$}-theory of bundle gerbes},
  Comm. Math. Phys. \textbf{228} (2002), no.~1, 17--45. \MR{1 911 247}

\bibitem{brown89:_build}
Kenneth~S. Brown, \emph{Buildings}, Springer-Verlag, New York, 1989.
  \MR{90e:20001}

\bibitem{bunke05:_twist_k_tqft}
U.~Bunke and I.~Schr{\"o}der, \emph{Twisted {$K$}-theory and {TQFT}},
  Mathematisches Institut, Georg-August-Universit\"at G\"ottingen: Seminars
  Winter Term 2004/2005, Universit\"atsdrucke G\"ottingen, G\"ottingen, 2005,
  pp.~33--80. \MR{MR2206878 (2007c:19008)}

\bibitem{dold63:_partit}
Albrecht Dold, \emph{Partitions of unity in the theory of fibrations}, Ann. of
  Math. (2) \textbf{78} (1963), 223--255. \MR{27 \#5264}

\bibitem{donovan70:_graded_brauer_k}
P.~Donovan and M.~Karoubi, \emph{Graded {B}rauer groups and {$K$}-theory with
  local coefficients}, Inst. Hautes \'Etudes Sci. Publ. Math. (1970), no.~38,
  5--25. \MR{43 \#8075}

\bibitem{douglas06:_k_lie}
Christopher~L. Douglas, \emph{On the twisted {$K$}-homology of simple {L}ie
  groups}, Topology \textbf{45} (2006), no.~6, 955--988. \MR{MR2263220}

\bibitem{duistermaat00:_lie}
J.~J. Duistermaat and J.~A.~C. Kolk, \emph{Lie groups}, Universitext,
  Springer-Verlag, Berlin, 2000. \MR{MR1738431 (2001j:22008)}

\bibitem{freed94:_higher}
Daniel~S. Freed, \emph{Higher algebraic structures and quantization}, Comm.
  Math. Phys. \textbf{159} (1994), no.~2, 343--398. \MR{MR1256993 (95c:58034)}

\bibitem{freed:_consis}
Daniel~S. Freed, Michael~J. Hopkins, and Constantin Teleman, \emph{{Consistent
  orientation of moduli spaces}}.

\bibitem{freed:_loop_group_twist_k_theor_ii}
\bysame, \emph{{Loop Groups and Twisted K-Theory II}}, arXiv:math.AT/0511232.

\bibitem{freed:_loop_group_twist_k_iii}
\bysame, \emph{{Loop Groups and Twisted K-theory III}}, arXiv:math.AT/0312155.

\bibitem{math.AT/0206257}
\bysame, \emph{{Twisted equivariant K-theory with complex coefficients}},
  arXiv:math.AT/0206257, accepted for publication by {T}opology.

\bibitem{helgason01:_differ_lie}
Sigurdur Helgason, \emph{Differential geometry, {L}ie groups, and symmetric
  spaces}, Graduate Studies in Mathematics, vol.~34, American Mathematical
  Society, Providence, RI, 2001, Corrected reprint of the 1978 original.
  \MR{2002b:53081}

\bibitem{joachim01:_k_o}
Michael Joachim, \emph{A symmetric ring spectrum representing {$K{\rm
  O}$}-theory}, Topology \textbf{40} (2001), no.~2, 299--308. \MR{MR1808222
  (2001k:55011)}

\bibitem{kitchloo:_domin_k_kac_moody}
Nitu Kitchloo, \emph{{D}ominant {$K$}-theory and highest weight representations
  of {K}ac-{M}oody groups}, preprint.

\bibitem{LMayS}
L.~G. Lewis, J.~P. May, and M.~Steinberger, \emph{Equivariant stable homotopy
  theory}, Lecture Notes in Mathematics, vol. 1213, Springer--Verlag, New York,
  1986.

\bibitem{minasian97:_k_ramon_ramon}
Ruben Minasian and Gregory Moore, \emph{{$K$}-theory and {R}amond-{R}amond
  charge}, J. High Energy Phys. (1997), no.~11, Paper 2, 7 pp. (electronic).
  \MR{2000a:81190}

\bibitem{mostow57:_equiv_euclid}
G.~D. Mostow, \emph{Equivariant embeddings in {E}uclidean space}, Ann. of Math.
  (2) \textbf{65} (1957), 432--446. \MR{MR0087037 (19,291c)}

\bibitem{palais61:_lie}
Richard~S. Palais, \emph{On the existence of slices for actions of non-compact
  {L}ie groups}, Ann. of Math. (2) \textbf{73} (1961), 295--323. \MR{MR0126506
  (23 \#A3802)}

\bibitem{rosenberg89:_contin}
Jonathan Rosenberg, \emph{Continuous-trace algebras from the bundle theoretic
  point of view}, J. Austral. Math. Soc. Ser. A \textbf{47} (1989), no.~3,
  368--381. \MR{MR1018964 (91d:46090)}

\bibitem{Segal:CSandSS}
G.~Segal, \emph{Classifying spaces and spectral sequences}, Inst. Hautes
  \'Etudes Sci. Publ. Math. \textbf{34} (1968), 105--112.

\bibitem{stolz04:_what}
Stephan Stolz and Peter Teichner, \emph{What is an elliptic object?}, Topology,
  geometry and quantum field theory, London Math. Soc. Lecture Note Ser., vol.
  308, Cambridge Univ. Press, Cambridge, 2004, pp.~247--343. \MR{MR2079378
  (2005m:58048)}

\bibitem{totaro04}
Burt Totaro, \emph{The resolution property for schemes and stacks}, J. Reine
  Angew. Math. \textbf{577} (2004), 1--22. \MR{MR2108211 (2005j:14002)}

\bibitem{tu04:_twist_k}
Jean-Louis Tu, Ping Xu, and Camille Laurent-Gengoux, \emph{Twisted {$K$}-theory
  of differentiable stacks}, Ann. Sci. \'Ecole Norm. Sup. (4) \textbf{37}
  (2004), no.~6, 841--910. \MR{MR2119241 (2005k:58037)}

\bibitem{verlinde88:_fusion_d}
Erik Verlinde, \emph{Fusion rules and modular transformations in {$2$}{D}
  conformal field theory}, Nuclear Phys. B \textbf{300} (1988), no.~3,
  360--376. \MR{89h:81238}

\bibitem{witten98:_d_k}
Edward Witten, \emph{D-branes and {$K$}-theory}, J. High Energy Phys. (1998),
  no.~12, Paper 19, 41 pp.\ (electronic). \MR{2000e:81151}

\end{thebibliography}

\providecommand{\bysame}{\leavevmode\hbox to3em{\hrulefill}\thinspace}
\providecommand{\MR}{\relax\ifhmode\unskip\space\fi MR }
\providecommand{\MRhref}[2]{%
  \href{http://www.ams.org/mathscinet-getitem?mr=#1}{#2}
}
\providecommand{\href}[2]{#2}

\end{document}